\documentclass[final,leqno,onefignum,onetabnum]{siamltex1213}
\usepackage{amsfonts}
\usepackage{amsmath,booktabs,ctable,threeparttable}
\usepackage{amssymb,amsfonts,boxedminipage}
\newtheorem{remark}{Remark}[section]
\newtheorem{prop}{Proposition}[section]

\title{Approximation Accuracy of the Krylov Subspaces for Linear Discrete
Ill-Posed Problems\thanks{This work was supported in part by
the National Science Foundation of China (No. 11771249).}}

\author{Zhongxiao Jia\thanks{Department of Mathematical Sciences, Tsinghua
University, 100084 Beijing, China. (\email{jiazx@tsinghua.edu.cn})}}

\begin{document}
\maketitle
\slugger{sirev}{xxxx}{xx}{x}{x--x}

\begin{abstract}
For the large-scale linear discrete ill-posed problem $\min\|Ax-b\|$ or $Ax=b$
with $b$ contaminated by Gaussian white noise, the Lanczos bidiagonalization
based Krylov solver LSQR and its mathematically equivalent CGLS, the
Conjugate Gradient (CG) method implicitly applied to $A^TAx=A^Tb$,
are most commonly used, and CGME, the CG method applied to
$\min\|AA^Ty-b\|$ or $AA^Ty=b$ with $x=A^Ty$, and
LSMR, which is equivalent to the minimal residual (MINRES) method
applied to $A^TAx=A^Tb$, have also been choices.
These methods exhibit typical semi-convergence feature,
and the iteration number $k$ plays the role of the regularization parameter.
However, there has been no definitive answer to the long-standing fundamental
question:
{\em Can LSQR and CGLS find 2-norm filtering best possible regularized solutions}?
The same question is for CGME and LSMR too.
At iteration $k$, LSQR, CGME and LSMR compute {\em different} iterates from the
{\em same} $k$ dimensional Krylov subspace.
A first and fundamental step towards to answering the above question
is to {\em accurately} estimate the accuracy of the underlying
$k$ dimensional Krylov subspace approximating the $k$ dimensional dominant right
singular subspace of $A$. Assuming that the singular values of $A$ are simple,
we present a general $\sin\Theta$ theorem for the 2-norm distances between these
two subspaces and derive accurate estimates on them for severely, moderately
and mildly ill-posed problems. We also
establish some relationships between the smallest Ritz values and
these distances. Numerical experiments justify the sharpness of our
results.
\end{abstract}


\pagestyle{myheadings}
\thispagestyle{plain}
\markboth{ZHONGXIAO JIA}{APPROXIMATION ACCURACY OF THE KRYLOV SUBSPACES}

\maketitle

\section{Introduction and Preliminaries}\label{intro}

Consider the linear discrete ill-posed problem
\begin{equation}
  \min\limits_{x\in \mathbb{R}^{n}}\|Ax-b\| \mbox{\,\ or \ $Ax=b$,}
  \ \ \ A\in \mathbb{R}^{m\times n}, \label{eq1}
  \ b\in \mathbb{R}^{m},
\end{equation}
where the norm $\|\cdot\|$ is the 2-norm of a vector or matrix, and
$A$ is extremely ill conditioned with its singular values decaying
to zero without a noticeable gap. Since the results in this paper
hold for both the overdetermined ($m\geq n$) and underdetermined ($m\leq n$) cases,
we assume that $m\geq n$ for brevity.
\eqref{eq1} typically arises
from the discretization of the first kind Fredholm integral equation
\begin{equation}\label{eq2}
Kx=(Kx)(t)=\int_{\Omega} k(s,t)x(t)dt=g(s)=g,\ s\in \Omega
\subset\mathbb{R}^q,
\end{equation}
where the kernel $k(s,t)\in L^2({\Omega\times\Omega})$ and
$g(s)$ are known functions, while $x(t)$ is the
unknown function to be sought. If $k(s,t)$ is non-degenerate
and $g(s)$ satisfies the Picard condition, there exists the unique square
integrable solution
$x(t)$; see \cite{engl00,hansen98,hansen10,kirsch,mueller}. Here for brevity
we assume that $s$ and $t$ belong to the same set $\Omega\subset
\mathbb{R}^q$ with $q\geq 1$.
Applications include image deblurring, signal processing, geophysics,
computerized tomography, heat propagation, biomedical and optical imaging,
groundwater modeling, and many others; see, e.g.,
\cite{aster,engl93,engl00,hansen10,ito15,kaipio,kern,kirsch,mueller,natterer,vogel02}.
The theory and numerical treatments of \eqref{eq2} can be found in \cite{kirsch,kythe}.
The right-hand side $b=b_{true}+e$ is noisy and assumed to be
contaminated by a Gaussian white noise $e$, caused by measurement, modeling
or discretization errors, where $b_{true}$
is noise-free and $\|e\|<\|b_{true}\|$.
Because of the presence of noise $e$ and the extreme
ill-conditioning of $A$, the naive
solution $x_{naive}=A^{\dagger}b$ of \eqref{eq1} bears no relation to
the true solution $x_{true}=A^{\dagger}b_{true}$, where
$\dagger$ denotes the Moore-Penrose inverse of a matrix.
Therefore, one has to use regularization to extract a
best possible approximation to $x_{true}$.

We assume that $b_{true}$ satisfies the discrete
Picard condition $\|A^{\dagger}b_{true}\|\leq C$ with some constant $C$
for $\|A^{\dagger}\|$  arbitrarily
large \cite{aster,gazzola15,hansen90,hansen90b,hansen98,hansen10,kern}.
It is a discrete analog of the Picard condition in the Hilbert space setting;
see, e.g., \cite{hansen90}, \cite[p.9]{hansen98},
\cite[p.12]{hansen10} and \cite[p.63]{kern}. Without loss of generality,
assume that $Ax_{true}=b_{true}$. Then for a Gaussian white noise $e$,
the two dominating regularization
approaches are to solve the following two equivalent problems:
\begin{equation}\label{posed}
\min\limits_{x\in \mathbb{R}^{n}}\|Lx\| \ \ \mbox{subject to}\ \
\|Ax-b\|\leq \tau\|e\|
\end{equation}
with $\tau>1$ slightly and general-form Tikhonov regularization
\begin{equation}\label{tikhonov}
  \min\limits_{x\in \mathbb{R}^{n}}\{\|Ax-b\|^2+\lambda^2\|Lx\|^2\}
\end{equation}
with $\lambda>0$ the regularization parameter
\cite{hansen98,hansen10,phillips,tikhonov63,tikhonov77},
where $L$ is a regularization matrix, and its suitable choice is based on
a-prior information on $x_{true}$. Typically, $L$ is either the identity matrix
$I$ or the scaled discrete approximation of a first or second order derivative
operator. If $L=I$, the identity matrix,
\eqref{posed} reduces to standard-form regularization
in 2-norm, and \eqref{tikhonov} is standard-form Tikhonov regularization,
both of which are {\em 2-norm filtering} regularization problems.

We are concerned with the case $L=I$ in this paper. If $L\not=I$,
\eqref{posed} and \eqref{tikhonov}, in principle,
can be transformed into standard-form
problems \cite{hansen98,hansen10}. In this case, the solutions of \eqref{eq1},
\eqref{posed} and \eqref{tikhonov} can be
fully analyzed by the singular value decomposition (SVD) of $A$. Let
\begin{equation}\label{eqsvd}
  A=U\left(\begin{array}{c} \Sigma \\ \mathbf{0} \end{array}\right) V^{T}
\end{equation}
be the SVD of $A$,
where $U = (u_1,u_2,\ldots,u_m)\in\mathbb{R}^{m\times m}$ and
$V = (v_1,v_2,\ldots,v_n)\in\mathbb{R}^{n\times n}$ are orthogonal,
$\Sigma = {\rm diag} (\sigma_1,\sigma_2,\ldots,\sigma_n)\in\mathbb{R}^{n\times n}$
with the singular values
$\sigma_1>\sigma_2 >\cdots >\sigma_n>0$ assumed to be simple
throughout this paper, and the superscript $T$
denotes the transpose of a matrix or vector. Then
\begin{equation}\label{eq4}
  x_{naive}=\sum\limits_{i=1}^{n}\frac{u_i^{T}b}{\sigma_i}v_i =
  \sum\limits_{i=1}^{n}\frac{u_i^{T}b_{true}}{\sigma_i}v_i +
  \sum\limits_{i=1}^{n}\frac{u_i^{T}e}{\sigma_i}v_i
  =x_{true}+\sum\limits_{i=1}^{n}\frac{u_i^{T}e}{\sigma_i}v_i
\end{equation}
with $\|x_{true}\|=\|A^{\dagger}b_{true}\|=
\left(\sum_{i=1}^n\frac{|u_i^Tb_{true}|^2}{\sigma_i^2}\right)^{1/2}\leq C$.

The discrete Picard condition means that, on average, the Fourier coefficient
$|u_i^{T}b_{true}|$ decays faster than $\sigma_i$ and enables
regularization to compute useful approximations to $x_{true}$. The
following popular simplifying
model is used throughout Hansen's books
\cite{hansen98,hansen10} and the references therein as well as the current paper:
\begin{equation}\label{picard}
  | u_i^T b_{true}|=\sigma_i^{1+\beta},\ \ \beta>0,\ i=1,2,\ldots,n,
\end{equation}
where $\beta$ is a model parameter that controls the decay rates of
$| u_i^T b_{true}|$.

The Gaussian white noise $e$ has a number of attractive properties which
play a critical role in the regularization analysis: Its covariance matrix
is $\eta^2 I$, the expected values $\mathcal{E}(\|e\|^2)=m \eta^2$ and
$\mathcal{E}(|u_i^Te|)=\eta,\,i=1,2,\ldots,n$, so that
$\|e\|\approx \sqrt{m}\eta$ and $|u_i^Te|\approx \eta,\
i=1,2,\ldots,n$; see, e.g., \cite[p.70-1]{hansen98} and \cite[p.41-2]{hansen10}.
The noise $e$ thus affects $u_i^Tb,\ i=1,2,\ldots,n,$ {\em more or less equally}.
With \eqref{picard}, relation \eqref{eq4} shows that for large singular values
the signal terms $|{u_i^{T}b_{true}}|/{\sigma_i}$ are dominant relative to
the noise terms $|u_i^{T}e|/{\sigma_i}$, that is, the $\sigma_i^{\beta}$ are
considerably bigger than the $\eta/\sigma_i$. Once
$| u_i^T b_{true}| \leq | u_i^T e|$ for small singular values, the noise
$e$ dominates $| u_i^T b|$, and the terms $\frac{| u_i^T b|}{\sigma_i}\approx
\frac{|u_i^{T}e|}{\sigma_i}$ overwhelm $x_{true}$
and thus must be filtered out. The transition or cutting-off point
$k_0$ is such that
\begin{equation}\label{picard1}
| u_{k_0}^T b|\approx | u_{k_0}^T b_{true}|> | u_{k_0}^T e|\approx
\eta, \ | u_{k_0+1}^T b|
\approx | u_{k_0+1}^Te|
\approx \eta;
\end{equation}
see \cite[p.42, 98]{hansen10} and a similar description \cite[p.70-1]{hansen98}.
In this sense, the $\sigma_k$ are divided into the $k_0$ large ones and the $n-k_0$
small ones.

The truncated SVD (TSVD) method \cite{hansen90,hansen98,hansen10} is a reliable
and effective method for solving small to medium sized \eqref{posed},
and it deals with a sequence of problems
\begin{equation}\label{tsvd}
\min\|x\| \ \ \mbox{subject to}\ \
 \|A_kx-b\|=\min
\end{equation}
starting with $k=1$ onwards, where $A_k=U_k\Sigma_k V_k^T$
is the best rank $k$ approximation to $A$ with respect to the 2-norm
with $U_k=(u_1,\ldots,u_k)$, $V_k=(v_1,\ldots,v_k)$ and $\Sigma_k=
{\rm diag}(\sigma_1,\ldots,\sigma_k)$; it holds that
$\|A-A_k\|=\sigma_{k+1}$ \cite[p.12]{bjorck96}.
The solution to \eqref{tsvd} is
$
x_{k}^{tsvd}=A_k^{\dagger}b,
$
called the TSVD regularized solution, which is the minimum 2-norm solution to
$
\min\limits_{x\in \mathbb{R}^{n}}\|A_kx-b\|
$
that replaces $A$ by $A_k$ in \eqref{eq1}.

Based on the above properties of the Gaussian white noise $e$, it is known
from \cite[p.70-1]{hansen98} and
\cite[p.71,86-8,95]{hansen10} that the TSVD solutions
\begin{equation}\label{solution}
  x^{tsvd}_k=A_k^{\dagger}b=\left\{\begin{array}{ll}
  \sum\limits_{i=1}^{k}\frac{u_i^{T}b}{\sigma_i}{v_i}\approx
  \sum\limits_{i=1}^{k}\frac{u_i^{T}b_{true}}
{\sigma_i}{v_i},\ \ \ &k\leq k_0;\\ \sum\limits_{i=1}^{k}\frac{u_i^{T}b}
{\sigma_i}{v_i}\approx
\sum\limits_{i=1}^{k_0}\frac{u_i^{T}b_{true}}{\sigma_i}{v_i}+
\sum\limits_{i=k_0+1}^{k}\frac{u_i^{T}e}{\sigma_i}{v_i},\ \ \ &k>k_0,
\end{array}\right.
\end{equation}
and $x_{k_0}^{tsvd}$ is the best TSVD regularized solution of \eqref{eq1};
we have
$\|x_{true}-x_{k_0}^{tsvd}\|=\min_{k=1,2,\ldots,n}\|x_{true}-x_k^{tsvd}\|$,
which balances the
regularization error $(A^{\dagger}-A_k^{\dagger})b_{true}$
and the perturbation error $A_k^{\dagger}e$ optimally, and
$\|Ax_k^{tsvd}-b\|\approx \|e\|$ stabilizes for $k$ not close to $n$
after $k>k_0$. The index $k$ plays the role of the
regularization parameter.

Tikhonov regularization \eqref{tikhonov} with $L=I$ is a filtered SVD method.
For each $\lambda$, the solution $x_{\lambda}$
satisfies $(A^TA+\lambda^2 I)x_{\lambda}=A^Tb$, which replaces
the ill-conditioned $A^TA$ in
normal equation of \eqref{eq1} by $A^TA+\lambda^2 I$,
and has a filtered SVD expansion
\begin{equation}\label{eqfilter}
  x_{\lambda} = \sum\limits_{i=1}^{n}f_i\frac{u_i^{T}b}{\sigma_i}v_i,
\end{equation}
where the $f_i=\frac{\sigma_i^2}{\sigma_i^2+\lambda^2}$ are filters.
The error
$x_{\lambda}-x_{true}$ can be written as the sum of the regularization error
$\left((A^TA+\lambda^2 I)^{-1}A^T-A^{\dagger}\right)b_{true}$ and
the perturbation error $(A^TA+\lambda^2 I)^{-1}A^Te$,
and an optimal $\lambda_{opt}$ is such that
$\|x_{true}-x_{\lambda_{opt}}\|=\min_{\lambda\geq 0}\|x_{true}-x_{\lambda}\|$
and balances these two errors \cite{hansen98,hansen10,kirsch,vogel02}.
In the spirit of $x_{k_0}^{tsvd}$, the best Tikhonov regularized
solution $x_{\lambda_{opt}}$ retains the $k_0$ dominant SVD components
and dampens the other $n-k_0$ small SVD components as much as
possible \cite{hansen98,hansen10}, that is, $\lambda_{opt}$ must be such that
$f_i=\mathcal{O}(1)$, $i=1,2,\ldots,k_0$
and $f_i/\sigma_i\approx 0$, $i=k_0+1,\ldots,n$.
Therefore, it is expected that $x_{k_0}^{tsvd}$ and
$x_{\lambda_{opt}}$ have very similar accuracy. Indeed,
it has been observed and justified that
these two regularized solutions essentially have the minimum 2-norm error;
see \cite{hansen90b}, \cite[p.109-11]{hansen98},
\cite[Sections 4.2 and 4.4]{hansen10} and \cite{varah79}.

As a matter of fact, there is solid mathematical
theory on the TSVD method and standard-form Tikhonov regularization:
for an underlying linear compact equation $Kx=g$, e.g., \eqref{eq2},
with the noisy $g$ and true solution
$x_{true}$, under the source condition that its solution $x_{true}
\in \mathcal{R}(K^*)$ or $x_{true}\in \mathcal{R}(K^*K)$, the range of
the adjoint $K^*$ of $K$ or that of $K^*K$,
the errors of $x_{k_0}^{tsvd}$ and
$x_{\lambda_{opt}}$ are {\em order optimal, i.e., the same order
as the worst-case error} \cite[p.13,18,20,32-40]{kirsch},
\cite[p.90]{natterer} and \cite[p.7-12]{vogel02}. These conclusions
carries over to \eqref{eq1} \cite[p.8]{vogel02}. Therefore,
when only deterministic 2-norm filtering methods are taken into account,
both $x_{\lambda_{opt}}$ and $x_{k_0}^{tsvd}$
are best possible solutions to \eqref{eq1} under the above
assumptions. More generally,
for $x_{true}\in \mathcal{R}((K^*K)^{\beta/2})$
with any $\beta>0$, the error of the TSVD solution $x_{k_0}^{tsvd}$ is always
order optimal, while $x_{\lambda_{opt}}$ is best possible only for
$\beta\leq 2$; see  \cite[Chap. 4-5]{engl00} for details.

As a consequence, we can take $x_{k_0}^{tsvd}$ as the reference standard
when assessing the accuracy of the 2-norm filtering
best regularized solution obtained
by a regularization method. In other words,
we take the TSVD method as reference standard
to evaluate the regularization ability of a deterministic 2-norm filtering
regularization method.

A number of parameter-choice methods have been developed for finding
$\lambda_{opt}$ or $k_0$, such as the discrepancy principle \cite{morozov},
the L-curve criterion, whose use goes back to Miller \cite{miller} and
Lawson and Hanson \cite{lawson} and is termed much later and studied in detail
in \cite{hansen92,hansen93}, the generalized cross validation
(GCV) \cite{golub79,wahba}, and the method based on error estimation
\cite{gfrerer87,raus85}; see, e.g.,
\cite{bauer11,engl00,hansen98,hansen10,kern,kilmer03,kindermann,neumaier98,vogel02}
for numerous comparisons. Each of these methods
has its own merits and disadvantages, and
none is absolutely reliable for all discrete ill-posed problems.
For example, some of them may fail to find accurate approximations to $\lambda_{opt}$;
see \cite{hanke96a,vogel96} for an analysis
on the L-curve criterion method and \cite{hansen98} for some other parameter-choice
methods.

For $A$ large, the TSVD method and the Tikhonov regularization
method are generally too demanding, and only iterative regularization
methods are computationally viable. Krylov solvers are
a major class of iterative methods for solving a large scale \eqref{eq1},
and they project \eqref{eq1}
onto a sequence of low dimensional Krylov subspaces
and computes iterates to approximate $x_{true}$
\cite{aster,engl00,gilyazov,hanke95,hansen98,hansen10,kirsch}.
Of them, the CGLS (or CGNR) method,
which implicitly applies the CG
method \cite{golub89,hestenes} to $A^TAx=A^Tb$,
and its mathematically equivalent LSQR algorithm \cite{paige82}
have been most commonly used. The Krylov solvers CGME
(or CGNE) \cite{bjorck96,bjorck15,craig,hanke95,hanke01} and
LSMR \cite{bjorck15,fong} have been also choices, which amount to the
CG method applied to $\min\|AA^Ty-b\|$ or $AA^Ty=b$
with $x=A^Ty$ and MINRES \cite{paige75}
applied to $A^TAx=A^Tb$, respectively. These Krylov solvers have been
intensively studied and known to have general regularizing
effects \cite{aster,eicke,gilyazov,hanke95,hanke01,hansen98,hansen10,hps16,hps09}
and exhibit semi-convergence \cite[p.89]{natterer};
see also \cite[p.314]{bjorck96}, \cite[p.733]{bjorck15},
\cite[p.135]{hansen98} and \cite[p.110]{hansen10}: the iterates
converge to $x_{true}$ in an initial stage; afterwards the
noise $e$ starts to deteriorate the iterates so that they start to diverge
from $x_{true}$ and instead converge to $x_{naive}$.
If we stop at the right time, then, in principle,
we have a regularization method, where the iteration number plays the
role of the regularization parameter. Semi-convergence is not only due to
the increasingly ill-conditioning of the projected problem,
but also to the fact that the noise progressively enters the approximation
subspace \cite{hps09}.

The behavior of ill-posed problems critically depends on the decay rate of
$\sigma_j$. The following characterization of the degree of ill-posedness
of \eqref{eq1} was introduced in \cite{hofmann86}
and has been widely used \cite{aster,engl00,hansen98,hansen10,mueller}:
if $\sigma_j=\mathcal{O}(\rho^{-j})$ with $\rho>1$,
$j=1,2,\ldots,n$, then \eqref{eq1} is severely ill-posed;
if $\sigma_j=\mathcal{O}(j^{-\alpha})$, then \eqref{eq1}
is mildly or moderately ill-posed for $\frac{1}{2}<\alpha\le1$ or $\alpha>1$.
Here for mildly ill-posed problems we add the requirement
$\alpha>\frac{1}{2}$, which does not appear in \cite{hofmann86}
but must be met for a linear compact operator equation \cite{hanke93,hansen98}.

The regularizing effects of CG type methods were
discovered in \cite{johnsson,squire,tal}.
Johnsson \cite{johnsson} had given a heuristic explanation on the success
of CGLS. Based on these works, on page 13 of \cite{bjorck79},
Bj\"{o}rck and Eld\'{e}n in their 1979 survey foresightedly expressed a fundamental
concern on CGLS (and LSQR): {\em More
research is needed to tell for which problems this approach will work, and
what stopping criterion to choose.} See also \cite[p.145]{hansen98}.
Hanke and Hansen \cite{hanke93} and Hansen \cite{hansen07} have addressed that
a strict proof of the regularizing properties of conjugate gradients is
extremely difficult.

An enormous effort has been made to the study of
regularizing effects of LSQR and CGLS; see
\cite{firro97,gilyza86,hanke95,hanke01,hansen98,hansen10,hps16,
hps09,huangjia,nemi,nolet,paige06,scales,vorst90}, many of which
concern the {\em asymptotic} behavior of the errors
of $x_{\lambda_{opt}}$ and $x_{k_0}^{tsvd}$ as the noise $e$, which assumes
no specific property, approaches zero
in the Hilbert space setting.
Our concern is to leave the {\em Gaussian white}
noise $e$ {\em fixed} and considers how the solution
by LSQR and CGLS behaves as the regularization parameter varies in
the {\em finite} dimensional space. Therefore, our analysis approach
and results are {\em non}-asymptotic and different.
It has long been well known \cite{hanke93,hansen98,hansen07,hansen10}
and will also be elaborated in this paper
that provided that the singular values of projected matrices
in LSQR, called the Ritz values, always approximate
the large singular values of $A$ in natural order until
semi-convergence, the best regularized solution
obtained by LSQR is as accurate as $x_{k_0}^{tsvd}$. Such
convergence is thus desirable. However, Hanke and Hansen \cite{hanke93},
Hansen \cite{hansen98,hansen07,hansen10} and some others,
e.g., Gazzola and Novati \cite{gazzola15},
address the difficulties to prove the convergence in this order.
Hitherto there has been no general
definitive and quantitative result on whether or not
the Ritz values converge in this order for the three kinds of ill-posed
problems.

Precisely, we now introduce such a definition:
For the 2-norm filtering regularization problem \eqref{posed},
if a regularized solution is as accurate as
$x_{k_0}^{tsvd}$,
then it is called a 2-norm filtering best possible regularized solution.
If a regularization method can compute such a best possible one,
then it is said to have
the {\em full} regularization in the sense of
2-norm filtering. Otherwise, it is said to have
only the {\em partial} regularization.

In order to obtain a 2-norm filtering best possible solution of \eqref{eq1},
LSQR and CGLS have been commonly combined with some explicit regularization
\cite{aster,berisha,gazzola18,hansen98,hansen10}. CGLS is combined with the
standard-form Tikhonov regularization, and it
solves $(A^TA+\lambda^2 I)x=A^Tb$ for several regularization
parameters $\lambda$ and
picks up a best solution \cite{aster,frommer99}.
The hybrid LSQR variants have been advocated by Bj\"{o}rck and Eld\'{e}n
\cite{bjorck79} and O'Leary and Simmons \cite{oleary81}, and improved and
developed by Bj\"orck \cite{bjorck88}, Bj\"{o}rck, Grimme and
van Dooren \cite{bjorck94}, and Renaut, Vatankhah, and Ardestani \cite{renaut}.
They first project \eqref{eq1} onto Krylov
subspaces and then regularize the projected problem explicitly at each iteration.
They aim to remove the effects
of small Ritz values and expands Krylov subspaces until they
captures all the needed right singular vectors of $A$.
The hybrid LSQR, CGME and LSMR have been intensively studied in, e.g.,
\cite{bazan10,bazan14,berisha,chung08,chung15,hanke01,hanke93,lewis09,neuman12,renaut}
and \cite{aster,hansen10,hansen13}. For further information on
hybrid methods, we refer to \cite{gazzola-online} and the references therein.

If an iterative solver is theoretically proved to
have the full regularization, one can stop it after
its semi-convergence is practically identified. To echo the concern of
Bj\"{o}rck and Eld\'{e}n, by the definition of
the full or partial regularization, our question is:
{\em  Do LSQR, CGLS, LSMR and CGME have the full or partial regularization for
severely, moderately and mildly ill-posed problems?}
As we have seen, there has been no definitive answer to this long-standing
fundamental question hitherto.

LSQR, CGME and LSMR are common in that, at iteration $k$, they
are mathematically based on the same $k$-step Lanczos bidiagonalization
process but compute {\em different} iterates from the {\em same}
$k$ dimensional right Krylov subspace. Remarkably, note that if the left and right
subspaces are the dominant left and right singular subspaces $span\{U_{k+1}\}$
(or $span\{U_k\}$) and $span\{V_k\}$ of $A$
then the Ritz values of $A$ with respect to them
are exactly the first $k$ large singular values of $A$. Therefore,
whether or not the Ritz values converge to the large singular values of $A$ in
natural order critically depends on
how the underlying $k$ dimensional right Krylov subspace approaches
$span\{V_k\}$. This paper concerns the fundamental problem
that these methods face: How does the underlying
$k$ dimensional right Krylov subspace approximate $span\{V_k\}$?
Accurate solutions of this problem play a central role in analyzing the
regularization ability of the mentioned four methods and in ultimately
determining if each method has the full regularization.
We will establish a general $\sin\Theta$ theorem
for the 2-norm distances between these two subspaces and
derive accurate estimates on them for
the three kinds of ill-posed problems. We
notice that the $\sin\Theta$ theorem involves some crucial quantities
used to study the regularizing effects of
LSQR \cite[p.150-2]{hansen98}, but there were no estimates for them
there and in the literature.

In Section \ref{methods}, we describe the
Lanczos bidiagonalization process and the LSQR method,
and make an introductory analysis.
In Section~\ref{argument} we make an analysis on the
regularizing effects of LSQR and establish a basic result on
its semi-convergence. In Section \ref{sine},
we establish the $\sin\Theta$ theorem for the 2-norm
distance between the underlying $k$ dimensional Krylov subspace and $span\{V_k\}$,
and derive accurate estimates on them for the
three kinds of ill-posed problems, which include accurate estimates for
those key quantities in \cite[p.150-2]{hansen98}.
In Section~\ref{manif} we consider the effects
of the $\sin\Theta$ theorem on the behavior of the
smallest Ritz values involved in LSQR.
We report a number of numerical examples to confirm our theory.
Finally, we summarize the paper in Section \ref{concl}.

Throughout the paper, we denote by
$\mathcal{K}_{k}(C, w)= span\{w,Cw,\ldots,C^{k-1}w\}$
the $k$ dimensional Krylov subspace generated
by the matrix $\mathit{C}$ and the vector $\mathit{w}$, and by $I$ and the
bold letter $\mathbf{0}$ the identity matrix
and the zero matrix with orders clear from the context, respectively.
For the matrix $B=(b_{ij})$, we define $|B|=(|b_{ij}|)$,
and for $|C|=(|c_{ij}|)$, $|B|\leq |C|$ means
$|b_{ij}|\leq |c_{ij}|$ componentwise.

\section{The LSQR algorithm}\label{methods}

The LSQR algorithm is based on the Lanczos bidiagonalization process,
which computes two orthonormal bases $\{q_1,q_2,\dots,q_k\}$ and
$\{p_1,p_2,\dots,p_{k+1}\}$  of $\mathcal{K}_{k}(A^{T}A,A^{T}b)$ and
$\mathcal{K}_{k+1}(A A^{T},b)$  for $k=1,2,\ldots,n$,
respectively. We describe the process as Algorithm 1.

{\bf Algorithm 1: \  $k$-step Lanczos bidiagonalization process}
\begin{itemize}
\item Take $ p_1=b/\|b\| \in \mathbb{R}^{m}$, and define $\beta_1{q_0}=\mathbf{0}$
with $\beta_1=\|b\|$.

\item For $j=1,2,\ldots,k$\\
 (i) $r = A^{T}p_j - \beta_j{q_{j-1}}$\\
(ii) $\alpha_j = \|r\|;q_j = r/\alpha_j$\\
(iii) $z = Aq_j - \alpha_j{p_{j}}$\\
(iv) $\beta_{j+1} = \|z\|;p_{j+1} = z/\beta_{j+1}.$
\end{itemize}

Algorithm 1 can be written in the matrix form
\begin{align}
  AQ_k&=P_{k+1}B_k,\label{eqmform1}\\
  A^{T}P_{k+1}&=Q_{k}B_k^T+\alpha_{k+1}q_{k+1}(e_{k+1}^{(k+1)})^{T},\label{eqmform2}
\end{align}
where $e_{k+1}^{(k+1)}$ is the $(k+1)$-th canonical basis vector of
$\mathbb{R}^{k+1}$, $P_{k+1}=(p_1,p_2,\ldots,p_{k+1})$,
$Q_k=(q_1,q_2,\ldots,q_k)$, and
\begin{equation}\label{bk}
  B_k = \left(\begin{array}{cccc} \alpha_1 & & &\\ \beta_2 & \alpha_2 & &\\ &
  \beta_3 &\ddots & \\& & \ddots & \alpha_{k} \\ & & & \beta_{k+1}
  \end{array}\right)\in \mathbb{R}^{(k+1)\times k}.
\end{equation}
It is known from \eqref{eqmform1} that
\begin{equation}\label{Bk}
B_k=P_{k+1}^TAQ_k.
\end{equation}
The singular values $\theta_i^{(k)},\ i=1,2,\ldots,k$ of $B_k$,
called the Ritz values of $A$ with
respect to the left and right subspaces $span\{P_{k+1}\}$ and $span\{Q_k\}$,
are all simple as $\alpha_i>0,\ \beta_{i+1}>0,\ i=1,2,\ldots,k$ provided that the
algorithm does not break down.

Write
$
\mathcal{V}_k^R=\mathcal{K}_k(A^TA,A^Tb)
$
and $\mathcal{V}_k=span\{V_k\}$. LSQR, CGME and LSMR  LSQR, CGME and LSMR
are based on the same Lanczos bidiagonalization process
but extract {\em different} iterates
from $\mathcal{V}_k^R$. We take LSQR as example.
At iteration $k$, LSQR solves the problem
$$
\|Ax_k^{lsqr}-b\|=\min_{x\in \mathcal{V}_k^R}
\|Ax-b\|
$$
for the iterate
\begin{equation}\label{yk}
x_k^{lsqr}=Q_ky_k^{lsqr} \ \ \mbox{with}\ \
y_k^{lsqr}=\arg\min\limits_{y\in \mathbb{R}^{k}}\|B_ky-\beta_1 e_1^{(k+1)}\|
  =\beta_1  B_k^{\dagger} e_1^{(k+1)},
\end{equation}
where $e_1^{(k+1)}$ is the first canonical basis vector of $\mathbb{R}^{k+1}$,
and the residual norm $\|Ax_k^{lsqr}-b\|=\|B_ky_k^{lsqr}-\beta_1 e_1^{(k+1)}\|$
and the solution norm $\|x_k^{lsqr}\|=\|y_k^{lsqr}\|$
decreases and increases monotonically with respect to $k$, respectively.

From $\beta_1 e_1^{(k+1)}=P_{k+1}^T b$ and \eqref{yk}, we have
\begin{equation}\label{xk}
x_k^{lsqr}=Q_k B_k^{\dagger} P_{k+1}^Tb,
\end{equation}
that is, $x_k^{lsqr}$ is the minimum 2-norm solution to the perturbed
problem that replaces $A$ in \eqref{eq1} by its rank $k$ approximation
$P_{k+1}B_k Q_k^T$. So LSQR solves a sequence of problems
\begin{equation}\label{lsqrreg}
\min\|x\| \ \ \mbox{ subject to }\ \ \|P_{k+1}B_kQ_k^Tx-b\|=\min
\end{equation}
for the regularized solutions $x_k^{lsqr}$ of \eqref{eq1} starting with $k=1$
onwards. Recall the TSVD method (cf. \eqref{tsvd}) and
that the best rank $k$ approximation $A_k$ to $A$ satisfies
$\|A-A_k\|=\sigma_{k+1}$. Consequently,
if $P_{k+1}B_k Q_k^T$ is a near best rank $k$ approximation
to $A$ with an approximate accuracy $\sigma_{k+1}$ and
the $k$ singular values of $B_k$ approximate the first $k$
large ones of $A$ in natural order for $k\leq k_0$,
then LSQR and the TSVD method are related naturally and
closely because (i) $x_k^{tsvd}$ and $x_k^{lsqr}$ are the
solutions to the two perturbed problems of \eqref{eq1} that replace $A$ by
its two rank $k$ approximations with the same quality, respectively;
(ii) $x_k^{tsvd}$ and $x_k^{lsqr}$ solve the two essentially
same regularization problems \eqref{tsvd} and \eqref{lsqrreg}, respectively.
As a consequence,
the LSQR iterate $x_{k_0}^{lsqr}$ is as accurate as $x_{k_0}^{tsvd}$,
and LSQR has the full regularization. Therefore,
that the $P_{k+1}B_k Q_k^T$ are near best rank $k$ approximations
to $A$  and the $k$ singular values of $B_k$ approximate the
large ones of $A$ in natural order for $k=1,2,\ldots,k_0$ are
{\em sufficient} conditions for which LSQR has the full regularization.

However, we {\em must remind} that the near best rank $k$ approximations
and the approximations of the singular values of $B_k$ to
the large singular values $\sigma_i$ in natural order are {\em not necessary}
conditions for the full regularization of LSQR. It is well possible that
LSQR has the full regularization even though these conditions
are not satisfied, as will be confirmed numerically later.

\section{An elementary analysis on the regularizing effects of LSQR}
\label{argument}

The following result (cf., e.g., van der Sluis and
van der Vorst \cite{vorst86}) has been widely used, e.g.,
in Hansen \cite{hansen98}, to illustrate the regularizing effects of LSQR
and CGLS.

\begin{prop}\label{help}
LSQR applied to \eqref{eq1} with the starting vector $p_1=b/\|b\|$ and CGLS
applied to $A^TAx=A^Tb$ with the zero starting vector generate the same iterates
\begin{equation}\label{eqfilter2}
  x_k^{lsqr}=\sum\limits_{i=1}^nf_i^{(k)}\frac{u_i^{T}b}{\sigma_i}v_i,\
  k=1,2,\ldots,n,
\end{equation}
where the filters
\begin{equation}\label{filter}
f_i^{(k)}=1-\prod\limits_{j=1}^k\frac{(\theta_j^{(k)})^2-\sigma_i^2}
{(\theta_j^{(k)})^2},\ i=1,2,\ldots,n,
\end{equation}
and the $\theta_j^{(k)}$ are the singular values of $B_k$
labeled as $\theta_1^{(k)}>\theta_2^{(k)}>\cdots>\theta_k^{(k)}$.
\end{prop}

This proposition shows that $x_k^{lsqr}$ has a filtered SVD expansion. It is
easily justified that if
the $k$ Ritz values $\theta_j^{(k)}$ approximate the first $k$ singular values
$\sigma_j$ of $A$ in natural order then the filters $f_i^{(k)}\approx 1$ for $
i=1,2,\ldots,k$ and $f_i^{(k)}$ monotonically approach zero
for $i=k+1,\ldots,n$. This indicates that if the $\theta_j^{(k)}$
approximate the first $k$ singular values
$\sigma_j$ of $A$ in natural order for $k=1,2,\ldots,k_0$ then
LSQR definitely has the full regularization.

Regarding the semi-convergence of LSQR and
the TSVD method, we present the following basic result.

\begin{theorem}\label{semicon}
The semi-convergence of LSQR occurs at some iteration
$$
k^*\leq k_0.
$$
If the Ritz values $\theta_j^{(k)}$ do not converge to the first
$k$ large singular values $\sigma_j$ of $A$ in natural order for some
$k\leq k^*$, then $k^*<k_0$, and vice versa.
\end{theorem}

{\em Proof}.
Applying the Cauchy's {\em strict}
interlacing theorem \cite[p.198, Corollary 4.4]{stewartsun} to the singular
values of $B_k$ and $B_n$, we always have
$$
\theta_i^{(k)}< \sigma_i,\ i=1,2,\ldots,k.
$$
Therefore, at iteration $k_0+1$ one must have $\theta_{k_0+1}^{(k_0+1)}<
\sigma_{k_0+1}$. As a result, if the $\theta_i^{(k)}$ approximate the large
$\sigma_i$ in natural order for $k=1,2,\ldots,k_0$, then by \eqref{eqfilter2}
and \eqref{filter} we have $f_{k_0+1}^{(k_0+1)}\approx 1$, meaning
that $x_{k_0+1}^{lsqr}$ must be deteriorated and the semi-convergence
of LSQR must occur at iteration $k^*=k_0$.

If the $\theta_j^{(k)}$ do not converge to the large singular values of
$A$ in natural order and $\theta_k^{(k)}<\sigma_{k_0+1}$ appears
for some iteration $k\leq k_0$ for the first time, then, by the strict
interlacing property of singular values of $B_{k-1}$ and $B_k$, it holds
that $\theta_{k-1}^{(k)}>\theta_{k-1}^{(k-1)}>\theta_k^{(k)}$, meaning
that $\theta_{k-1}^{(k)}>\theta_{k-1}^{(k-1)}>\sigma_{k_0+1}$. In this case,
$x_k^{lsqr}$ is already deteriorated by the noise $e$ before iteration $k$:
Suppose that $\sigma_{j^*}<\theta_k^{(k)}<\sigma_{k_0+1}$ with
$j^*$ the smallest integer
$j^*>k_0+1$. Then we can easily justify from \eqref{filter} that
$f_i^{(k)}\in (0,1)$ and tends to zero
monotonically for $i=j^*,j^*+1,\ldots,n$, but
we have
$$
\prod\limits_{j=1}^k\frac{(\theta_j^{(k)})^2-\sigma_i^2}
{(\theta_j^{(k)})^2}=\frac{(\theta_k^{(k)})^2-\sigma_i^2}
{(\theta_k^{(k)})^2}\prod\limits_{j=1}^{k-1}
\frac{(\theta_j^{(k)})^2-\sigma_i^2}{(\theta_j^{(k)})^2}\leq 0,
\ i=k_0+1,\ldots,j^*-1
$$
since the first factor is non-positive and the second factor is positive by
noticing that $\theta_j^{(k)}>\sigma_i$, $j=1,2,\ldots,k-1$ for
$i=k_0+1,\ldots,j^*-1$.
Hence $f_i^{(k)}\geq 1$ for $i=k_0+1,\ldots, j^*-1$, showing
that $x_k^{lsqr}$ has been deteriorated by the noise $e$ and
the semi-convergence of LSQR has occurred at some iteration $k^*<k_0$.

On the other hand, if the semi-convergence of LSQR occurs at iteration
$k^*<k_0$, then the $k^*$ Ritz values
$\theta_j^{(k^*)}$ must not approximate the first $k^*$ large singular
values $\sigma_j$ of $A$ in natural order. Otherwise, notice that
$\theta_{k^*}^{(k^*)}\in (\sigma_{k^*+1},\sigma_{k^*})$
means $\theta_{k^*}^{(k^*)}>\sigma_{k^*+1}\geq \sigma_{k_0}$, which indicates
that the semi-convergence of LSQR
does not yet occur at iteration $k^*$, a contradiction.
\qquad\endproof

If the semi-convergence of LSQR occurs at iteration $k^*<k_0$,
the regularizing effects of LSQR is much more complicated, and
there has been no definitive result on the full or partial regularization of
LSQR. This problem will be our future concern.

The standard $k$-step Lanczos bidiagonalization method \cite{bjorck96,bjorck15}
computes the $k$ Ritz values $\theta_j^{(k)}$, which is mathematically
equivalent to the symmetric Lanczos method for
the eigenvalue problem of $A^TA$ starting with $q_1=A^Tb/\|A^Tb\|$;
see \cite{bai,bjorck96,bjorck15,parlett,vorst02} or \cite{jia03,jia10} for
several variations that are based on standard, harmonic,
and refined projection \cite{bai,stewart01,vorst02}
or a combination of them \cite{jia05}. An attractive feature is that,
for general singular value distribution and $b$, some
$\theta_i^{(k)}$ become good approximations to the extreme
singular values of $A$ as $k$ increases.
If large singular values are well separated but
small singular values are clustered, the large $\theta_j^{(k)}$
converge fast but small $\theta_j^{(k)}$ show up late and converge slowly.

For \eqref{eq1}, since the singular values $\sigma_j$ of $A$ decay
to zero, $A^Tb=\sum_{j=1}^n\sigma_j (u_j^Tb) v_j$
contains more information on dominant right
singular vectors than on the ones corresponding to small singular values.
Therefore, the Krylov subspace $\mathcal{V}_k^R$ with $A^Tb$ as the starting
vector is expected to contain richer
information on the first $k$ right singular vectors $v_j$ than on the
other $n-k$ ones. In the meantime,
notice that $A$ has many small singular values clustered at zero.
Due to these two basic facts, all the $\theta_j^{(k)}$ are expected to
approximate the large singular values of $A$ in natural order until some
iteration $k$. In this case, the iterates $x_k^{lsqr}$ mainly consists
of the $k$ dominant SVD components of $A$.
This is why LSQR and CGLS have general regularizing effects; see,
e.g., \cite{aster,hansen98,hansen08,hansen10}.

Unfortunately, the above arguments are purely qualitative and not rigorous.
They do not give any hints on the size of $k$ for possible desired
convergence in natural order for any given \eqref{eq1}.
As has been addressed previously, proving how the Ritz values converge
is extremely difficult.

\section{The $\sin\Theta$ theorem and its estimates
for the 2-norm distances between $\mathcal{V}_k^R$
and $\mathcal{V}_k$} \label{sine}

As can be seen from Sections~\ref{methods}--\ref{argument},
a complete understanding of the regularization
of LSQR includes accurate solutions of the following problems:
(i) How accurately does the $k$ dimensional right
Krylov subspace $\mathcal{V}_k^R$
approximate the $k$ dimensional dominant right singular subspace
$\mathcal{V}_k$
of $A$? (ii) How accurate is the rank $k$ approximation
$P_{k+1}B_kQ_k^T$ to $A$?
(iii) When do the $\theta_j^{(k)}$
approximate the large $\sigma_j$ in natural order? (iv)
When does at least
a small Ritz value appear, i.e., $\theta_k^{(k)}<\sigma_{k+1}$ for some
$k\leq k^*$?
(v) Does LSQR have the full or partial regularization
when the $k$ Ritz values $\theta_j^{(k)}$ do not approximate the large
$\sigma_j$ in natural order for some $k\leq k^*$?
(vi) What are the counterparts of Problems (i)-(v)
in the case that $A$ has multiple singular values, and what are
the accurate solutions or definitive answers correspondingly?

Problem~(i) is the starting and key point of the other problems, and its
accurate solutions form an absolutely necessary
basis of dealing with the others.
In this paper, we focus on Problem (i) and present accurate results
for the three kinds of ill-posed problems. By them, we will make
an elementary exploration on Problem (iv).
In-depth treatments of Problem (iv) and
the others will given in separate papers.

In terms of the canonical angles $\Theta(\mathcal{X},\mathcal{Y})$ between
two subspaces $\mathcal{X}$ and $\mathcal{Y}$ of equal
dimension (cf. \cite[p.74-5]{stewart01} and \cite[p.43]{stewartsun}),
we first present a general $\sin\Theta$
theorem which measures the 2-norm distance between
$\mathcal{V}_k^R$ and $\mathcal{V}_k$.

\begin{lemma}\label{lemma1}
For $k=1,2,\ldots,n-1$ we have
\begin{align}
\|\sin\Theta(\mathcal{V}_k,\mathcal{V}_k^R)\|&=
\frac{\|\Delta_k\|}{\sqrt{1+\|\Delta_k\|^2}}
\label{deltabound}
\end{align}
with $\Delta_k \in \mathbb{R}^{(n-k)\times k}$ defined by \eqref{defdelta}.
\end{lemma}

{\em Proof.}
Let $U_n=(u_1,u_2,\ldots,u_n)$ whose columns are the
first $n$ left singular vectors of $A$ defined by \eqref{eqsvd}.
Then the Krylov subspace $\mathcal{K}_{k}(\Sigma^2,
\Sigma U_n^Tb)=span\{DT_k\}$ with
\begin{equation*}
  D={\rm diag}(\sigma_j u_j^Tb)\in\mathbb{R}^{n\times n},\ \
  T_k=\left(\begin{array}{cccc} 1 &
  \sigma_1^2&\ldots & \sigma_1^{2k-2}\\
1 &\sigma_2^2 &\ldots &\sigma_2^{2k-2} \\
\vdots & \vdots&&\vdots\\
1 &\sigma_n^2 &\ldots &\sigma_n^{2k-2}
\end{array}\right).
\end{equation*}
Partition the diagonal matrix $D$ and the matrix $T_k$ as
\begin{equation}\label{tk12p}
  D=\left(\begin{array}{cc} D_1 & 0 \\ 0 & D_2 \end{array}\right),\ \ \
  T_k=\left(\begin{array}{c} T_{k1} \\ T_{k2} \end{array}\right),
\end{equation}
where $D_1, T_{k1}\in\mathbb{R}^{k\times k}$. Since $T_{k1}$ is
a Vandermonde matrix with $\sigma_j$ distinct for $j=1,2,\ldots,k$,
it is nonsingular. Therefore, from $\mathcal{K}_{k}(A^{T}A, A^{T}b)=span\{VDT_k\}$
we have
\begin{equation}\label{kry}
\mathcal{V}_k^R=\mathcal{K}_{k}(A^{T}A, A^{T}b)=span
  \left\{V\left(\begin{array}{c} D_1T_{k1} \\ D_2T_{k2} \end{array}\right)\right\}
  =span\left\{V\left(\begin{array}{c} I \\ \Delta_k \end{array}\right)\right\},
\end{equation}
where
\begin{equation}\label{defdelta}
\Delta_k=D_2T_{k2}T_{k1}^{-1}D_1^{-1}\in \mathbb{R}^{(n-k)\times k}.
\end{equation}
Write $V=(V_k, V_k^{\perp})$, and define
\begin{equation}\label{zk}
Z_k=V\left(\begin{array}{c} I \\ \Delta_k \end{array}\right)
=V_k+V_k^{\perp}\Delta_k.
\end{equation}
Then $Z_k^TZ_k=I+\Delta_k^T\Delta_k$, and the columns of
$\hat{Z}_k=Z_k(Z_k^TZ_k)^{-\frac{1}{2}}$
form an orthonormal basis of $\mathcal{V}_k^R$. As a result, from \eqref{zk}
we get an orthogonal direct sum decomposition
\begin{equation}\label{decomp}
\hat{Z}_k=(V_k+V_k^{\perp}\Delta_k)(I+\Delta_k^T\Delta_k)^{-\frac{1}{2}}.
\end{equation}
By the definition of $\Theta(\mathcal{V}_k,\mathcal{V}_k^R)$
and \eqref{decomp}, we obtain
$$
   \|\sin\Theta(\mathcal{V}_k,\mathcal{V}_k^R)\|
   =\|(V_k^{\perp})^T\hat{Z}_k\|
   =\|\Delta_k(I+\Delta_k^T\Delta_k)^{-\frac{1}{2}}\|
   =\frac{\|\Delta_k\|}{\sqrt{1+\|\Delta_k\|^2}},
$$
which proves \eqref{deltabound}.
\qquad\endproof

We remark that it is direct from \eqref{deltabound} to get
\begin{equation}\label{tangent}
\|\tan\Theta(\mathcal{V}_k,\mathcal{V}_k^R)\|=
\|\Delta_k\|.
\end{equation}

\eqref{deltabound} has been established in
\cite[Theorem 2.1]{huangjia}, and we include the proof here
for completeness and for the introduction of notation $\|\Delta_k\|$,
which will be used later.

We now establish accurate estimates for $\|\Delta_k\|$ for severely ill-posed
problems.

\begin{theorem}\label{thm2}
Let the SVD of $A$ be as \eqref{eqsvd}, and assume that \eqref{eq1} is severely
ill-posed with $\sigma_j=\mathcal{O}(\rho^{-j})$ and $\rho>1$, $j=1,2,\ldots,n$.
Then
\begin{align}
\|\Delta_1\|&\leq \frac{\sigma_{2}}{\sigma_1}\frac{\max_{2\leq i\leq n}|u_i^Tb|}{|u_1^Tb|}
\left(1+\mathcal{O}(\rho^{-2})\right),\label{k1}\\
  \|\Delta_k\|&\leq
  \frac{\sigma_{k+1}}{\sigma_k}\frac{\max_{k+1\leq i\leq n}|u_i^Tb|}
  {\min_{1\leq i\leq k}|u_i^Tb|}
  \left(1+\mathcal{O}(\rho^{-2})\right)
  |L_{k_1}^{(k)}(0)|,\ k=2,3,\ldots,n-1,\label{eqres1}
\end{align}
where
\begin{equation}\label{lk}
|L_{k_1}^{(k)}(0)|=\max_{j=1,2,\ldots,k}|L_j^{(k)}(0)|,
\ |L_j^{(k)}(0)|=\prod\limits_{i=1,i\ne j}^k\frac{\sigma_i^2}{|\sigma_j^2-
\sigma_i^2|},\,j=1,2,\ldots,k.
\end{equation}
\end{theorem}

{\em Proof.}
For $k=2,3,\ldots,n-1$,
it is easily justified that the $j$-th column of $T_{k1}^{-1}$ consists of
the coefficients of the $j$-th Lagrange polynomial
\begin{equation*}\label{}
  L_j^{(k)}(\lambda)=\prod\limits_{i=1,i\neq j}^k
  \frac{\lambda-\sigma_i^2}{\sigma_j^2-\sigma_i^2}
\end{equation*}
that interpolates the elements of the $j$-th canonical basis vector
$e_j^{(k)}\in \mathbb{R}^{k}$ at the abscissas $\sigma_1^2,\sigma_2^2
\ldots, \sigma_k^2$. Consequently, the $j$-th column of $T_{k2}T_{k1}^{-1}$ is
\begin{equation}\label{tk12i}
  T_{k2}T_{k1}^{-1}e_j^{(k)}=(L_j^{(k)}(\sigma_{k+1}^2),\ldots,L_j^{(k)}
  (\sigma_{n}^2))^T, \ j=1,2,\ldots,k,
\end{equation}
from which we obtain
\begin{equation}\label{tk12}
  T_{k2}T_{k1}^{-1}=\left(\begin{array}{cccc} L_1^{(k)}(\sigma_{k+1}^2)&
  L_2^{(k)}(\sigma_{k+1}^2)&\ldots & L_k^{(k)}(\sigma_{k+1}^2)\\
L_1^{(k)}(\sigma_{k+2}^2)&L_2^{(k)}(\sigma_{k+2}^2) &\ldots &
L_k^{(k)}(\sigma_{k+2}^2) \\
\vdots & \vdots&&\vdots\\
L_1^{(k)}(\sigma_{n}^2)&L_2^{(k)}(\sigma_{n}^2) &\ldots &L_k^{(k)}(\sigma_{n}^2)
\end{array}\right)\in \mathbb{R}^{(n-k)\times k}.
\end{equation}
Since $|L_j^{(k)}(\lambda)|$ is monotonically
decreasing for $0\leq \lambda<\sigma_k^2$, it is bounded by $|L_j^{(k)}(0)|$.
With this property and the definition of $L_{k_1}^{(k)}(0)$, we obtain
\begin{align}
|\Delta_k|&=|D_2T_{k2}T_{k1}^{-1}D_1^{-1}| \notag \\
&\leq
\left(\begin{array}{cccc}
\frac{\sigma_{k+1}}{\sigma_1}\left|\frac{u_{k+1}^Tb}
{u_1^Tb}\right||L_{k_1}^{(k)}(0)| &\frac{\sigma_{k+1}}{\sigma_2}
\left|\frac{u_{k+1}^Tb}
{u_2^Tb}\right||L_{k_1}^{(k)}(0)| &\ldots&\frac{\sigma_{k+1}}{\sigma_k}
\left|\frac{u_{k+1}^Tb}{u_k^Tb}\right||L_{k_1}^{(k)}(0)| \\
\frac{\sigma_{k+2}}{\sigma_1}\left|\frac{ u_{k+2}^Tb}
{u_1^Tb}\right| |L_{k_1}^{(k)}(0)| &\frac{\sigma_{k+2}}{\sigma_2}
\left|\frac{u_{k+2}^Tb}
{u_2^Tb}\right| |L_{k_1}^{(k)}(0)|&
\ldots &\frac{\sigma_{k+2}}{\sigma_k}\left|\frac{u_{k+2}^Tb}
{u_k^Tb}\right| |L_{k_1}^{(k)}(0)| \\
\vdots &\vdots & &\vdots\\
\frac{\sigma_n}{\sigma_1}\left|\frac{u_n^Tb}
{u_1^Tb}\right| |L_{k_1}^{(k)}(0)| &\frac{\sigma_n}{\sigma_2}\left|\frac{u_n^Tb}
{u_2^Tb}\right| |L_{k_1}^{(k)}(0)|& \ldots &
\frac{\sigma_n}{\sigma_k}\left|\frac{u_n^Tb}{u_k^Tb}\right| |L_{k_1}^{(k)}(0)|
\end{array}
\right) \notag\\
&= |L_{k_1}^{(k)}(0)||\tilde\Delta_k|, \label{amplify}
\end{align}
where
\begin{equation}
|\tilde\Delta_k|=\left|(\sigma_{k+1} u_{k+1}^T b,\sigma_{k+2}u_{k+2}^Tb,
\ldots,\sigma_n u_n^T b)^T
\left(\frac{1}{\sigma_1 u_1^Tb},\frac{1}{\sigma_2 u_2^Tb},\ldots,
\frac{1}{\sigma_k u_k^Tb}\right)\right|  \label{delta1}
\end{equation}
is a rank one matrix. Therefore, by $\|C\|\leq \||C|\|$
(cf. \cite[p.53]{stewart98}), we have
\begin{align}
\|\Delta_k\| &\leq \||\Delta_k|\|\leq |L_{k_1}^{(k)}(0)|
\left\||\tilde\Delta_k|\right\| \notag\\
&=|L_{k_1}^{(k)}(0)|\left(\sum_{j=k+1}^n\sigma_j^2| u_j^Tb|^2\right)^{1/2}
\left(\sum_{j=1}^k \frac{1}{\sigma_j^2| u_j^Tb|^2}\right)^{1/2}.
\label{delta2}
\end{align}

In the following we bound the above two square root factors
separately.

From $\sigma_j=\mathcal{O}(\rho^{-j}),\ j=1,2,\ldots,n$,
for $k=1,2,\ldots,n-1$ we obtain
\begin{align}
\left(\sum_{j=k+1}^n\sigma_j^2| u_j^Tb|^2\right)^{1/2}
&= \sigma_{k+1}\max_{k+1\leq i\leq n}| u_i^Tb| \left(\sum_{j=k+1}^n
\frac{\sigma_j^2| u_j^Tb|^2}{\sigma_{k+1}^2\max_{k+1\leq i\leq n}| u_i^Tb|}\right)^{1/2}
\notag\\
&\leq \sigma_{k+1}\max_{k+1\leq i\leq n}| u_i^Tb| \left(\sum_{j=k+1}^n
\frac{\sigma_j^2}{\sigma_{k+1}^2}\right)^{1/2} \notag\\
&=\sigma_{k+1}\max_{k+1\leq i\leq n}| u_i^Tb|\left(1+\sum_{j=k+2}^n\mathcal{O}
(\rho^{2(k-j)+2})\right)^{1/2}
\notag \\
&=\sigma_{k+1}\max_{k+1\leq i\leq n}| u_i^Tb|\left(1+\mathcal{O}\left(\sum_{j=k+2}^n
\rho^{2(k-j)+2}\right)\right)^{1/2}
\notag \\
&=\sigma_{k+1}\max_{k+1\leq i\leq n}| u_i^Tb|\left(1+ \mathcal{O}\left(\frac{\rho^{-2}}
    {1-\rho^{-2}}\left(1-\rho^{-2(n-k-1)}\right)\right)\right)^{1/2}
\notag \\
&=\sigma_{k+1}\max_{k+1\leq i\leq n}| u_i^Tb|\left(1+\mathcal{O}(\rho^{-2})\right)^{1/2}\notag\\
&=\sigma_{k+1}\max_{k+1\leq i\leq n}| u_i^Tb| \left(1+\mathcal{O}(\rho^{-2})\right)
\label{severe1}
\end{align}
with $1+\mathcal{O}(\rho^{-2})=1$ for $k=n-1$.

For $k=2,3,\ldots,n-1$, we get
\begin{align*}
\left(\sum_{j=1}^k \frac{1}{\sigma_j^2| u_j^Tb|^2}\right)^{1/2}
&=\frac{1}{\sigma_k \min_{1\leq i\leq k}| u_i^T b|}\left(\sum_{j=1}^k\frac{\sigma_k^2
\min_{1\leq i\leq k}| u_i^T b|}
{\sigma_j^2| u_j^Tb|^2}\right)^{1/2}\\
&\leq \frac{1}{\sigma_k \min_{1\leq i\leq k}| u_i^T b|}\left(\sum_{j=1}^k\frac{\sigma_k^2}
{\sigma_j^2}\right)^{1/2} \\
&=\frac{1}{\sigma_k \min_{1\leq i\leq k}| u_i^T b|}\left(1+\mathcal{O}\left(\sum_{j=1}^{k-1}
\rho^{2(j-k)}\right)\right)^{1/2} \\
&=\frac{1}{\sigma_k \min_{1\leq i\leq k}| u_i^T b|}\left(1+\mathcal{O}(\rho^{-2})\right).
\end{align*}
From the above and \eqref{delta2}--\eqref{severe1},
we finally obtain \eqref{eqres1}.

Note that the Lagrange polynomials $L_j^{(k)}(\lambda)$ require $k\geq 2$.
Therefore, we need to treat the case $k=1$ separately.
Note that
$$
T_{k2}=(1,1,\ldots,1)^T,\ D_2T_{k2}=(\sigma_2u_2^Tb,\sigma_3 u_3^Tb,
\ldots,\sigma_n u_n^Tb)^T,\
T_{k1}^{-1}=1,\ D_1^{-1}=\frac{1}{\sigma_1 u_1^Tb}.
$$
Therefore, from \eqref{defdelta} we have
\begin{equation}\label{deltaexp}
\Delta_1=(\sigma_2u_2^Tb,\sigma_3 u_3^Tb,
\ldots,\sigma_n u_n^Tb)^T\frac{1}{\sigma_1 u_1^Tb},
\end{equation}
from which and \eqref{severe1} it is direct to get \eqref{k1}.
\qquad\endproof

A crucial step in proving \eqref{k1}--\eqref{lk}
is to first derive \eqref{amplify}--\eqref{delta1}
and then bound the resulting {\em rank one} matrix accurately.
Huang and Jia~\cite{huangjia} simply bounded
$$
\|\Delta_k\|\leq \|\Delta_k\|_F \leq \|D_2\| \|T_{k2}T_{k1}^{-1}\|_F \|D_1^{-1}\|
$$
with $\|\cdot\|_F$ the F-norm of a matrix, which led to a too pessimistic
overestimate
$$
\|\Delta_k\|\leq \frac{\sigma_{k+1}}{\sigma_k}\frac{\max_{k+1\leq i\leq n}| u_i^T b|}
{\min_{1\leq i\leq k}|u_i^T b|}
\sqrt{k(n-k)}|L_{k_1}^{(k)}(0)|
$$
due to the excessive factor $\sqrt{k(n-k)}$, which ranges from $\sqrt{n-1}$
to $\frac{n}{2}$ for $n$ even and $\frac{\sqrt{n^2-1}}{2}$ for $n$ odd.

$\|\Delta_k\|$ and $|L_j^{(k)}(0)|,\ j=1,2,\ldots,k$
are used to study the regularizing effects of LSQR in \cite[p.150-2]{hansen98},
but there have been no estimates on them for any kind of ill-posed problem.
We next give accurate estimates for $|L_j^{(k)}(0)|,\ j=1,2,\ldots,k$
and get insight into them for severely ill-posed problems.

\begin{theorem}\label{estlk}
For the severely ill-posed problem with the singular values
$\sigma_j=\mathcal{O}(\rho^{-j})$ and
suitable $\rho>1$, $j=1,2,\ldots,n$ and $k=2,3,\ldots,n-1$, we have
\begin{align}
|L_k^{(k)}(0)|&=1+\mathcal{O}(\rho^{-2}), \label{lkkest}\\
|L_j^{(k)}(0)|&=\frac{1+\mathcal{O}(\rho^{-2})}
{\prod\limits_{i=j+1}^k\left(\frac{\sigma_{j}}{\sigma_i}\right)^2}
=\frac{1+\mathcal{O}(\rho^{-2})}{\mathcal{O}(\rho^{(k-j)(k-j+1)})},
\ j=1,2,\ldots,k-1, \label{lj0}\\
|L_{k_1}^{(k)}(0)|&=\max_{j=1,2,\ldots,k}|L_j^{(k)}(0)|
=1+\mathcal{O}(\rho^{-2}). \label{lkk}
\end{align}
\end{theorem}

{\em Proof.}
Exploiting the Taylor series expansion and
$\sigma_i=\mathcal{O}(\rho^{-i})$ with suitable $\rho>1$, $i=1,2,\ldots,n$,
by definition, for $j=1,2,\ldots,k-1$ we have
\begin{align}
|L_j^{(k)}(0)|&=\prod\limits_{i=1,i\neq j}^k
\left|\frac{\sigma_i^2}{\sigma_i^2-\sigma_j^2}\right|
  =\prod\limits_{i=1}^{j-1}\frac{\sigma_i^2}{\sigma_i^2-\sigma_j^2}
   \cdot\prod\limits_{i=j+1}^{k}\frac{\sigma_i^2}{\sigma_j^2-\sigma_{i}^2}
   \notag\\
& =\prod\limits_{i=1}^{j-1}\frac{1}
{1-\mathcal{O}(\rho^{-2(j-i)})}
\prod\limits_{i=j+1}^{k}\frac{1}
{1-\mathcal{O}(\rho^{-2(i-j)})}\frac{1}
{\prod\limits_{i=j+1}^{k}\mathcal{O}(\rho^{2(i-j)})} \notag\\
&=\frac{\left(1+\sum\limits_{i=1}^j \mathcal{O}(\rho^{-2i})\right)
\left(1+\sum\limits_{i=1}^{k-j+1} \mathcal{O}(\rho^{-2i})\right)}
{\prod\limits_{i=j+1}^{k}\mathcal{O}(\rho^{2(i-j)})} \label{lik}
\end{align}
by absorbing the higher order small terms into $\mathcal{O}(\cdot)$
in the numerator. For $j=k$, we obtain
\begin{align*}
|L_k^{(k)}(0)|&=\prod\limits_{i=1}^{k-1}
\left|\frac{\sigma_i^2}{\sigma_i^2-\sigma_{k}^2}\right|
=\prod\limits_{i=1}^{k-1}\frac{1}
{1-\mathcal{O}(\rho^{-2(k-i)})}=
\prod\limits_{i=1}^{k-1}\frac{1}
{1-\mathcal{O}(\rho^{-2i})}\\
&=1+\sum\limits_{i=1}^k \mathcal{O}(\rho^{-2i})
=1+\mathcal{O}\left(\sum\limits_{i=1}^k\rho^{-2i}\right)\\
&=1+ \mathcal{O}\left(\frac{\rho^{-2}}
    {1-\rho^{-2}}(1-\rho^{-2k})\right)
=1+\mathcal{O}(\rho^{-2}),
\end{align*}
which proves \eqref{lkkest}.

For the numerator of \eqref{lik} we have
  $$
  1+\sum\limits_{i=1}^j \mathcal{O}(\rho^{-2i})
    =1+ \mathcal{O}\left(\sum\limits_{i=1}^j\rho^{-2i}\right)
    =1+ \mathcal{O}\left(\frac{\rho^{-2}}
    {1-\rho^{-2}}(1-\rho^{-2j})\right),
  $$
  and
  $$
    1+\sum\limits_{i=1}^{k-j+1} \mathcal{O}(\rho^{-2i})
    =1+ \mathcal{O}\left(\sum\limits_{i=1}^{k-j+1}\rho^{-2i}\right)
    =1+ \mathcal{O}\left(\frac{\rho^{-2}}{1-\rho^{-2}}
    (1-\rho^{-2(k-j+1)})\right),
  $$
whose product for any $k$ is
  $$
  1+ \mathcal{O}\left(\frac{2\rho^{-2}}{1-\rho^{-2}}\right)
  +\mathcal{O}\left(\left(\frac{\rho^{-2}}{1-\rho^{-2}}\right)^2\right)=
  1+ \mathcal{O}\left(\frac{2\rho^{-2}}{1-\rho^{-2}}\right)
  =  1+\mathcal{O}(\rho^{-2}).
  $$
On the other hand, note that the denominator of \eqref{lik} is defined by
$$
\prod\limits_{i=j+1}^k\left(\frac{\sigma_{j}}{\sigma_i}\right)^2
=\prod\limits_{i=j+1}^{k}\mathcal{O}(\rho^{2(i-j)})
=\mathcal{O}((\rho\cdot\rho^2\cdots\rho^{k-j})^2)
=\mathcal{O}(\rho^{(k-j)(k-j+1)}),
$$
which, together with the above estimate
for the numerator of \eqref{lik}, proves \eqref{lj0}.
Since the above quantity
is always {\em bigger than one} for $j=1,2,\ldots,k-1$,
for any $k$, combining \eqref{lkkest} with \eqref{lj0} gives \eqref{lkk}.
\qquad\endproof

\begin{remark}\label{severerem}
\eqref{lkk} indicates that the bounds \eqref{k1} and \eqref{eqres1}
can be unified as
\begin{align}
\|\Delta_k\|&\leq\frac{\sigma_{k+1}}{\sigma_k}\frac{\max_{k+1\leq i\leq n}| u_i^Tb|}
{\min_{1\leq i\leq k}|u_i^Tb|}
\left(1+\mathcal{O}(\rho^{-2})\right),\ k=1,2,\ldots, n-1.
\label{case3}
\end{align}
\end{remark}

\begin{remark}\label{feature}
(i) \eqref{lj0} shows that $|L_j^{(k)}(0)|$ exhibits monotonic increasing
property with $j$ for a fixed $k$, and $k_1$ in \eqref{lkk}
is close to $k$; (ii) $|L_j^{(k)}(0)|$ decays monotonically with $k$
for a fixed $j$; (iii) \eqref{lkk} indicates
$|L_{k_1}^{(k)}(0)|$ almost remains a constant close to one
with $k$ for suitable $\rho>1$.
\end{remark}

By taking the equalities in \eqref{k1} and \eqref{eqres1} as
estimates for $\|\Delta_k\|$, we substitute them into \eqref{deltabound}
and compute the corresponding estimates for
$\|\sin\Theta(\mathcal{V}_k,\mathcal{V}_k^R)\|$.
We next illustrate the sharpness of our estimates
and justify Remark~\ref{feature}.

Before proceeding, we make some necessary comments.
By the discrete Picard condition \eqref{picard}, \eqref{picard1} and the
properties of $e$, it is known from \cite[p.70-1]{hansen98}
and \cite[p.41-2]{hansen10} that
$
| u_j^T b|\approx| u_j^T b_{true}|=\sigma_j^{1+\beta}>\eta
$
monotonically decreases with $j=1,2,\ldots,k_0$,
and $| u_j^T b|\approx | u_j^T e |$ with the expected values
$
\mathcal{E}(| u_j^T e |)=\eta
$
for $j>k_0$. Therefore, we must have
\begin{align}
\frac{\max_{k+1\leq i\leq n}|u_i^Tb|}{\min_{1\leq i\leq k} |u_i^Tb|}
&\approx \frac{|u_{k+1}^Tb|}{| u_k^Tb|}\approx \frac{\sigma_{k+1}^{1+\beta}}
{\sigma_k^{1+\beta}}<1, \ k=1,2,\ldots,k_0,\label{replace}\\
\frac{\max_{k+1\leq i\leq n}|u_i^Tb|}{\min_{1\leq i\leq k} |u_i^Tb|}
&\approx \frac{|u_{k+1}^Tb|}{| u_k^Tb|}\approx \frac{\eta}{\eta}=1, \
k=k_0+1,\ldots,n-1.\label{replace1}
\end{align}
In numerical justifications, we will
use $\frac{|u_{k+1}^Tb|}{| u_k^Tb|}$ to replace the left-hand sides
of \eqref{replace} and \eqref{replace1}.
It is known from \cite{stewart01,stewartsun}
that if the ratio $\sigma_1/\sigma_k=
\mathcal{O}(\frac{1}{\epsilon_{\rm mach}})$ then both $\sigma_k$ and
$(u_k,v_k)$ are generally
computed with no accuracy in finite precision arithmetic,
where $\epsilon_{\rm mach}=2.22\times 10^{-16}$
is the machine precision, since $\sigma_k$ is very close to
its neighbors for ill-posed problems. Thus,
our above treatment is not only reasonable but also avoids the
{\em intrinsic} difficulty to accurately compute the left-hand sides
of \eqref{replace} and \eqref{replace1}
which involve all the left singular vectors, including
those associated with small singular values clustered at zero that are
computed with no accuracy.

In the meantime, the above also tells us that it is
{\em unreliable} to compute $\Delta_k$ defined by \eqref{defdelta}
and $\|\Delta_k\|$ because, in finite precision arithmetic,
we cannot compute $T_{k2}$ in \eqref{tk12p} reliably
due to the high inaccuracy of the computed small singular values
and the possible underflows of $\sigma_i^{2j-2}$ for $i$ big
and $j=1,2,\ldots,k$. Fortunately,
we can use the Matlab built-in function {\sf subspace.m}, which
computes the maximum of the canonical angles $\Theta(\mathcal{V}_k,
\mathcal{V}_k^R)$, to calculate
$\|\sin\Theta(\mathcal{V}_k,\mathcal{V}_k^R)\|$ accurately
when the first $k$ singular triplets $(\sigma_i,u_i,v_i)$ are computed
accurately.

It appears hard to find a 2D real-life severely ill-posed problem for
justifying the sharpness of our estimates in Theorems~\ref{thm2}--\ref{estlk}.
Gazzola, Hansen and Nagy \cite{gazzola18} have
very recently presented a number of 2D test problems, where
the image deblurring problem {\sf PRblurgauss}, the inverse diffusion problem
{\sf PRdiffusion} and the nuclear magnetic resonance (NMR) relaxometry
problem {\sf PRnmr} are severely ill-posed. But the latter two matrices are
only available as a function handle, for which we cannot compute their SVDs.
Setting the parameter {\sf options.BlurLevel='severe'}, we
have computed the SVD of {\sf PRblurgauss} with $m=n=10000$ and
found $\sigma_1/\sigma_{1500}\approx 1.99\times 10^{14}=
\mathcal{O}(\frac{1}{\epsilon_{\rm mach}})$.
Unfortunately, we have found out that at least half of the first 1500
singular values are (genuinely or at least numerically) {\em multiple}.
For example, among the first 40 singular values,
$\sigma_3,\sigma_6, \sigma_8, \sigma_{11}, \sigma_{13}, \sigma_{15}, \sigma_{17},
\sigma_{19}, \sigma_{22}, \sigma_{24}, \sigma_{25}, \sigma_{26}, \sigma_{28},
\sigma_{30}, \sigma_{33}, \sigma_{35},
\sigma_{37}$ and $\sigma_{39}$ are multiple. As a result,
Theorems~\ref{thm2}--\ref{estlk} cannot apply here because of multiple singular
values. In addition, we have found that the average decay rate of the
{\em distinct} ones among the first 1500
singular value is approximately $\rho^{-1}=0.9697$, i.e., $\rho=1.0312$,
which is {\em fairly} close to one and means that the problem is only {\em
slightly} severely ill-posed.

We should point out that, for the purpose of justifying the sharpness
of our estimates in Theorems~\ref{thm2}--\ref{estlk},
it is enough to test any severely ill-posed problem with the discrete
Picard condition satisfied. To this end, we take
the 1D severely ill-posed problem {\sf shaw} of $m=n=10240$
from \cite{hansen07} with $\sigma_k=
\mathcal{O}({\bf e}^{-2k})$, where ${\bf e}$ is the natural constant.
We take $\rho={\bf e}^2$, and compute the estimate \eqref{deltabound} for
$\|\sin\Theta(\mathcal{V}_k,\mathcal{V}_k^R)\|$
by taking the equalities in \eqref{k1} and \eqref{eqres1},
$1+\mathcal{O}(\rho^{-2})=1+2\rho^{-2}$ and
$\left(1+\mathcal{O}(\rho^{-2})\right)|L_{k_1}^{(k)}(0)|=1+3\rho^{-2}$.
For $k>1$, we use \eqref{lk} to
compute $| L_j^{(k)}(0)|, \ j=1,2,\ldots,k$ and $| L_{k_1}^{(k)}(0)|$ so
as to confirm Theorem~\ref{estlk} and Remark~\ref{feature}.
We generate $b=b_{true}+e$
with $\varepsilon=\frac{\|e\|}{\|b_{true}\|}=10^{-3}$ and $e$
the Gaussian white noise with zero mean. In all the figures, we will
abbreviate $\|\sin\Theta(\mathcal{V}_k,\mathcal{V}_k^R)\|$
by $\sin\Theta_k$.

\begin{figure}[!htp]
\begin{minipage}{0.48\linewidth}
  \centerline{\includegraphics[width=6.0cm,height=3.5cm]{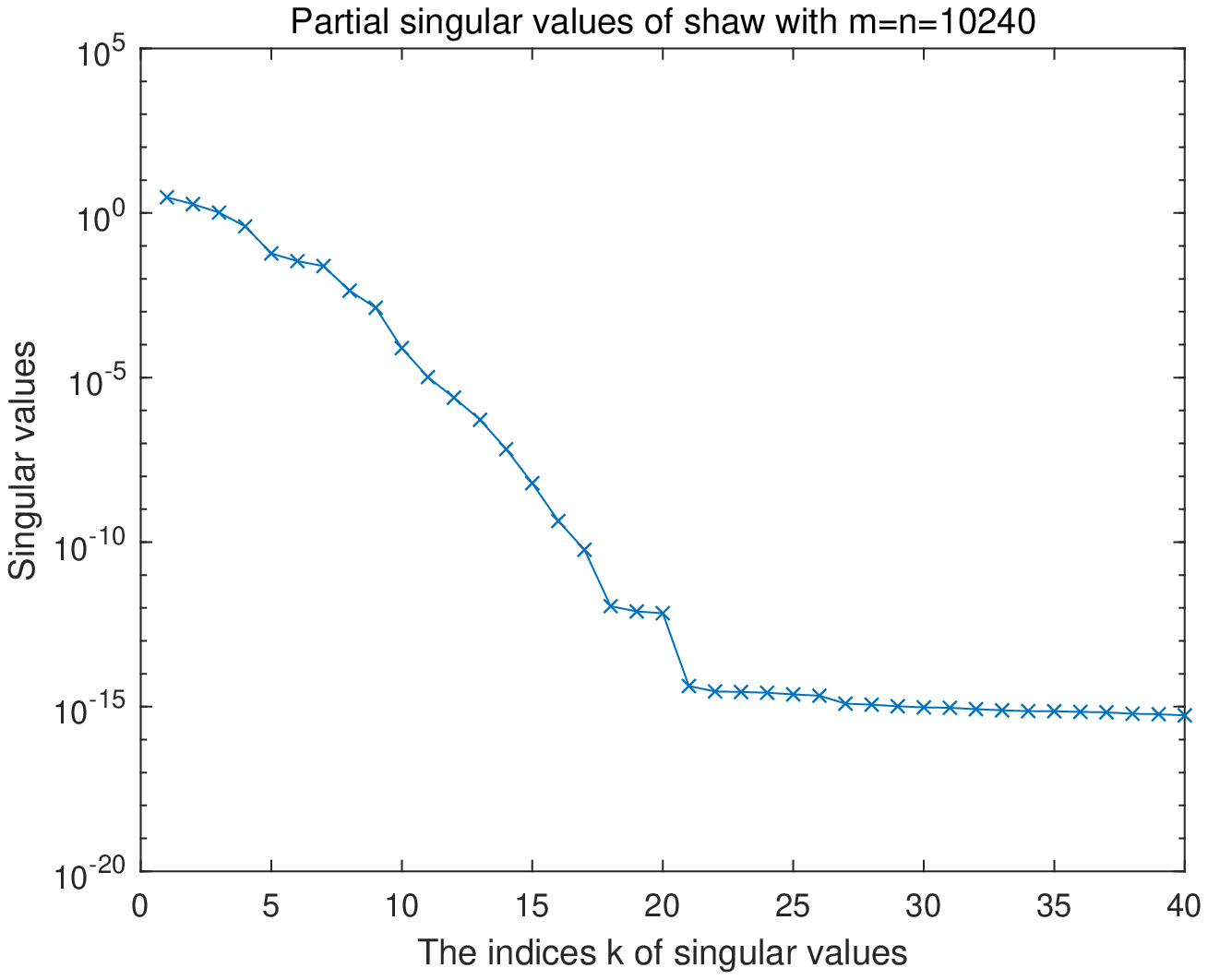}}
  \centerline{(a)}
\end{minipage}
\hfill
\begin{minipage}{0.48\linewidth}
  \centerline{\includegraphics[width=6.0cm,height=3.5cm]{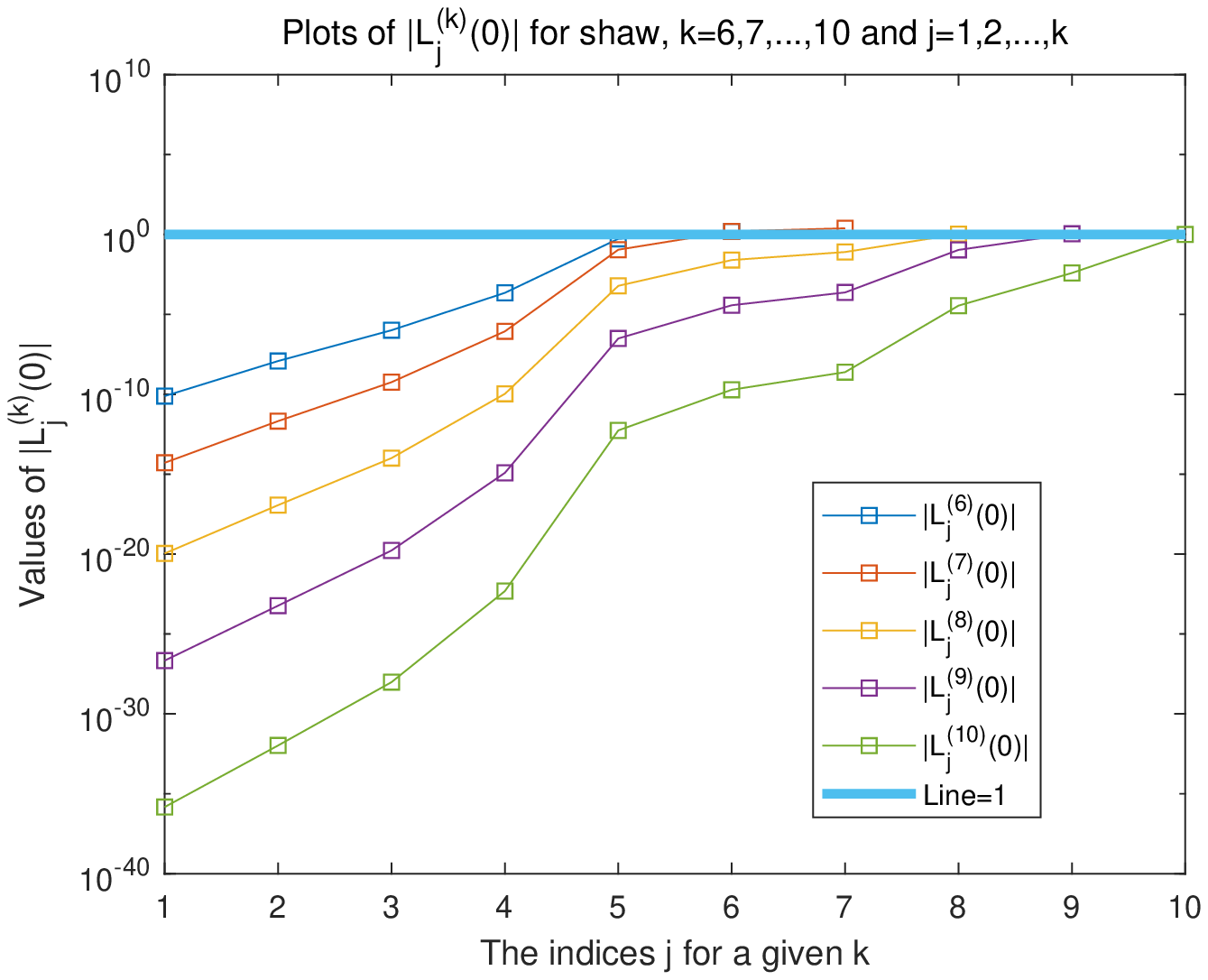}}
  \centerline{(b)}
\end{minipage}
\begin{minipage}{0.48\linewidth}
  \centerline{\includegraphics[width=6.0cm,height=3.5cm]{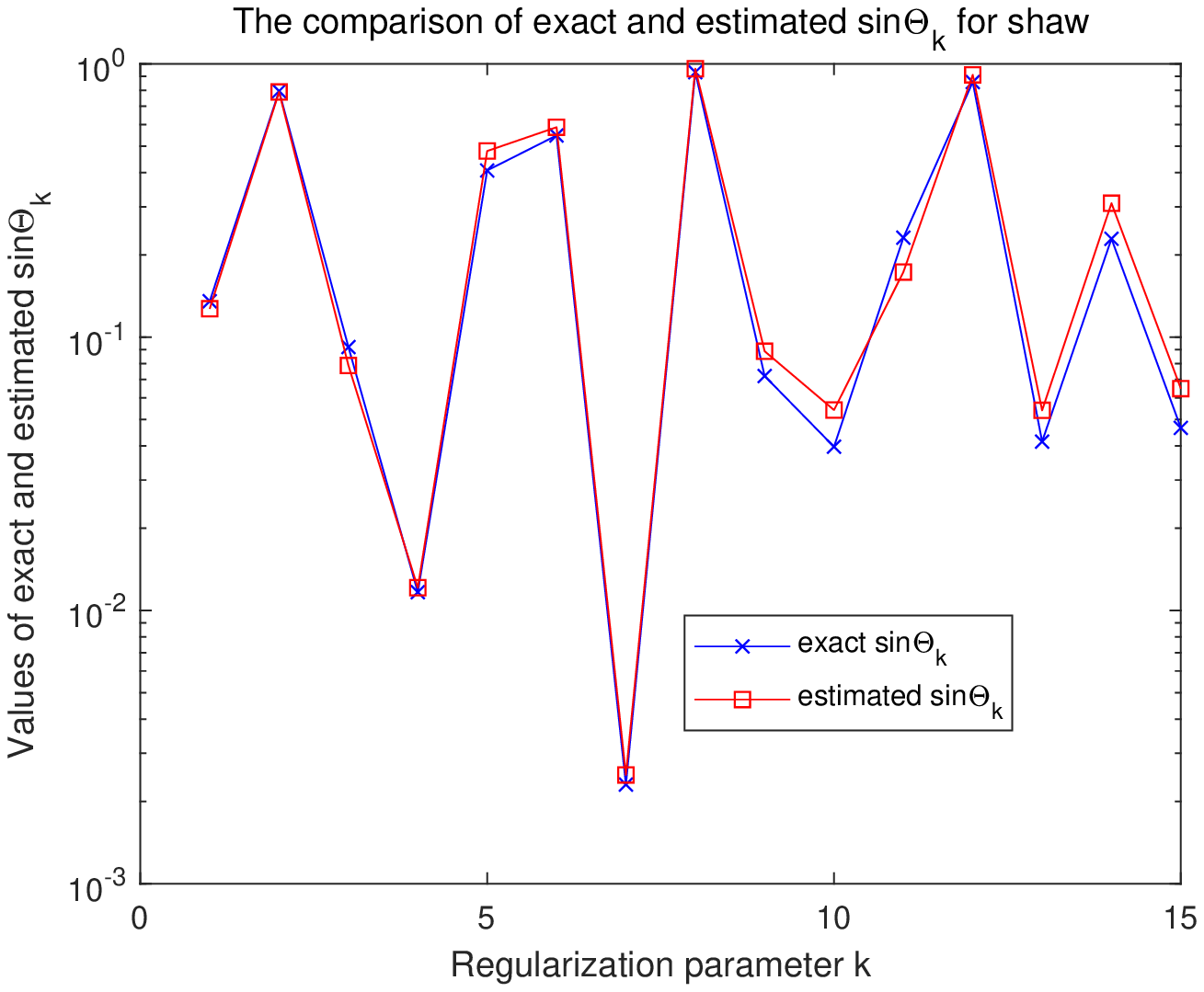}}
  \centerline{(c)}
\end{minipage}
\hfill
\begin{minipage}{0.48\linewidth}
  \centerline{\includegraphics[width=6.0cm,height=3.5cm]{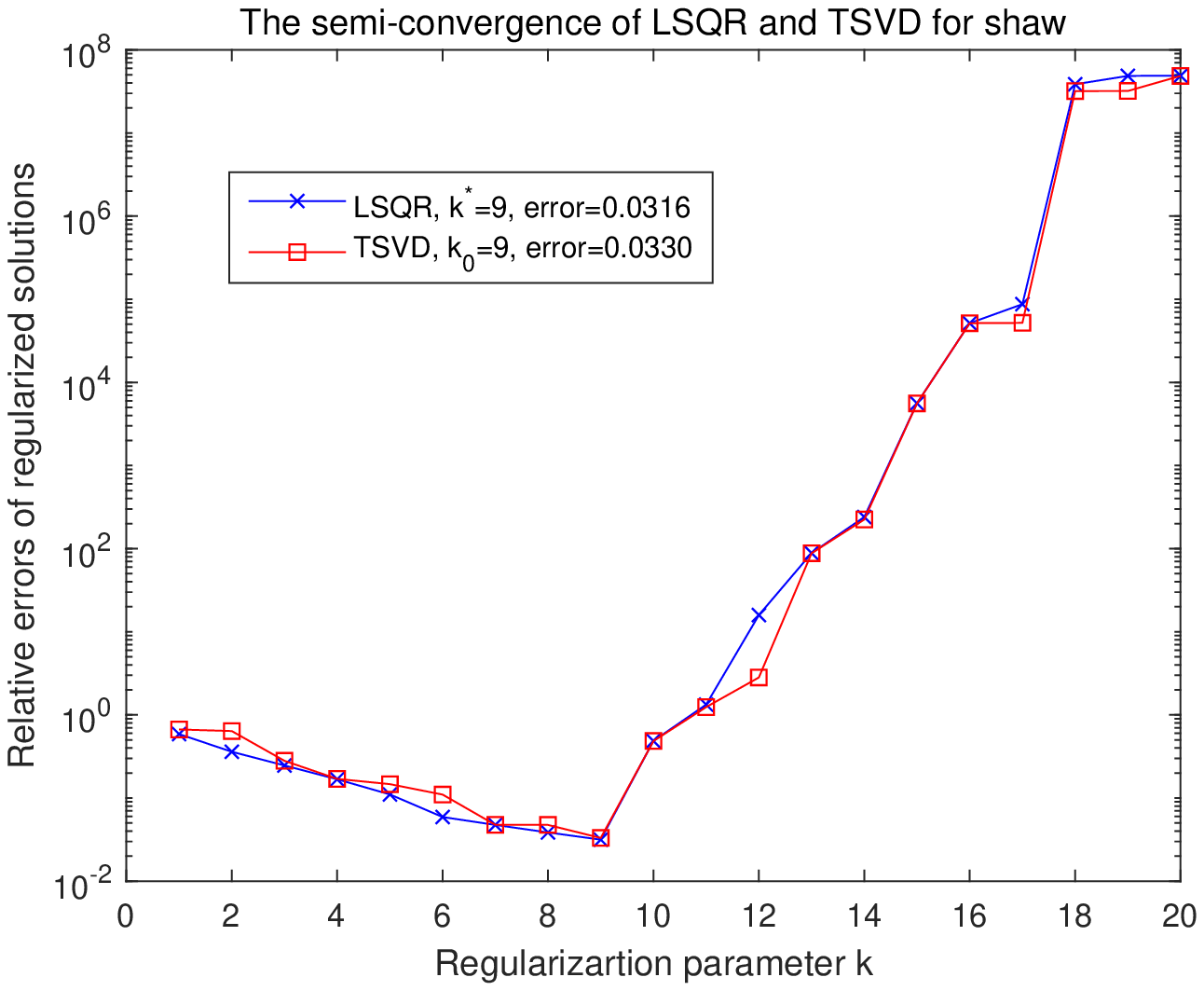}}
  \centerline{(d)}
\end{minipage}
\caption{(a): Partial singular values of {\sf shaw};
(b): plots of $|L_j^{(k)}(0)|$ for $k=6,7,8,9,10$;
(c): the exact and estimated $\|\sin\Theta(\mathcal{V}_k,\mathcal{V}_k^R)\|$;
(d): the semi-convergence process of LSQR and TSVD.}
\label{figshaw1}
\end{figure}

Figure~\ref{figshaw1} (a) plots the first 40 singular values $\sigma_k$
of ${\sf shaw}$, and we find $\sigma_1/\sigma_{21}\approx 2.4\times
10^{15}=\mathcal{O}(\frac{1}{\epsilon_{\rm mach}})$, meaning
that the $(\sigma_k,u_k,v_k)$ are generally computed with no accuracy for $k\geq 21$.
Figure~\ref{figshaw1} (b) clearly confirms each of the three points in
Remark~\ref{feature}.
Moreover, we see that the $|L_j^{(k)}(0)|$ become tiny swiftly for $j$ small
when $k$ increases.
Figure~\ref{figshaw1} (c) indicates that our estimates
for $\|\sin\Theta(\mathcal{V}_k,\mathcal{V}_k^R)\|$ match the exact ones quite
well for $k=1,2,\ldots,15$.
We have found that the maximum and minimum of ratios of the estimated and exact
ones are 1.3924 and 0.7485, respectively,
and the geometric mean of these ratios is 0.9146.
Precisely, the fifteen ratios are
  {\sf
   0.9386,
   0.9924,
   0.8564,
   1.0382,
   1.1781,
   1.0719,
   1.0851,
   1.0302,
   1.2323,
   1.3630,
   0.7485,
   1.0624,
   1.3013,
   1.3488,
   1.3924,
   }
respectively. Figure~\ref{figshaw1} (d) draws the semi-convergence process
of LSQR and the TSVD method. It shows that they
compute the best regularized solutions at the same iterations $k^*=k_0=9$
and the best LSQR solution $x_{k^*}^{lsqr}$ is as accurate
as the best TSVD solution $x_{k_0}^{tsvd}$.


Next we estimate $\|\sin\Theta(\mathcal{V}_k,\mathcal{V}_k^R)\|$ for
moderately and mildly ill-posed problems.

\begin{theorem}\label{moderate}
For a moderately or mildly ill-posed \eqref{eq1} with the singular values
$\sigma_j=\zeta j^{-\alpha},\ j=1,2,\ldots,n$,
where $\alpha>\frac{1}{2}$ and $\zeta>0$ is some constant,
we have
\begin{align}
\|\Delta_1\|&\leq \frac{\max_{2\leq i\leq n}|u_i^Tb|}{| u_1^Tb|}
\sqrt{\frac{1}{2\alpha-1}},\label{mod1}\\
\|\Delta_k\|&\leq \frac{\max_{k+1\leq i\leq n}|u_i^Tb|}{
\min_{1\leq i\leq k}| u_i^T b|}
\sqrt{\frac{k^2}{4\alpha^2-1}+\frac{k}{2\alpha-1}}|L_{k_1}^{(k)}(0)|,
\ k=2,3,\ldots,n-1. \label{modera2}
\end{align}
\end{theorem}

{\em Proof.}
We only need to
accurately bound the right-hand side of \eqref{delta2}. For $k=1,2,\ldots, n-1$,
we obtain
\begin{align}
\left(\sum_{j=k+1}^n\sigma_j^2| u_j^Tb|^2\right)^{1/2}
&= \sigma_{k+1}\max_{k+1\leq i\leq n}| u_i^Tb| \left(\sum_{j=k+1}^n
\frac{\sigma_j^2| u_j^Tb|^2}{\sigma_{k+1}^2\max_{k+1\leq i\leq n}| u_i^Tb|^2}\right)^{1/2}
\notag\\
&\leq \sigma_{k+1}\max_{k+1\leq i\leq n}| u_i^Tb| \left(\sum_{j=k+1}^n
\frac{\sigma_j^2}{\sigma_{k+1}^2}\right)^{1/2} \notag\\
&= \sigma_{k+1}\max_{k+1\leq i\leq n}| u_i^Tb| \left(\sum_{j=k+1}^n \left(\frac{j}{k+1}
\right)^{-2\alpha}\right)^{1/2} \notag \\
&=\sigma_{k+1}\max_{k+1\leq i\leq n}| u_i^Tb|
\left((k+1)^{2\alpha}\sum_{j=k+1}^n \frac{1}{j^{2\alpha}}\right)^{1/2}
\notag\\
&< \sigma_{k+1}\max_{k+1\leq i\leq n}| u_i^Tb| (k+1)^{\alpha}\left(\int_k^{\infty}
\frac{1}{x^{2\alpha}} dx\right)^{1/2}
\notag \\
&= \sigma_{k+1}\max_{k+1\leq i\leq n}| u_i^Tb|\left(\frac{k+1}{k}\right)^{\alpha}
\sqrt{\frac{k}{2\alpha-1}} \notag\\
&=\sigma_{k+1}\max_{k+1\leq i\leq n}| u_i^Tb|\frac{\sigma_k}{\sigma_{k+1}} \sqrt{\frac{k}
{2\alpha-1}} \notag\\
&=\sigma_k \max_{k+1\leq i\leq n}| u_i^Tb|\sqrt{\frac{k}
{2\alpha-1}}.
\label{modeest}
\end{align}

Since the function $x^{2\alpha}$ with any $\alpha> \frac{1}{2}$
is convex over the interval $[0,1]$, for $k=2,\ldots, n-1$,
we obtain
\begin{align}
\left(\sum_{j=1}^k \frac{1}{\sigma_j^2| u_j^Tb|^2}\right)^{1/2}
&=\frac{1}{\sigma_k \min_{1\leq i\leq k}| u_i^T b|}\left(\sum_{j=1}^k
\frac{\sigma_k^2 \min_{1\leq i\leq k}| u_i^T b|^2}{\sigma_j^2| u_j^Tb|^2}\right)^{1/2}\notag\\
&\leq \frac{1}{\sigma_k \min_{1\leq i\leq k}| u_i^T b|}
\left(\sum_{j=1}^k \frac{\sigma_k^2}{\sigma_j^2}\right)^{1/2}\notag\\
&=\frac{1}{\sigma_k \min_{1\leq i\leq k}| u_i^T b|}\left(\sum_{j=1}^k \left(\frac{j}{k}
\right)^{2\alpha }\right)^{1/2} \notag\\
&=\frac{1}{\sigma_k \min_{1\leq i\leq k}| u_i^T b|}
\left(k\sum_{j=1}^{k} \frac{1}{k}\left(\frac{j-1}{k}
\right)^{2\alpha }+1\right)^{1/2} \notag \\
&< \frac{1}{\sigma_k \min_{1\leq i\leq k}| u_i^T b|} \left(k\int_0^1
x^{2\alpha }dx+1\right)^{1/2} \notag\\
& \leq
\frac{1}{\sigma_k\min_{1\leq i\leq k}| u_i^T b|} \sqrt{\frac{k}{2\alpha+1}+1}.
\label{estimate2}
\end{align}
Substituting the above and \eqref{modeest} into \eqref{delta2} yields
\eqref{modera2}. For $k=1$, \eqref{mod1} follows from \eqref{deltaexp}
and \eqref{modeest}.
\qquad\endproof

\begin{remark}
For a purely technical reason,
we have used the simplifying singular value
model $\sigma_j=\zeta j^{-\alpha}$ to replace the general form
$\sigma_j=\mathcal{O}(j^{-\alpha})$. This simplifying model
can avoid some troublesome derivations and
non-transparent formulations.
\end{remark}

In the following we estimate $|L_j^{(k)}(0)|, j=1,2,\ldots,k$
for moderately and mildly ill-posed problems. As will be seen
from the proof, it turns out impossible to bound
them from above both elegantly and accurately unless $\alpha>1$ sufficiently.

\begin{theorem}\label{estlk2}
For a moderately ill-posed problem with
$\sigma_j=\zeta j^{-\alpha},\ j=1,2,\ldots,n$ and $\alpha>1$,
if $\alpha>1$ suitably, then for $k=2,3,\ldots,n-1$ we have
\begin{align}
|L_j^{(k)}(0) |&\approx\left(1+\frac{j}{2\alpha+1}\right)
\prod_{i=j+1}^k\left(\frac{j}{i}\right)^{2\alpha},\,j=1,2,\ldots,k-1,
\label{lkjmod1}\\
\frac{k}{2\alpha+1}&<|L_{k_1}^{(k)}(0)|\approx 1+\frac{k}{2\alpha+1}
\label{lk1size}
\end{align}
with the lower bound requiring that $k$ satisfies $\frac{2\alpha+1}{k}\leq 1$;
for a mildly ill-posed problem with
$\sigma_j=\zeta j^{-\alpha},\ j=1,2,\ldots,n$ and $\frac{1}{2}<\alpha\leq 1$,
if $k$ satisfies $\frac{2\alpha+1}{k}\leq 1$,
we have
\begin{equation}\label{lk1sizemild}
\frac{k}{2\alpha+1}<|L_{k_1}^{(k)}(0)|.
\end{equation}
\end{theorem}

{\em Proof.}
Exploiting the first order Taylor expansion, we obtain
\begin{align}
|L_k^{(k)}(0)|&=
\prod\limits_{i=1}^{k-1}\frac{\sigma_i^2}{\sigma_i^2-\sigma_k^2}
=\prod\limits_{i=1}^{k-1}\frac{1}
{1-(\frac{i}{k})^{2\alpha}} \notag\\
&\approx 1+\sum\limits_{i=1}^{k-1}\left(\frac{i}{k}\right)^{2\alpha}
=1+k\sum\limits_{i=1}^k\frac{1}{k}\left(\frac{i-1}{k}\right)^{2\alpha}
\notag\\
&\approx1+k\int_0^1 x^{2\alpha}dx=1+\frac{k}{2\alpha+1}. \label{lkkmoderate}
\end{align}

For $j=1,2,\ldots,k-1$, by the definition of $\sigma_i$,
since $\alpha >\frac{1}{2}$, we have
\begin{align*}
|L_j^{(k)}(0)|&=\prod\limits_{i=1,i\neq j}^k
\left|\frac{\sigma_i^2}{\sigma_i^2-\sigma_j^2}\right|
  =\prod\limits_{i=1}^{j-1}\frac{\sigma_i^2}{\sigma_i^2-\sigma_j^2}
   \cdot\prod\limits_{i=j+1}^{k}\frac{\sigma_i^2}{\sigma_j^2-\sigma_{i}^2}
   \\
& =\prod\limits_{i=1}^{j-1}\frac{1}
{1-\left(\frac{i}{j}\right)^{2\alpha}}
\prod\limits_{i=j+1}^{k}\frac{1}
{1-\left(\frac{j}{i}\right)^{2\alpha}}\frac{1}
{\prod\limits_{i=j+1}^{k}\left(\frac{i}{j}\right)^{2\alpha}} \\
&\approx \left(1+\sum\limits_{i=1}^{j-1}\left(\frac{i}{j}\right)^{2\alpha}\right)
\left(1+\sum\limits_{i=j+1}^{k} \left(\frac{j}{i}\right)^{2\alpha}\right)
{\prod\limits_{i=j+1}^{k}\left(\frac{j}{i}\right)^{2\alpha}}\\
&\leq \left(1+\int_0^1 x^{2\alpha} dx\right)\left(1+j^{2\alpha}\int_j^k \frac{1}
{x^{2\alpha}}dx\right){\prod\limits_{i=j+1}^{k}\left(\frac{j}{i}\right)^{2\alpha}}\\
&= \left(1+\frac{j}{2\alpha+1}\right)\left(1+\frac{j-j^{2\alpha}
k^{-2\alpha+1}}{2\alpha-1}\right)
\prod_{i=j+1}^k\left(\frac{j}{i}\right)^{2\alpha}.
\end{align*}
Note that $\prod_{i=j+1}^k\left(\frac{j}{i}\right)^{2\alpha}$ are always smaller
than one for $j=1,2,\ldots,k-1$, and the smaller $j$ is, the smaller it is.
Furthermore, exploiting
$$
\left(\frac{j}{k}\right)^{k-j}<\prod_{i=j+1}^k\frac{j}{i}
<\left(\frac{j}{j+1}\right)^{k-j},
$$
by some elementary manipulation, for suitable $\alpha>1$ we can justify the
estimates
\begin{equation*}
\frac{j-j^{2\alpha}
k^{-2\alpha+1}}{2\alpha-1}
\prod_{i=j+1}^k\left(\frac{j}{i}\right)^{2\alpha}\approx 0, \ j=1,2,\ldots,k-1.
\end{equation*}
As a result, for suitable $\alpha>1$ we have
$$
|L_j^{(k)}(0) |\approx\left(1+\frac{j}{2\alpha+1}\right)
\prod_{i=j+1}^k\left(\frac{j}{i}\right)^{2\alpha},\,j=1,2,\ldots,k-1,
$$
which proves \eqref{lkjmod1}. The right-hand side of \eqref{lk1size}
follows from the monotonic increasing property of the right-hand side
of \eqref{lkjmod1} with respect to $j$.

On the other hand, once $k$ is such that $\frac{2\alpha+1}{k}\leq 1$, we
always have
\begin{align}
 |L_{k_1}^{(k)}(0)|&\geq |L_{k}^{(k)}(0)|=
\prod\limits_{i=1}^{k-1}\frac{\sigma_i^2}{\sigma_i^2-\sigma_{k}^2}
=\prod\limits_{i=1}^{k-1}\frac{1}
{1-(\frac{i}{k})^{2\alpha}} \notag\\
&> 1+\sum\limits_{i=1}^{k-1}\left(\frac{i}{k}\right)^{2\alpha}
>1+k\int_0^{\frac{k-1}{k}} x^{2\alpha}dx \notag\\
&=1+\frac{k\left(\frac{k-1}{k}\right)^{2\alpha+1}}{2\alpha+1}
\approx 1+\frac{k}{2\alpha+1}\left(1-\frac{2\alpha+1}{k}\right)=
\frac{k}{2\alpha+1}, \label{lowerbound}
\end{align}
which yields the lower bound of \eqref{lk1size} and \eqref{lk1sizemild}.
\qquad\endproof

\begin{remark}\label{feature2}
The inaccuracy source of \eqref{lkjmod1} and \eqref{lkkmoderate}
consists in using the summations $1+\Sigma$ to replace the corresponding
products $\Pi$ in the proof. They
are considerable underestimates for $\frac{1}{2}<\alpha\leq 1$ but are accurate
for suitable $\alpha>1$;
the bigger $\alpha$ is, the more accurate the estimates
\eqref{lkjmod1} and \eqref{lk1size}
are. The derivation of \eqref{lkkmoderate} and \eqref{lowerbound}
indicates that $|L_{k_1}^{(k)}(0)|$ is substantially bigger than
$\frac{k}{2\alpha+1}$ and cannot be bounded from above effectively
for $\alpha>\frac{1}{2}$ not big enough.
\end{remark}

\begin{remark} \label{feature3}
\eqref{lkjmod1} shows that the first two points in
Remark~\ref{feature} apply here, and \eqref{lk1size}--\eqref{lk1sizemild}
indicate that $|L_{k_1}^{(k)}(0)|$ has increasing
tendency with respect to $k$.
\end{remark}

\begin{remark}
For severely ill-posed problems,
we have $\frac{\sigma_{k+1}}{\sigma_k}
\sim \rho^{-1}$, $\frac{| u_{k+1}^T b|}{| u_k^T b|}\sim
\rho^{-1-\beta}<1$ for $k\leq k_0$. Therefore, we see
from \eqref{case3}--\eqref{replace1}
that $\|\sin\Theta(\mathcal{V}_k,\mathcal{V}_k^R)\|$
exhibits neither increasing
nor decreasing tendency for $k=1,2,\ldots,k_0$ and $k=k_0+1,\ldots,n-1$,
respectively, as is numerically justified by Figure~\ref{figshaw1} (c).
However, for moderately ill-posed
problems, notice from \eqref{picard}
that $\frac{|u_{k+1}^Tb|}{|u_k^Tb|}\approx
\left(\frac{k}{k+1}\right)^{\alpha(1+\beta)}$ increases slowly;
\eqref{lk1size} indicate that $\sqrt{\frac{k^2}{4\alpha^2-1}+\frac{k}{2\alpha-1}}
|L_{k_1}^{(k)}(0)|$ increases as $k$ grows. As a result, \eqref{deltabound}
amd \eqref{modera2} illustrate that $\|\sin\Theta(\mathcal{V}_k,\mathcal{V}_k^R)\|$
exhibits increasing tendency with $k$, meaning
that $\mathcal{V}_k^R$ cannot capture
$\mathcal{V}_k$ so well as it does for severely ill-posed problems as $k$
increases. In fact, $\|\sin\Theta(\mathcal{V}_k,\mathcal{V}_k^R)\|$
starts to approach one as $k$ increases,
meaning that $\mathcal{V}_k^R$ will contain substantial
information on the right singular vectors corresponding to the $n-k$ small
singular values of $A$.
\end{remark}

\begin{remark}\label{mildrem}
Regarding mildly ill-posed problems,
\eqref{lowerbound} and the comment on it
indicate that $|L_{k_1}^{(k)}(0)|$ is substantially bigger than one
for $\frac{1}{2}<\alpha\leq 1$. Consequently,
the bound \eqref{modera2} thus becomes increasingly large as $k$
increases, causing that $\|\Delta_k\|$ is large and
$\|\sin\Theta(\mathcal{V}_k,\mathcal{V}_k^R)\|\approx 1$ soon.
\end{remark}

\begin{figure}[!htp]
\hfill
\begin{minipage}{0.48\linewidth}
  \centerline{\includegraphics[width=6.0cm,height=3.5cm]{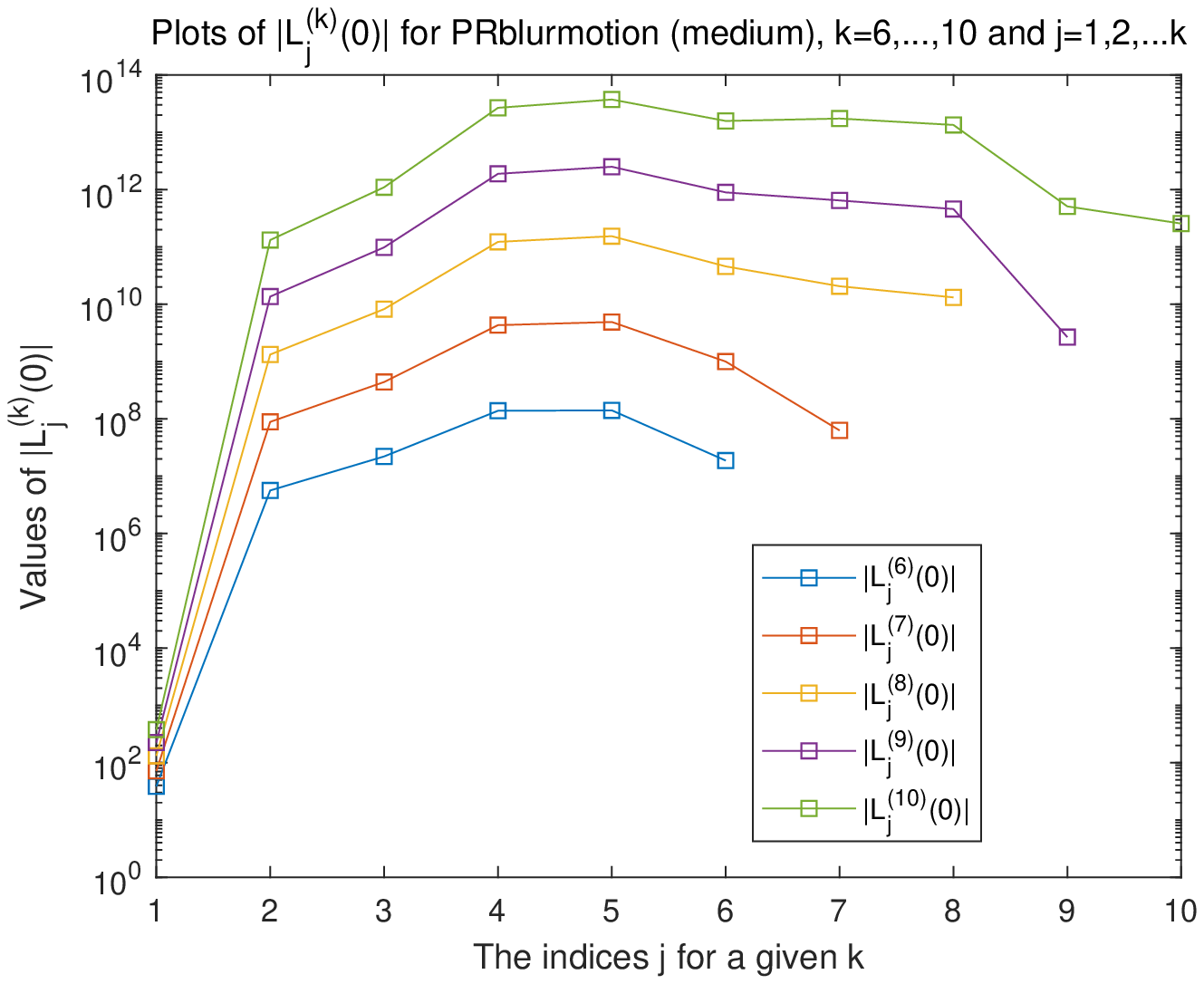}}
\centerline{(a)}
\end{minipage}
\hfill
\begin{minipage}{0.48\linewidth}
  \centerline{\includegraphics[width=6.0cm,height=3.5cm]{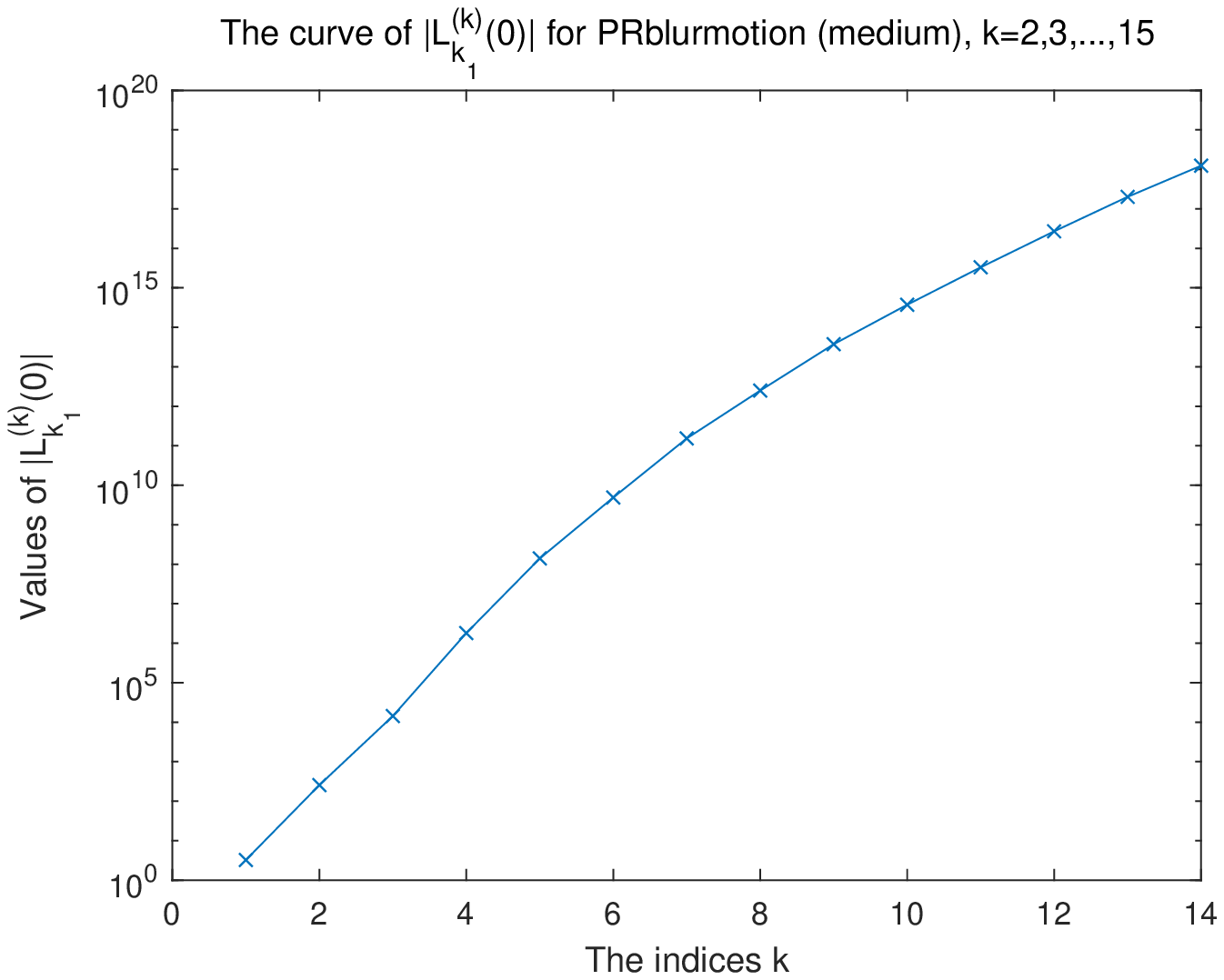}}
\centerline{(b)}
\end{minipage}
\hfill
\begin{minipage}{0.48\linewidth}
  \centerline{\includegraphics[width=6.0cm,height=3.5cm]{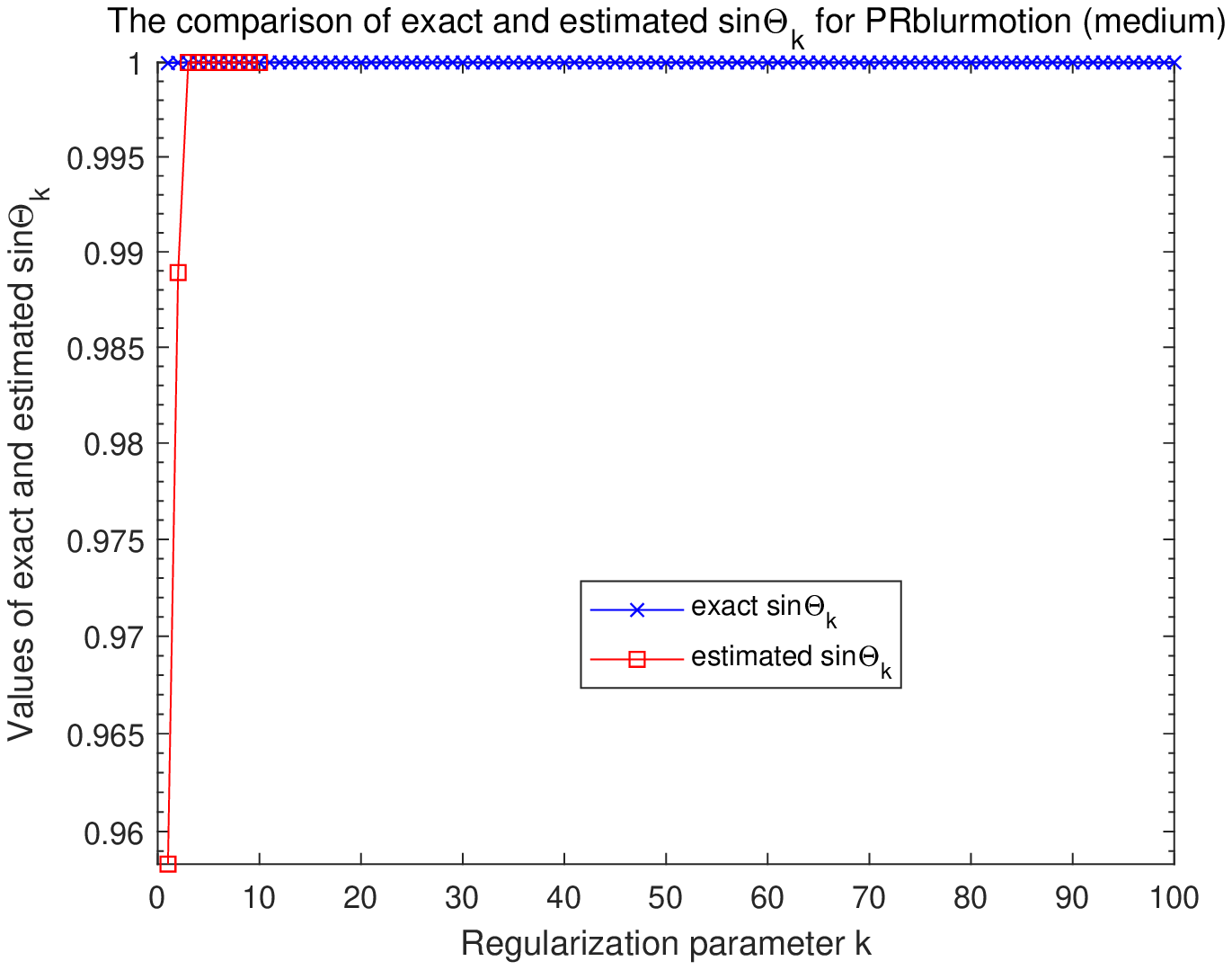}}
\centerline{(c)}
\end{minipage}
\hfill
\begin{minipage}{0.48\linewidth}
  \centerline{\includegraphics[width=6.0cm,height=3.5cm]{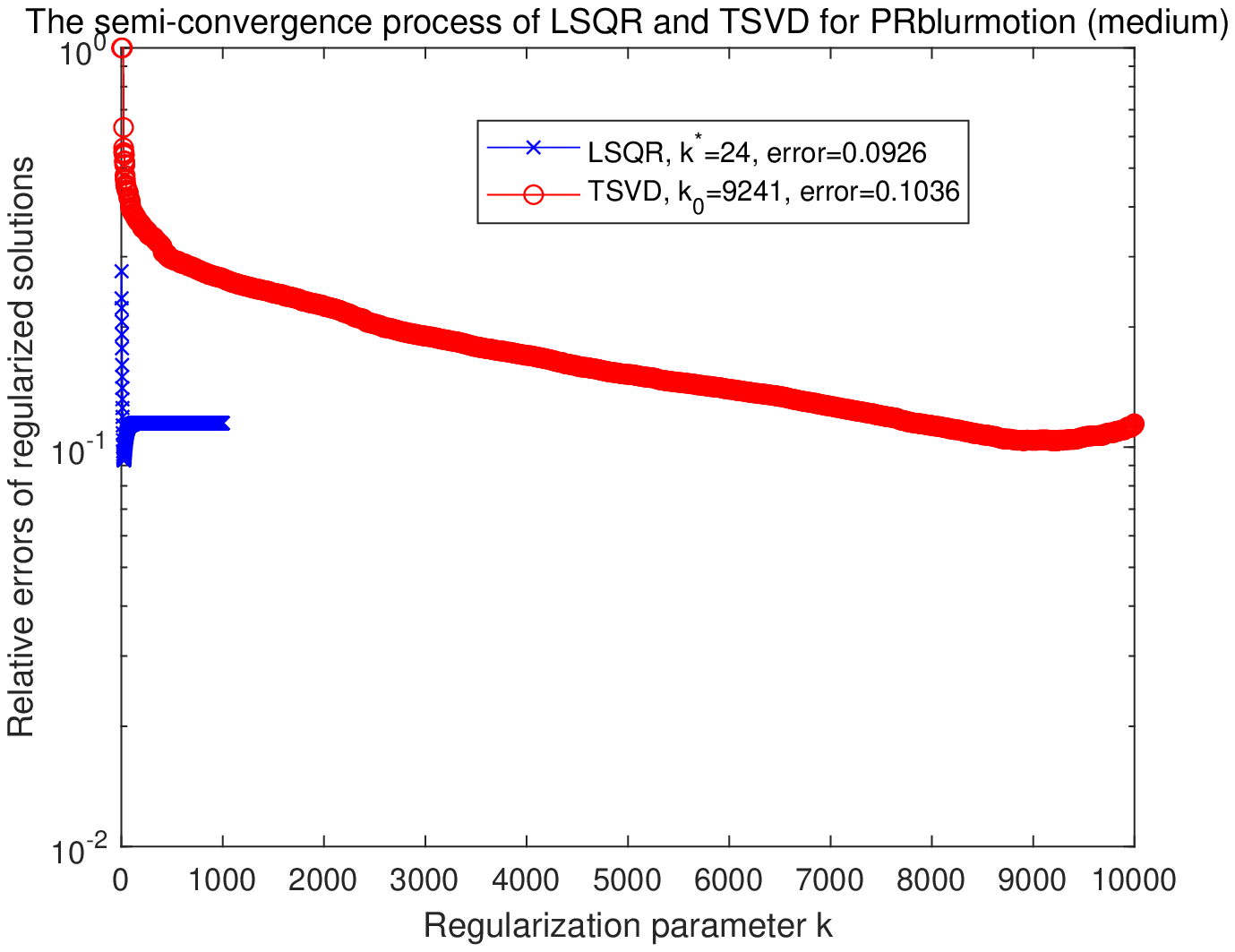}}
\centerline{(d)}
\end{minipage}
\caption{(a): plots of $|L_j^{(k)}(0)|$ for $k=6,7,\ldots,10$;
(b): the curve of $|L_{k_1}^{(k)}(0)|$ for $k=2,3,..,15$;
(c): the exact and estimated $\|\sin\Theta(\mathcal{V}_k,\mathcal{V}_k^R)\|$;
(d): the semi-convergence process of LSQR and TSVD.}
\label{figPRblur}
\end{figure}

In the following we justify our results numerically.
A number of 2D real-life mildly ill-posed problems are presented
in \cite{berisha,gazzola18,hansen07}, which are from image deblurring,
seismic and computerized tomography, inverse diffusion and inverse interpolation,
etc. However, we do not find a 2D moderately ill-posed
one in \cite{berisha,gazzola18}.
{\sf PRblurmotion} from \cite{gazzola18} is a mildly ill-posed
2D image deblurring problem. We use it to show the
effectiveness of our estimates.
Taking {\sf options.BlurLevel='severe', 'medium', 'mild'},
we compute the singular values of three corresponding matrices of $m=n
=10000$, and
find that $\frac{\sigma_1}{\sigma_n}=905.3448,\ 81.1847,\ 28.7967$, respectively.
In the test, except the matrix order, we take all the other parameters as
defaults, by which it is fairly reasonable to use $\alpha=0.6$.
We shall illustrate the sharpness of the estimates for \eqref{deltabound}
when inserting \eqref{mod1}--\eqref{modera2} into it,
where we take the equalities in \eqref{mod1} and \eqref{modera2}. As we have
commented, since the problem is only mildly ill-posed,
we cannot bound $|L_{k_1}^{(k)}(0)|$ from above, and instead
compute it accurately by definition.
We add a Gaussian white noise $e$ to $b$ with the relative noise level
$\varepsilon=0.01$. Figure~\ref{figPRblur} plots the results, where
Figure~\ref{figPRblur} (a)--(b) depict the curves of $|L_j^{(k)}(0)|$
for $k=6,7,8,9,10$ and of $|L_{k_1}^{(k)}(0)|$ for $k=2,3,..,15$, respectively,
Figure~\ref{figPRblur} (c) draws the curves of the exact and estimated
$\|\sin\Theta(\mathcal{V}_k,\mathcal{V}_k^R)\|$, and Figure~\ref{figPRblur} (d)
exhibits the semi-convergence process of LSQR and the TSVD method.

Figure~\ref{figPRblur} (a) justifies that $|L_j^{(k)}(0)|$ increases
with $k$ for a fixed $j$, and and Figure~\ref{figPRblur} (b)
indicates that $|L_{k_1}^{(k)}(0)|$ increases very quickly with $k$.
The ratios of the estimated $\|\sin\Theta(\mathcal{V}_k,\mathcal{V}_k^R)\|$ and
true ones are
{\sf
0.9584,
   0.9889,
   1.0000,
   1.0000,
   1.0000,
   1.0000,
   1.0000,
   1.0000,
   1.0000,
   1.0000,}
respectively, and the geometric mean is 0.9946. In fact,
$\|\sin\Theta(\mathcal{V}_k,\mathcal{V}_k^R)\|$ approaches one from the very
first iteration, and the its first four values are 0.999967738200681,
   0.999986335237377,
   0.999999995155997,
   0.999999999987418,
respectively, which confirm our theory that $\mathcal{V}_k^R$
captures $\mathcal{V}_k$ very poorly and deviates completely from the
latter very soon. In any event, however, our estimates for
$\|\sin\Theta(\mathcal{V}_k,\mathcal{V}_k^R)\|$ match the true ones
very well, as is also seen from Figure~\ref{figPRblur} (c).

Very surprisingly, for this mildly ill-posed problem,
it is completely beyond one's common expectation
that the LSQR best solution $x_{k^*}^{lsqr}$  with the relative error
0.0926 at
semi-convergence is at least
as accurate as the TSVD best solution $x_{k_0}^{tsvd}$  with
the relative error 0.1036;
see Figure~\ref{figPRblur} (d), where LSQR finds its best regularized
solution at iteration $k^*=24$, much more early than the TSVD method,
which computes the best regularized solution at $k_0=9241$, quite close
to $n$. Based on Theorem~\ref{semicon},
this indicates that the $k$ Ritz values $\theta_i^{(k)}$ must not
approximate the large singular values of $A$ in natural order for some
$k\leq k^*$, and we will report numerical results in the next section.
Nevertheless, the results illustrate that LSQR still has the
full regularization. This demonstrates
that the approximations of the $\theta_i^{(k)}$ to the large
$\sigma_i$ in natural order until the occurrence of semi-convergence
of LSQR are {\em not necessary} conditions for the
full regularization of LSQR.
We have also used CGME and LSMR to solve this problem and found that
LSMR has the full regularization too, but CGME computes a considerably
less accurate regularized solution at its semi-convergence and thus
has only the partial regularization; here we omit details on the numerical
results obtained by CGME and LSMR.

In the following we test the 1D moderately ill-posed {\sf heat}
of $n=10240$ from \cite{hansen07},
and illustrate the sharpness of the estimates for \eqref{deltabound}
when inserting \eqref{mod1} and \eqref{modera2} into it,
where we again take the equalities in \eqref{mod1} and \eqref{modera2}.
Regarding the determination of $\alpha$, we compare
the first 1000 singular values
of {\sf heat} with the model singular values $\sigma_1/k^3$,
and we find that the model singular values first decay somewhat faster
than those of {\sf heat} and the rest ones decay more slowly than
those of {\sf heat}. As a result, we take $\alpha=3$ as a rough estimate,
and use it in our estimates.

\begin{figure}[!htp]
\begin{minipage}{0.48\linewidth}
  \centerline{\includegraphics[width=6.0cm,height=3.5cm]{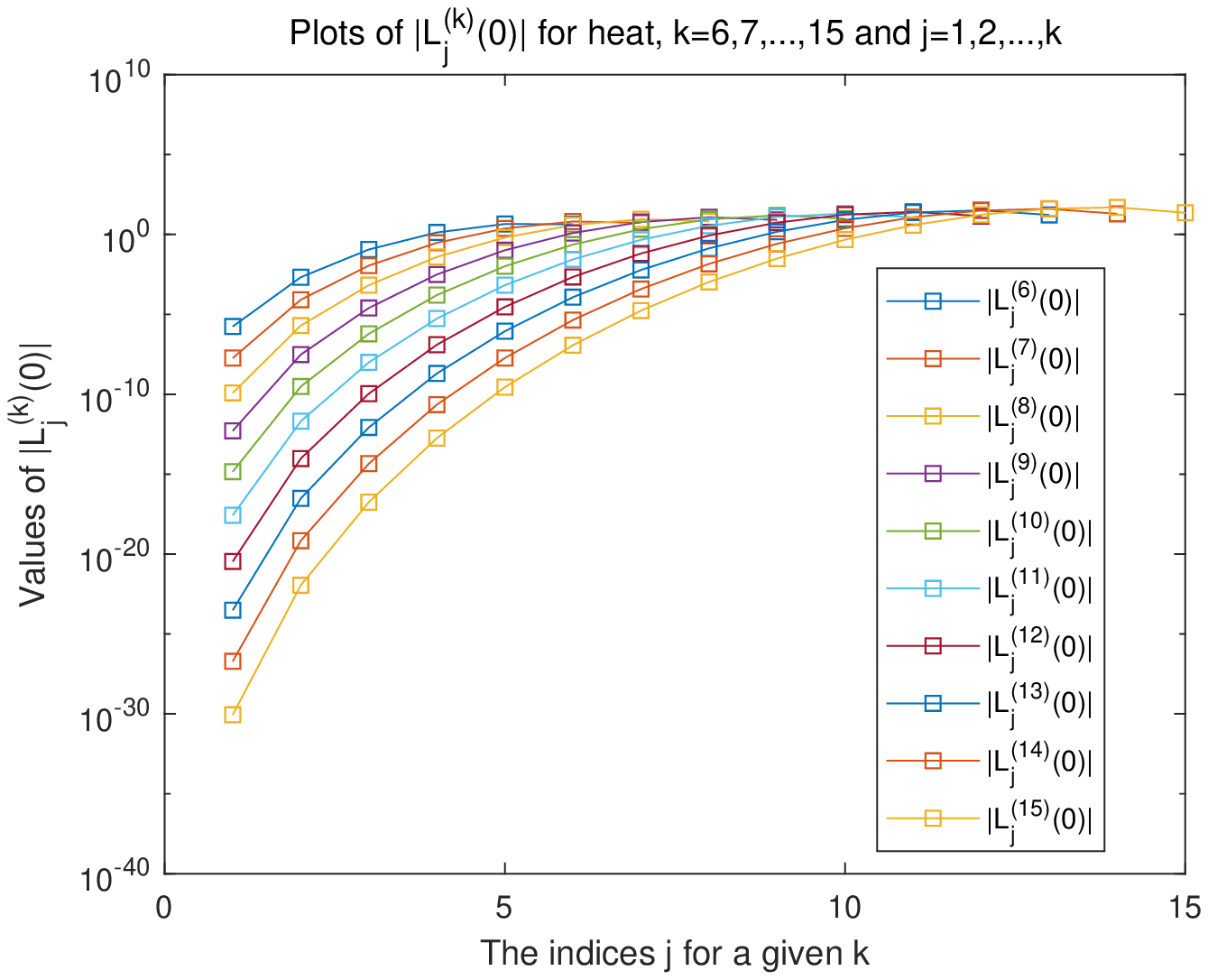}}
\centerline{(a)}
\end{minipage}
\hfill
\begin{minipage}{0.48\linewidth}
  \centerline{\includegraphics[width=6.0cm,height=3.5cm]{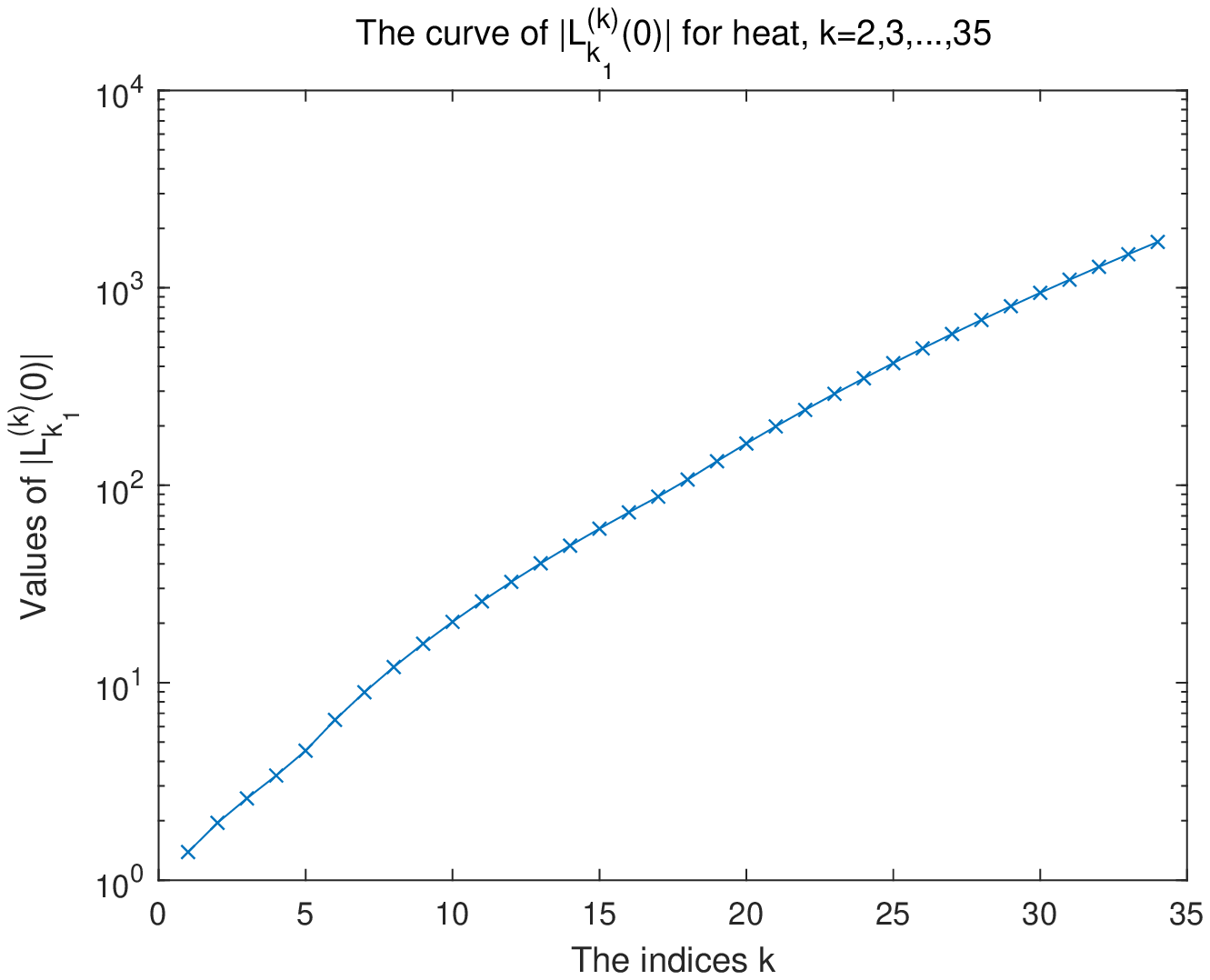}}
\centerline{(b)}
\end{minipage}
\hfill
\begin{minipage}{0.48\linewidth}
  \centerline{\includegraphics[width=6.0cm,height=3.5cm]{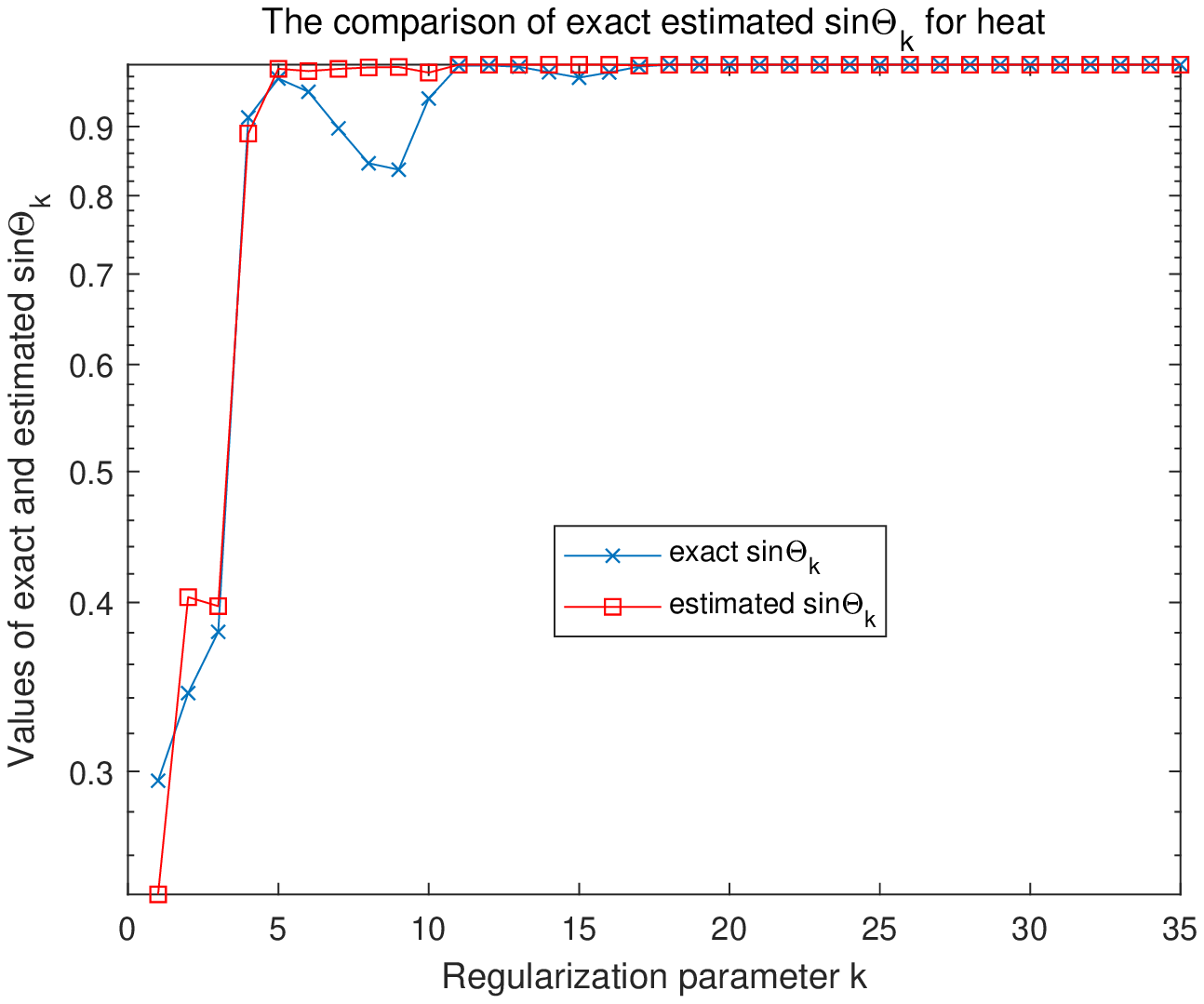}}
\centerline{(c)}
\end{minipage}
\hfill
\begin{minipage}{0.48\linewidth}
  \centerline{\includegraphics[width=6.0cm,height=3.5cm]{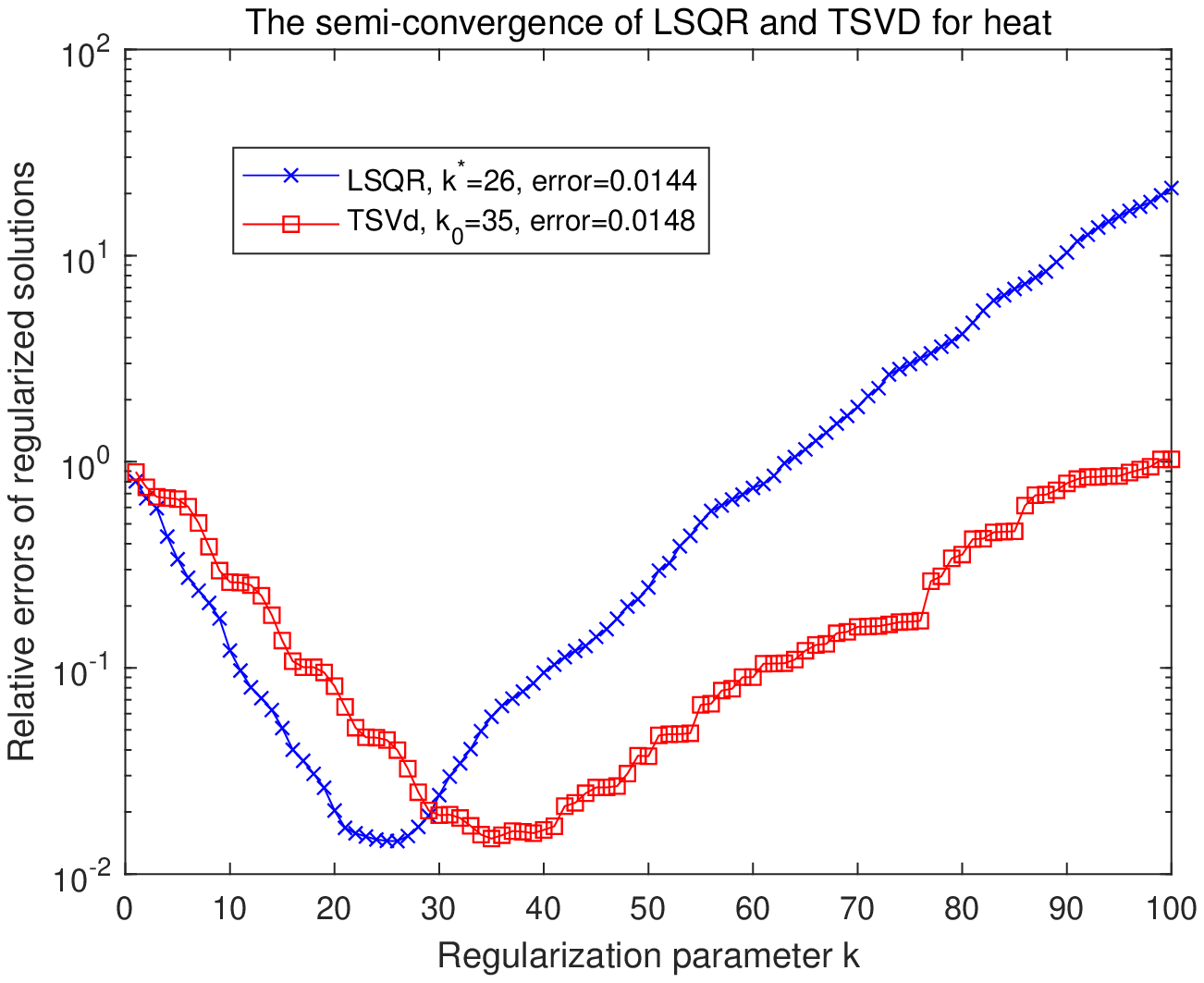}}
\centerline{(d)}
\end{minipage}
\caption{(a): plots of $|L_j^{(k)}(0)|$ for $k=6,7,\ldots,15$;
(b): the curve of $|L_{k_1}^{(k)}(0)|$ for $k=2,3,..,35$;
(c): the exact and estimated $\|\sin\Theta(\mathcal{V}_k,\mathcal{V}_k^R)\|$;
(d): the semi-convergence process of LSQR and TSVD.}
\label{figheat}
\end{figure}

Figure~\ref{figheat} (a) plots the curves of
$| L_j^{(k)}(0)|, \ j=1,2,\ldots,k$ and
$k=6,7,\ldots,15$. It is clear
that $| L_j^{(k)}(0)|$ increases with $k$ for a given $j$ and
exhibits an apparent increasing tendency with $j$ for a given $k$.
Figure~\ref{figheat} (b) shows that, unlike for the severely ill-posed
problem {\sf shaw},  $|L_{k_1}^{(k)}(0)|$ now increases substantially with
$k$ and $\max_{k=2,3,\ldots,35} |L_{k_1}^{(k)}(0)|\approx 1708$,
considerably bigger than one, but it increases much more slowly
than it does for the mildly ill-posed problem {\sf PRblurmotion}.
Figure~\ref{figheat} (c) indicates that
our estimates for $\|\sin\Theta(\mathcal{V}_k,\mathcal{V}_k^R)\|$
match the exact ones quite well for $k=1,2,\ldots,35$.
We have found that the maximum and
minimum of ratios of the estimated and true ones are 1.1911 and 0.8241,
respectively,
and the geometric mean of the ratios is 1.0167. All these results
indicate that our estimates for $\|\sin\Theta(\mathcal{V}_k,\mathcal{V}_k^R)\|$
are sharp. We also plot the semi-convergence
process of LSQR and the TSVD method; see  Figure~\ref{figheat} (d),
where the transition point
$k_0=35$ of the TSVD method but the semi-convergence of LSQR
occurs at $k^*=26$, considerably smaller than $k_0$.
It follows from Theorem~\ref{semicon}
that the $k$ Ritz values $\theta_i^{(k)}$ must not
approximate the large singular values of $A$ in natural order for some
$k\leq k^*$. We will report numerical justifications in the next section.
Remarkably, we see that the best LSQR solution $x_{k^*}^{lsqr}$ is
as accurate as the best TSVD solution $x_{k_0}^{tsvd}$. Again, this indicates
that the approximations of the $\theta_i^{(k)}$ to the large
$\sigma_i$ in natural order until the occurrence of semi-convergence
of LSQR are {\em not necessary} conditions for the
full regularization of LSQR. We have also tested CGME and LSMR,
and found LSMR has the full regularization but CGME has only the partial
regularization; the details are omitted here.

Finally, we pay special attention to $|L_j^{(k)}(0)|$ and give more transparent
numerical supports for Theorem~\ref{estlk2} and Remark~\ref{feature2}.
Precisely, for $k=2,3,\ldots,10$
we compute $|L_j^{(k)}(0)|$ and their rough upper bounds
$1+\frac{k}{2\alpha+1}$ when $\alpha>1$ for the model singular
values $\sigma_i=\frac{1}{i^{\alpha}}$ with $\alpha=0.6, 1, 3$ and 4.
Together with Figure~\ref{figheat} (b) for {\sf heat} and the
observations on them, we have found that (i) the smaller $\alpha$ is,
the bigger $|L_{k_1}^{(k)}(0)|$ is for the same $k$ and (ii)
the bigger $\alpha$ is, the smaller $|L_j^{(k)}(0)|$ is for a
fixed $k$ and the same small $j$. Moreover, we have found that,
for $k=10$, $|L_{k_1}^{(k)}(0)|\approx 3962.7$
for $\alpha=0.6$, $|L_{k_1}^{(k)}(0)|\approx 199.88$ for $\alpha=1$,
$|L_{k_1}^{(k)}(0)|\approx 3.5103$ for $\alpha=3$,
and $|L_{k_1}^{(k)}(0)|\approx 2.2877$ for $\alpha=4$.
Actually, the approximate upper bounds $1+\frac{k}{2\alpha+1}$
are 2.4286 and 2.1111 for $\alpha=3$ and $4$,
respectively. Therefore, $1+\frac{k}{2\alpha+1}$ is indeed
a reasonably good estimate for $|L_{k_1}^{(k)}(0)|$ for suitable $\alpha>1$;
the bigger $\alpha$, the more accurate
$ 1+\frac{k}{2\alpha+1}$ is as an estimate for $|L_{k_1}^{(k)}(0)|$.
On the other hand, if $\alpha$ is small, $1+\frac{k}{2\alpha+1}$
underestimates $|L_{k_1}^{(k)}(0)|$ very considerably.
For $k=2,3,\ldots,10$, we have
also observed that $|L_{k_1}^{(k)}(0)|> \frac{k}{2\alpha+1}$
always holds, confirming the low bounds in
\eqref{lk1size} and \eqref{lk1sizemild}.

\section{The effects of $\sin\Theta$ theorem on
the smallest Ritz value $\theta_k^{(k)}$}\label{manif}

In this section, we investigate how
$\|\sin\Theta(\mathcal{V}_k,\mathcal{V}_k^R)\|$ affects the smallest Ritz
value $\theta_k^{(k)}$. We aim at achieving a manifestation
that (i) we {\em may} have $\theta_k^{(k)}>\sigma_{k+1}$ for suitable
$\|\sin\Theta(\mathcal{V}_k,\mathcal{V}_k^R)\|<1$, and (ii) we {\em must} have
$\theta_k^{(k)}<\sigma_{k+1}$, that is, the $k$ Ritz values
$\theta_i^{(k)}$ do not approximate the large singular values $\sigma_i$ of $A$
in natural order when
$\|\sin\Theta(\mathcal{V}_k,\mathcal{V}_k^R)\|$ is sufficiently close to one.
As it will turn out, the occurrence of (i) or (ii) has different
requirements on the size of
$\|\sin\Theta(\mathcal{V}_k,\mathcal{V}_k^R)\|$ for the three kinds of
ill-posed problems.

\begin{theorem}\label{initial}
Let $\|\sin\Theta(\mathcal{V}_k,\mathcal{V}_k^R)\|^2=1-\varepsilon_k^2$ with
$0< \varepsilon_k< 1$, $k=1,2,\ldots,n-1$, and let the unit-length
$\tilde{q}_k\in\mathcal{V}_k^R$
be the vector that has the smallest angle with $span\{V_k^{\perp}\}$, i.e.,
the closest to $span\{V_k^{\perp}\}$, where $V_k^{\perp}$ is the matrix consisting
of the last $n-k$ columns of $V$ defined by \eqref{eqsvd}. Then it holds that
\begin{equation}\label{rqi}
\varepsilon_k^2\sigma_k^2+
(1-\varepsilon_k^2)\sigma_n^2< \tilde{q}_k^TA^TA\tilde{q}_k<
\varepsilon_k^2\sigma_{k+1}^2+
(1-\varepsilon_k^2)\sigma_1^2.
\end{equation}
If $\varepsilon_k\geq \frac{\sigma_{k+1}}{\sigma_k}$,
then
\begin{equation}
\sqrt{\tilde{q}_k^TA^TA\tilde{q}_k}>\sigma_{k+1};
\label{est1}
\end{equation}
if $\varepsilon_k^2\leq\frac{\delta}
{(\frac{\sigma_1}{\sigma_{k+1}})^2-1}$ for a given arbitrarily small
$\delta>0$, then
\begin{equation}\label{thetasigma}
\theta_k^{(k)}<(1+\delta)^{1/2}\sigma_{k+1},
\end{equation}
meaning that $\theta_k^{(k)}<\sigma_{k+1}$
once $\varepsilon_k$ is sufficiently small, i.e.,
$\|\sin\Theta(\mathcal{V}_k,\mathcal{V}_k^R)\|$ is sufficiently close to
one.
\end{theorem}

{\em Proof.}
Since the columns of $Q_k$ generated by Lanczos bidiagonalization form an
orthonormal basis of $\mathcal{V}_k^R$, by definition and the assumption on
$\tilde{q}_k$ we have
\begin{align}
\|\sin\Theta(\mathcal{V}_k,\mathcal{V}_k^R)\|&=\|(V_k^{\perp})^TQ_k\|
=\|V_k^{\perp}(V_k^{\perp})^TQ_k\| \notag\\
&=\max_{\|c\|=1}\|V_k^{\perp}(V_k^{\perp})^TQ_kc\|
=\|V_k^{\perp}(V_k^{\perp})^T Q_kc_k\| \notag\\
&=\|V_k^{\perp}(V_k^{\perp})^T\tilde{q}_k\|=\|(V_k^{\perp})^T\tilde{q}_k\|
=\sqrt{1-\varepsilon_k^2}
\label{qktilde}
\end{align}
with $\tilde{q}_k=Q_kc_k\in\mathcal{V}_k^R$ and $\|c_k\|=1$.

Expand $\tilde{q}_k$ as the following orthogonal direct sum decomposition:
\begin{equation}\label{decompqk}
\tilde{q}_k=V_k^{\perp}(V_k^{\perp})^T\tilde{q}_k+V_kV_k^T\tilde{q}_k.
\end{equation}
Then from $\|\tilde{q}_k\|=1$ and \eqref{qktilde} we obtain
\begin{align}\label{angle2}
\|V_k^T\tilde{q}_k\|&=\|V_kV_k^T\tilde{q}_k\|=
\sqrt{1-\|V_k^{\perp}(V_k^{\perp})^T\tilde{q}_k\|^2}=\sqrt{1-(1-\varepsilon_k^2)}
=\varepsilon_k.
\end{align}
Keep in mind \eqref{decompqk}. We next bound the Rayleigh quotient of $\tilde{q}_k$
with respect to $A^TA$ from below. By the SVD \eqref{eqsvd} of $A$ and
$V=(V_k,V_k^{\perp})$, we partition
$$
\Sigma=\left(\begin{array}{cc}
\Sigma_k &\\
&\Sigma_k^{\perp}
\end{array}
\right),
$$
where $\Sigma_k={\rm diag}(\sigma_1,\sigma_2,\ldots,\sigma_k)$ and
$\Sigma_k^{\perp}={\rm diag}(\sigma_{k+1},\sigma_{k+2},\ldots,\sigma_n)$.
Making use of $A^TAV_k=V_k\Sigma_k^2$ and $A^TAV_k^{\perp}=
V_k^{\perp}(\Sigma_k^{\perp})^2$ as well as $V_k^TV_k^{\perp}=\mathbf{0}$,
we obtain
\begin{align}
\tilde{q}_k^TA^TA\tilde{q}_k&=\left(V_k^{\perp}(V_k^{\perp})^T\tilde{q}_k+V_kV_k^T
\tilde{q}_k\right)^TA^TA \left(V_k^{\perp}(V_k^{\perp})^T\tilde{q}_k+
V_kV_k^T\tilde{q}_k\right) \notag\\
&=\left(\tilde{q}_k^TV_k^{\perp}(V_k^{\perp})^T+\tilde{q}_k^TV_kV_k^T\right)
\left(V_k^{\perp}(\Sigma_k^{\perp})^2(V_k^{\perp})^T\tilde{q}_k+V_k\Sigma_k^2V_k^T
\tilde{q}_k\right) \notag\\
&=\tilde{q}_k^TV_k^{\perp}(\Sigma_k^{\perp})^2(V_k^{\perp})^T\tilde{q}_k
+\tilde{q}_k^TV_k\Sigma_k^2V_k^T\tilde{q}_k. \label{expansion}
\end{align}

$(V_k^{\perp})^T\tilde{q}_k$ and
$V_k^T\tilde{q}_k$ are unlikely to be the eigenvectors
of $(\Sigma_k^{\perp})^2$
and $\Sigma_k^2$ associated with their respective smallest eigenvalues
$\sigma_n^2$ and $\sigma_k^2$ simultaneously, which are
the $(n-k)$-th canonical vector $e_{n-k}^{(n-k)}$ of $\mathbb{R}^{n-k}$ and
the $k$-th canonical vector $e_k^{(k)}$ of $\mathbb{R}^{k}$, respectively;
otherwise, $\tilde{q}_k=v_n$ and
$\tilde{q}_k=v_k$ simultaneously, which are impossible as $k<n$. Therefore,
from \eqref{expansion}, \eqref{qktilde} and \eqref{angle2},
we obtain the strict inequality
\begin{align*}
\tilde{q}_k^TA^TA\tilde{q}_k&> \|(V_k^{\perp})^T\tilde{q}_k\|^2
\sigma_n^2+\|V_k^T\tilde{q}_k\|^2\sigma_k^2
=(1-\varepsilon_k^2)\sigma_n^2+\varepsilon_k^2 \sigma_k^2,
\end{align*}
from which it follows that the lower bound of \eqref{rqi} holds. By a
similar argument, from \eqref{expansion}  and \eqref{qktilde}, \eqref{angle2}
we obtain the upper bound of \eqref{rqi}:
$$
\tilde{q}_k^TA^TA\tilde{q}_k <\|(V_k^{\perp})^T\tilde{q}_k\|^2
\|(\Sigma_k^{\perp})^2\|+\|V_k^T\tilde{q}_k\|^2\|\Sigma_k^2\|
=(1-\varepsilon_k^2)\sigma_{k+1}^2+\varepsilon_k^2 \sigma_1^2.
$$

From the lower bound of \eqref{rqi}, we see that if
$\varepsilon_k$ satisfies $\varepsilon_k^2 \sigma_k^2\geq \sigma_{k+1}^2$,
i.e., $\varepsilon_k\geq \frac{\sigma_{k+1}}{\sigma_k}$,
then $\sqrt{\tilde{q}_k^TA^TA\tilde{q}_k}>\sigma_{k+1}$, i.e.,
\eqref{est1} holds.

From \eqref{Bk}, we obtain $B_k^TB_k=Q_k^TA^TAQ_k$.
Note that $(\theta_k^{(k)})^2$ is the smallest eigenvalue
of the symmetric positive definite matrix $B_k^TB_k$.
Therefore, we have
\begin{equation}\label{rqi2}
(\theta_k^{(k)})^2=\min_{\|c\|=1} c^TQ_k^TA^TAQ_kc=
\min_{q\in \mathcal{V}_k^R,\ \|q\|=1} q^TA^TAq
=\hat{q}_k^TA^TA\hat{q}_k,
\end{equation}
where $\hat{q}_k$ is, in fact, the Ritz vector of $A^TA$ from
$\mathcal{V}_k^R$ corresponding to the smallest Ritz value $(\theta_k^{(k)})^2$.
Therefore, we have
\begin{equation}\label{thetak}
\theta_k^{(k)}\leq \sqrt{\tilde{q}_k^TA^TA\tilde{q}_k},
\end{equation}
from which and \eqref{rqi} it follows that
$(\theta_k^{(k)})^2<(1-\varepsilon_k^2)\sigma_{k+1}^2+\varepsilon_k^2 \sigma_1^2$.
For any $\delta>0$, we choose $\varepsilon_k\geq 0$ such that
$$
(\theta_k^{(k)})^2<(1-\varepsilon_k^2)\sigma_{k+1}^2+\varepsilon_k^2 \sigma_1^2
\leq (1+\delta)\sigma_{k+1}^2,
$$
i.e., \eqref{thetasigma} holds,
solving which for $\varepsilon_k^2$ gives $\varepsilon_k^2\leq\frac{\delta}
{(\frac{\sigma_1}{\sigma_{k+1}})^2-1}$.
\qquad\endproof

In the sense of \eqref{rqi2}, $\hat{q}_k\in \mathcal{V}_k^R$ is the optimal vector
that extracts the least information from $\mathcal{V}_k$ and the richest
information from $span\{V_k^{\perp}\}$. From the assumption on $\tilde{q}_k$,
since $\mathcal{V}_k$ is the orthogonal complement
of $span\{V_k^{\perp}\}$, we know that $\tilde{q}_k\in \mathcal{V}_k^R$
has the largest acute angle with $\mathcal{V}_k$, that is,
it contains the least information from $\mathcal{V}_k$ and the richest information
from $span\{V_k^{\perp}\}$. Therefore,
$\hat{q}_k$ and $\tilde{q}_k$ have a similar optimality, and consequently
\begin{equation}\label{approxeq}
\theta_k^{(k)}\approx \sqrt{\tilde{q}_k^TA^TA\tilde{q}_k}.
\end{equation}
Combining this estimate with \eqref{est1} and \eqref{thetak}, we {\em may} have
$\theta_k^{(k)}>\sigma_{k+1}$ if
$\varepsilon_k\geq \frac{\sigma_{k+1}}{\sigma_k}$.

We analyze $\theta_k^{(k)}$ and
inspect the condition $\varepsilon_k\geq\frac{\sigma_{k+1}}
{\sigma_k}$ for \eqref{est1}. It is known that
$\varepsilon_k\geq \frac{\sigma_{k+1}}{\sigma_k}\sim \rho^{-1}$
for severely ill-posed problems, meaning that
$\|\sin\Theta(\mathcal{V}_k,\mathcal{V}_k^R)\|$ is approximately
smaller than $1-\frac{1}{2}\rho^{-2}$. For moderately ill-posed
problems, the lower bound $\sigma_{k+1}/\sigma_k$
increases with $k$, and it cannot be close to one for suitable $\alpha>1$;
for mildly ill-posed problems, the lower
bound for $\varepsilon_k$ increases faster than it does for moderately
ill-posed problems since $\alpha\leq 1$,
and, furthermore, it may well approach one for $k$ small. In conclusion,
the condition $\varepsilon_k\geq \frac{\sigma_{k+1}}{\sigma_k}$
for \eqref{est1} requires that
$\|\sin\Theta(\mathcal{V}_k,\mathcal{V}_k^R)\|$
be not close to one for severely and moderately ill-posed problems
with suitable $\alpha>1$,
but it must be fairly small for mildly ill-posed problems.

We now investigate if the true, i.e., actual $\varepsilon_k$ resulting
from the three kinds of ill-posed problems satisfies
 the condition $\varepsilon_k\geq\frac{\sigma_{k+1}}
{\sigma_k}$ for \eqref{est1}. In view of \eqref{deltabound} and
$\|\sin\Theta(\mathcal{V}_k,\mathcal{V}_k^R)\|^2=1-\varepsilon_k^2$,
we have $\|\Delta_k\|^2=\frac{1-\varepsilon_k^2}{\varepsilon_k^2}$.
Thus, the condition $\varepsilon_k\geq\frac{\sigma_{k+1}}
{\sigma_k}$ for \eqref{est1} amounts to requiring
that $\|\Delta_k\|$ cannot
be large for severely and moderately ill-posed problems but
it must be fairly small for mildly ill-posed problems. Unfortunately,
Theorems~\ref{thm2}--\ref{moderate} and the remarks on them
indicate that $\|\Delta_k\|$ is approximately $\rho^{-(2+\beta)}$
by \eqref{case3} and \eqref{replace} for $k\leq k_0$,
considerably smaller than one for a
severely ill-posed problem with $\rho>1$ not close to one,
it is modest and increases slowly with $k$ for a moderately
ill-posed problem with suitable $\alpha>1$, and it
increases with $k$ and is generally
large for a mildly ill-posed problem.
Consequently, for mildly ill-posed problems,
the actual $\|\Delta_k\|$ can hardly be small and is generally
large, namely, the true $\varepsilon_k$ is small, which
causes that the condition $\varepsilon_k\geq\frac{\sigma_{k+1}}{\sigma_k}$
fails to meet soon as $k$ increases, while it is satisfied for severely
or moderately ill-posed problems with suitable $\rho>1$ or $\alpha>1$
for $k$ small.

We report numerical experiments to confirm Theorem~\ref{initial}
and the above remarks. Besides the previous severely, moderately and
mildly problems {\sf shaw}, {\sf heat}, {\sf PRblurmotion},
we also test the 1D moderately ill-posed problem {\sf deriv2} of $n=10000$
\cite{hansen07} and the 2D mildly ill-posed problems
{\sf PRspherical} of $m=14100, n=10000$ and {\sf PRseismic} of
$m=20000, n=10000$, which are from seismic travel-time tomography and
spherical means tomography  \cite{gazzola18}, respectively.
For the latter three problems, except the matrix orders,
we take all the other parameter(s) as default(s). The singular values $\sigma_k$
of {\sf deriv2} decay exactly like $\frac{1}{k^2}$ (cf. \cite[p.21]{hansen10});
the singular values of {\sf PRspherical} and {\sf PRseismic} decay roughly
like $\sigma_1/k^{0.6}$ in an initial
stage, we add Gaussian white noises with the relative level $\varepsilon=0.05$ to
the right-hand sides $b_{true}$ of {\sf PRspherical} and {\sf PRseismic},
respectively. The noise level $\varepsilon=0.001$ is used
in {\sf shaw}, {\sf heat} and {\sf deriv2},
and $\varepsilon=0.01$ is used in {\sf PRblurmotion}, the same as that
in the last section.

\begin{figure}[!htp]
\begin{minipage}{0.48\linewidth}
  \centerline{\includegraphics[width=6.0cm,height=3.5cm]{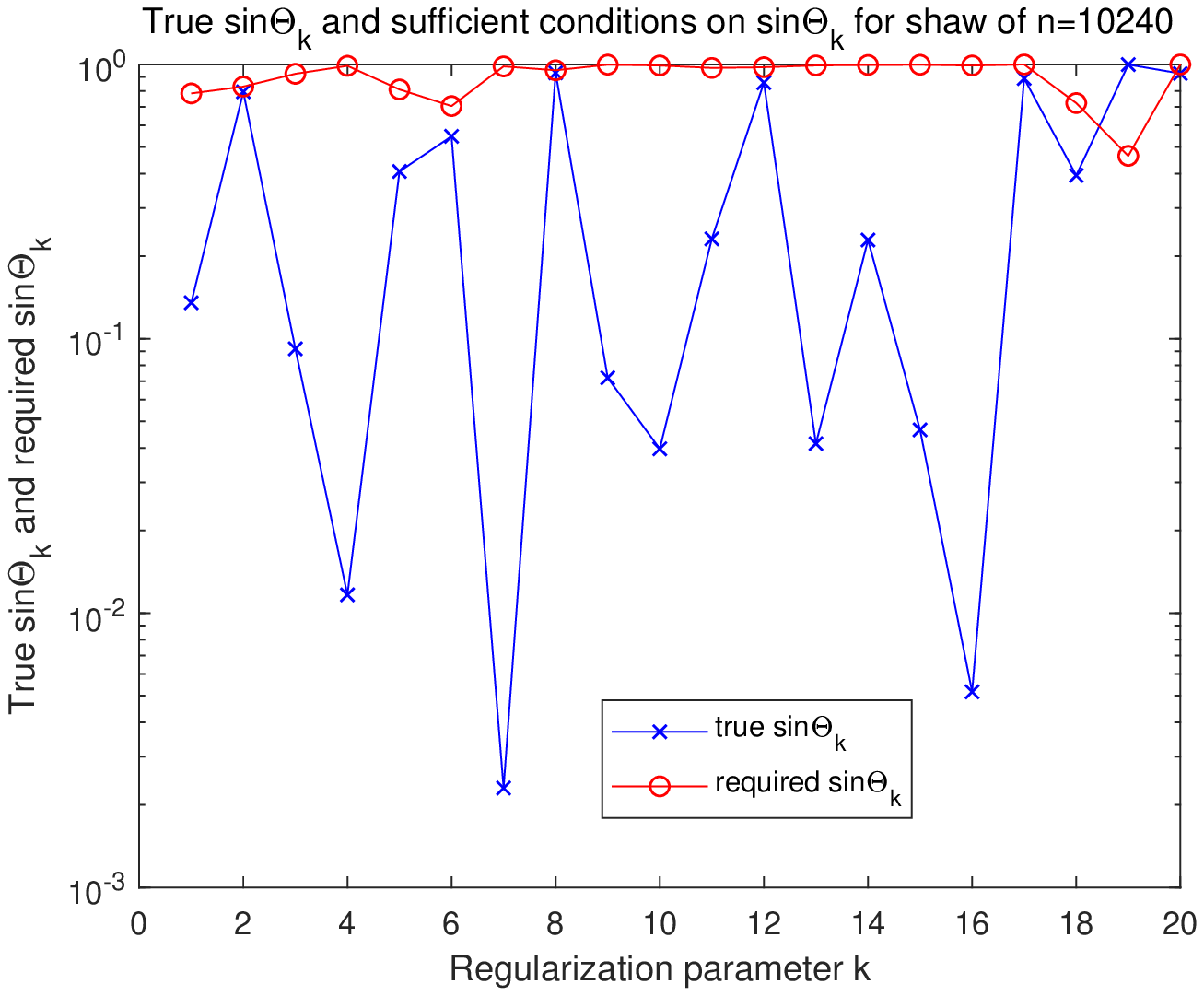}}
\centerline{(a)}
\end{minipage}
\hfill
\begin{minipage}{0.48\linewidth}
  \centerline{\includegraphics[width=6.0cm,height=3.5cm]{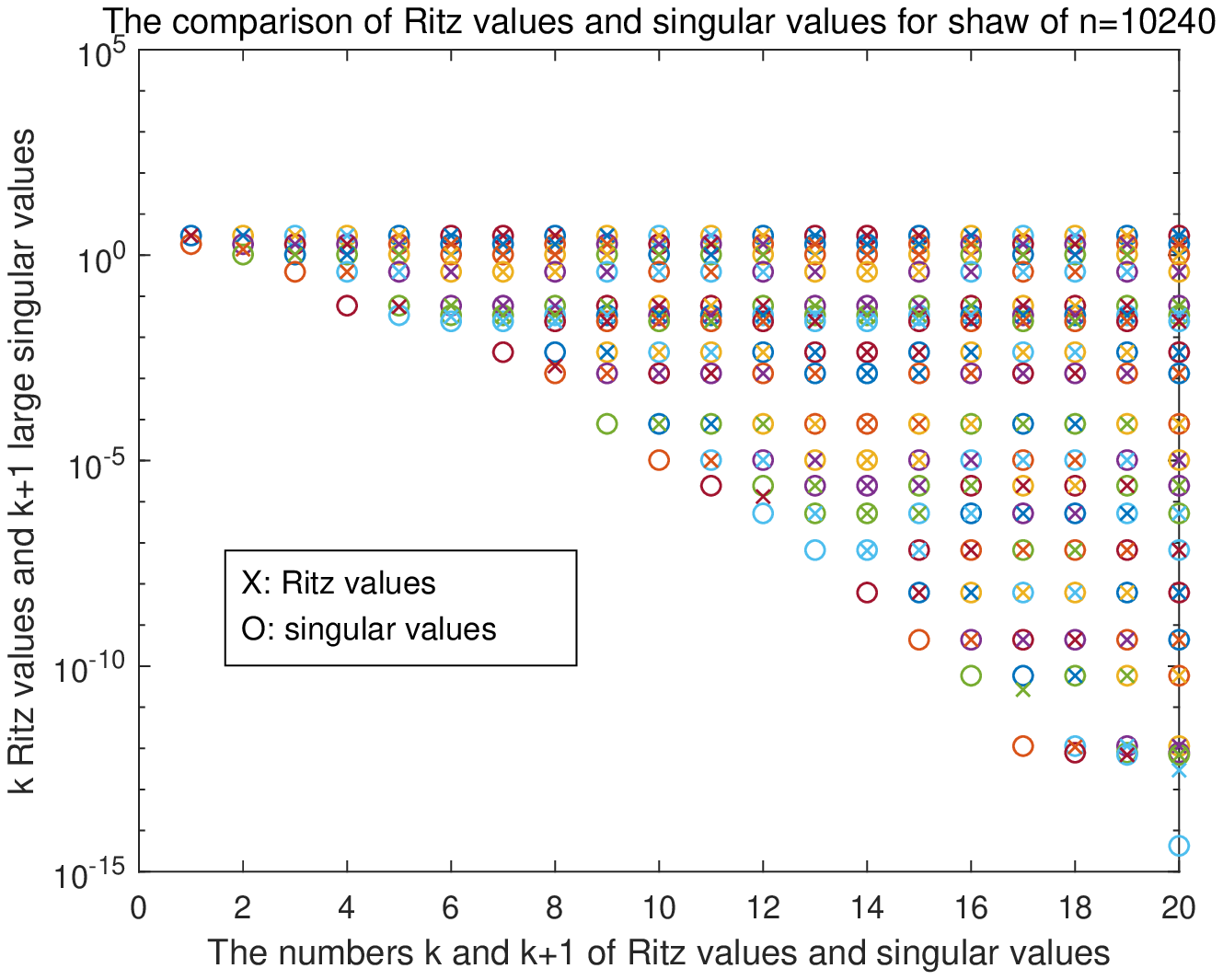}}
\centerline{(b)}
\end{minipage}
\caption{(a): The true $\|\sin\Theta(\mathcal{V}_k,\mathcal{V}_k^R)\|$
and the required sufficient conditions on them; (b): $k$ Ritz values
and the first $k+1$ large singular values
of {\sf shaw}, $k=1,2,\ldots,20$.}
\label{figshaw2}
\end{figure}

\begin{figure}[!htp]
\begin{minipage}{0.48\linewidth}
  \centerline{\includegraphics[width=6.0cm,height=3.5cm]{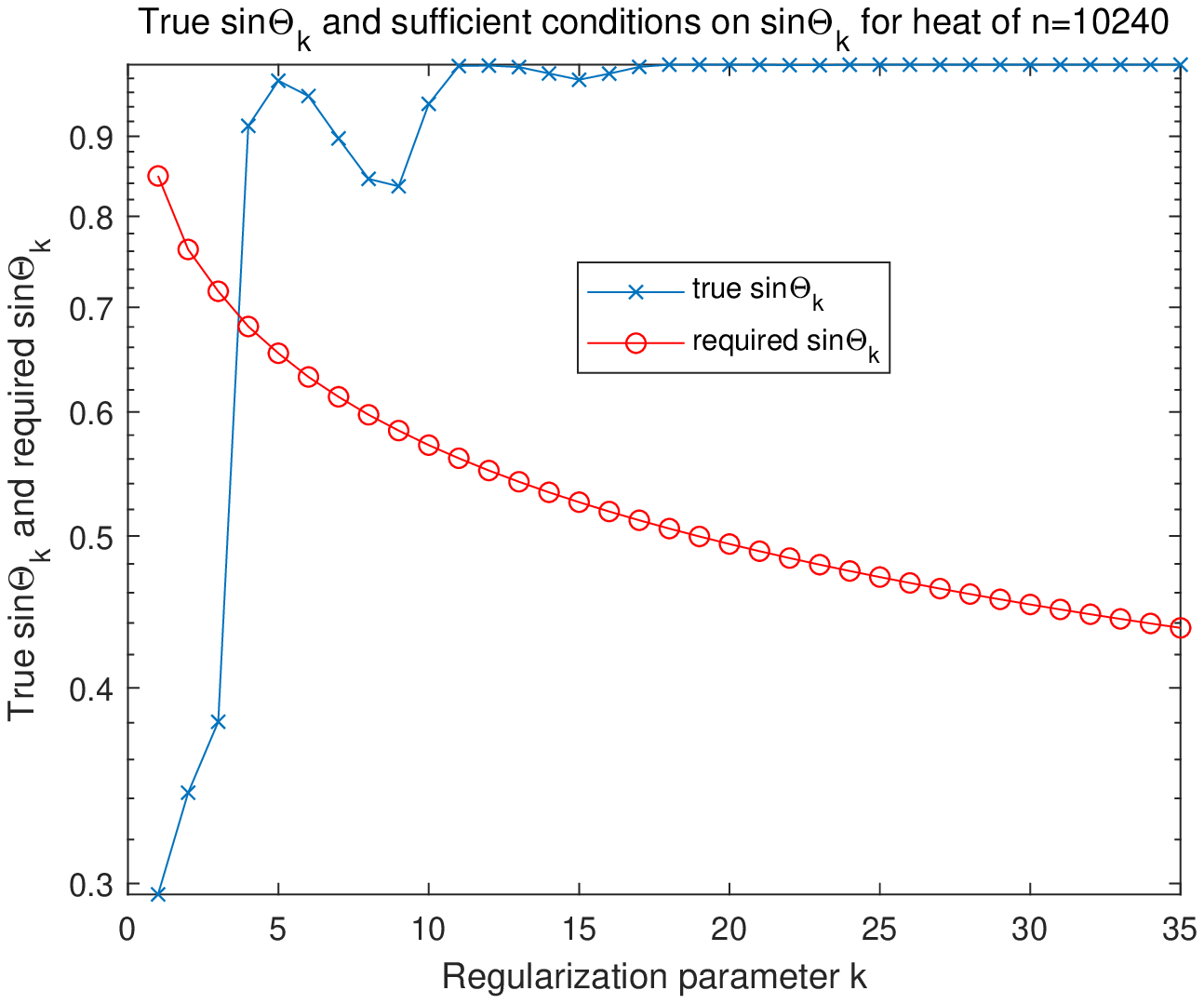}}
\centerline{(a)}
\end{minipage}
\hfill
\begin{minipage}{0.48\linewidth}
  \centerline{\includegraphics[width=6.0cm,height=3.5cm]{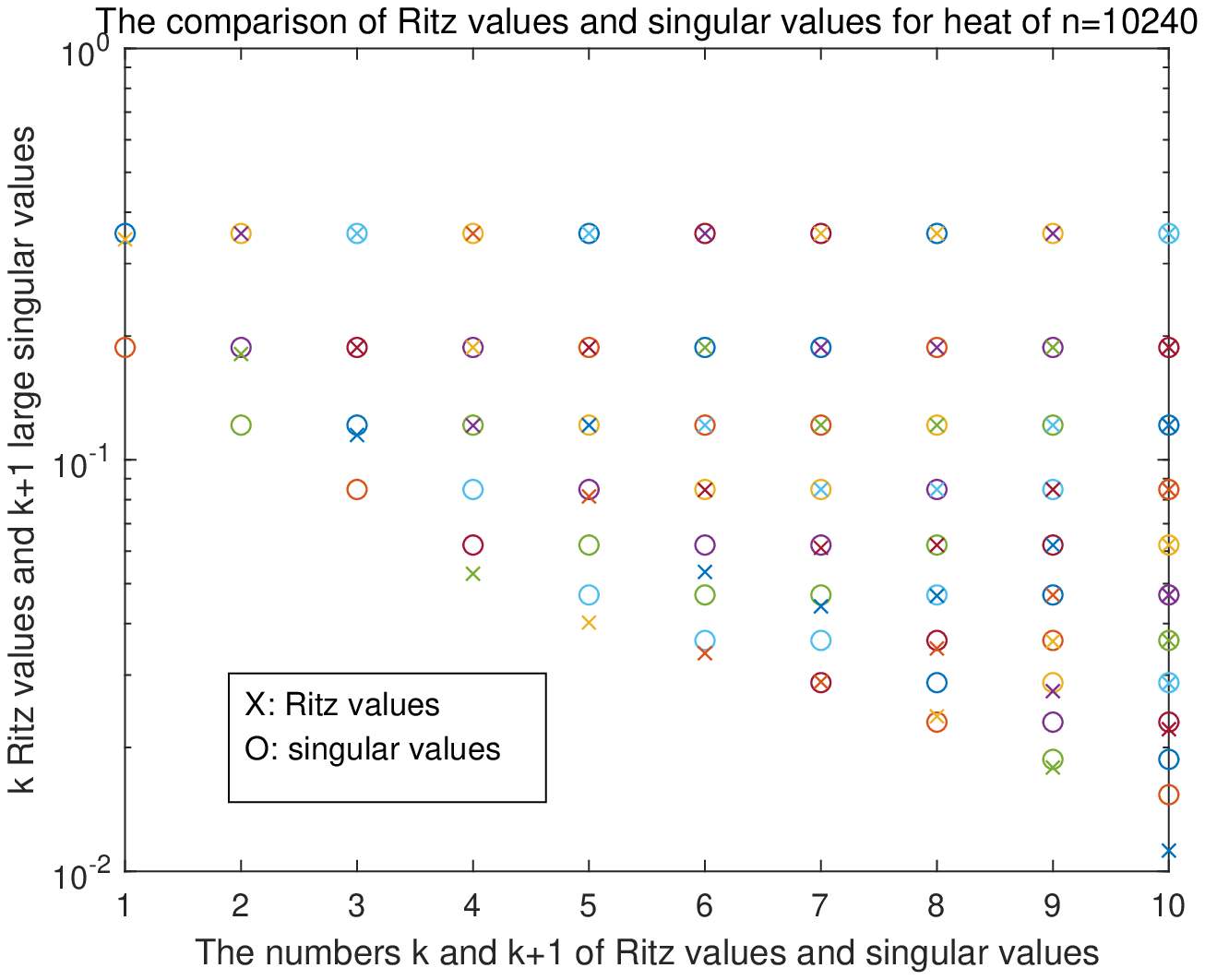}}
\centerline{(b)}
\end{minipage}
\caption{(a) The true  $\|\sin\Theta(\mathcal{V}_k,\mathcal{V}_k^R)\|$
and the required sufficient conditions on them; (b): $k$ Ritz values
and the first $k+1$ large singular values
of {\sf heat}, $k=1,2,\ldots,10$.}
\label{figheat2}
\end{figure}

\begin{figure}[!htp]
\begin{minipage}{0.48\linewidth}
  \centerline{\includegraphics[width=6.0cm,height=3.5cm]{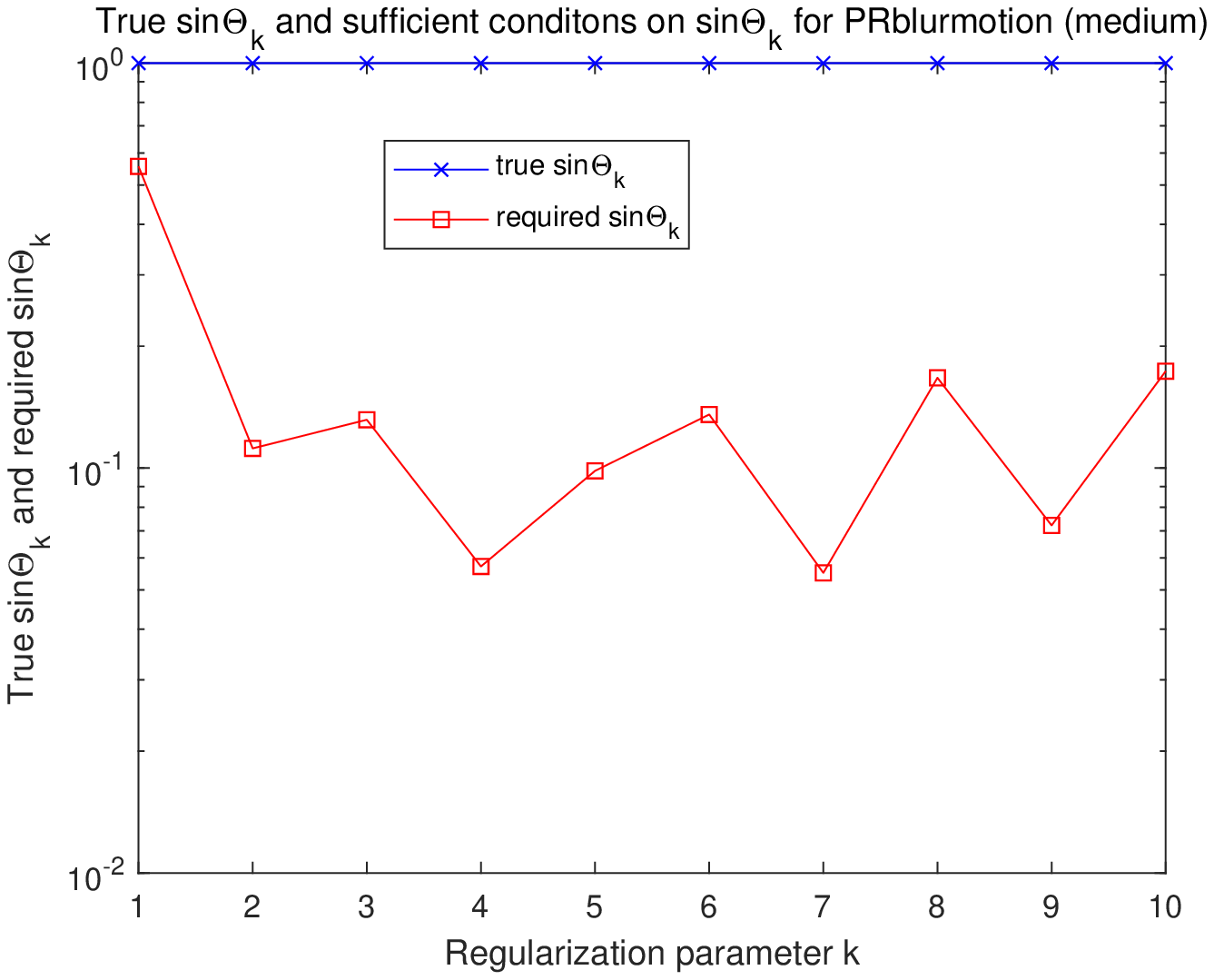}}
\centerline{(a)}
\end{minipage}
\hfill
\begin{minipage}{0.48\linewidth}
  \centerline{\includegraphics[width=6.0cm,height=3.5cm]{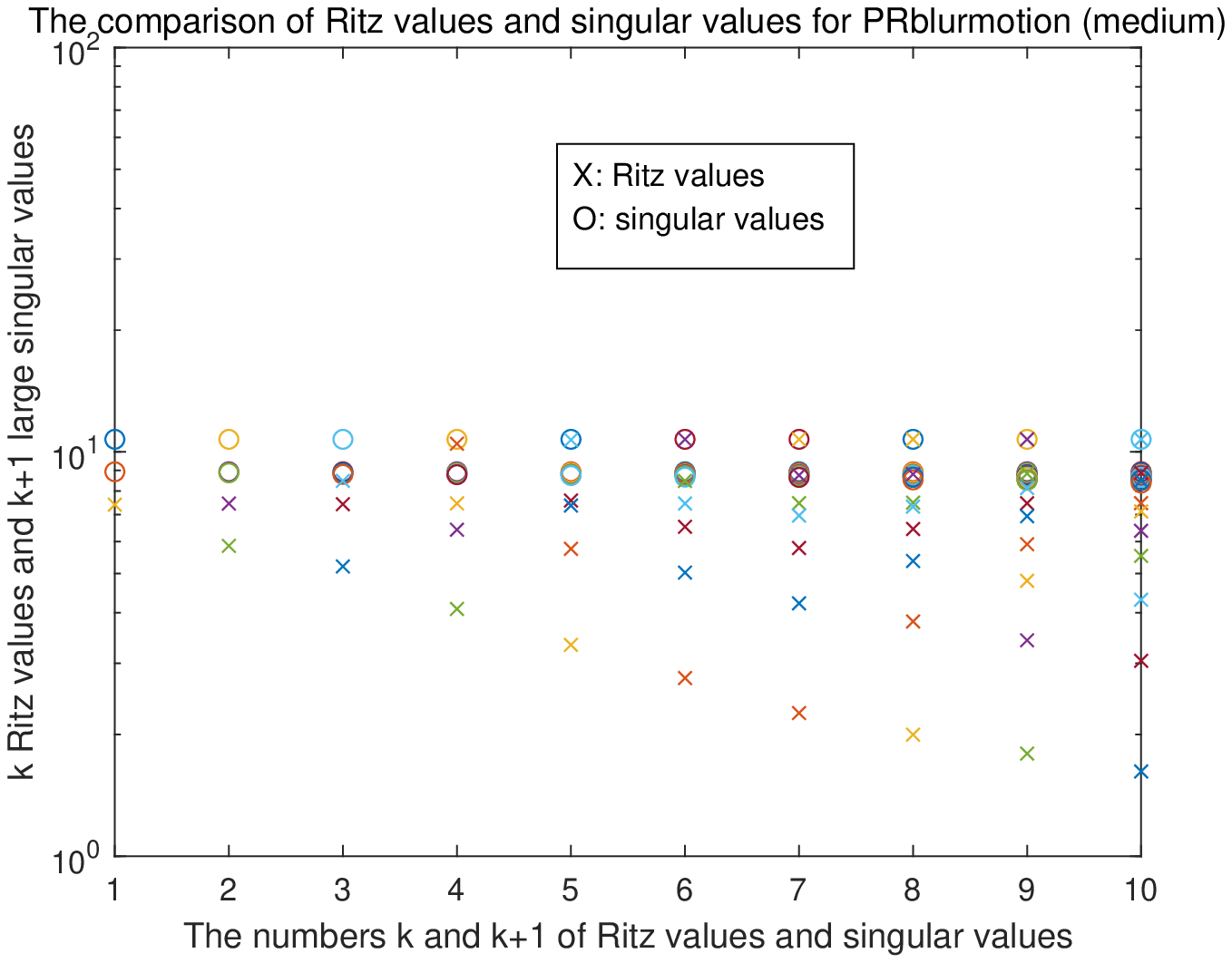}}
\centerline{(b)}
\end{minipage}
\caption{(a): The true $\|\sin\Theta(\mathcal{V}_k,\mathcal{V}_k^R)\|$
and the required sufficient conditions on them; (b): $k$ Ritz values
and the first $k+1$ large singular values
of {\sf PRblurmotion}, $k=1,2,\ldots,10$.}
\label{figPRblur2}
\end{figure}

\begin{figure}[!htp]
\begin{minipage}{0.48\linewidth}
  \centerline{\includegraphics[width=6.0cm,height=3.5cm]{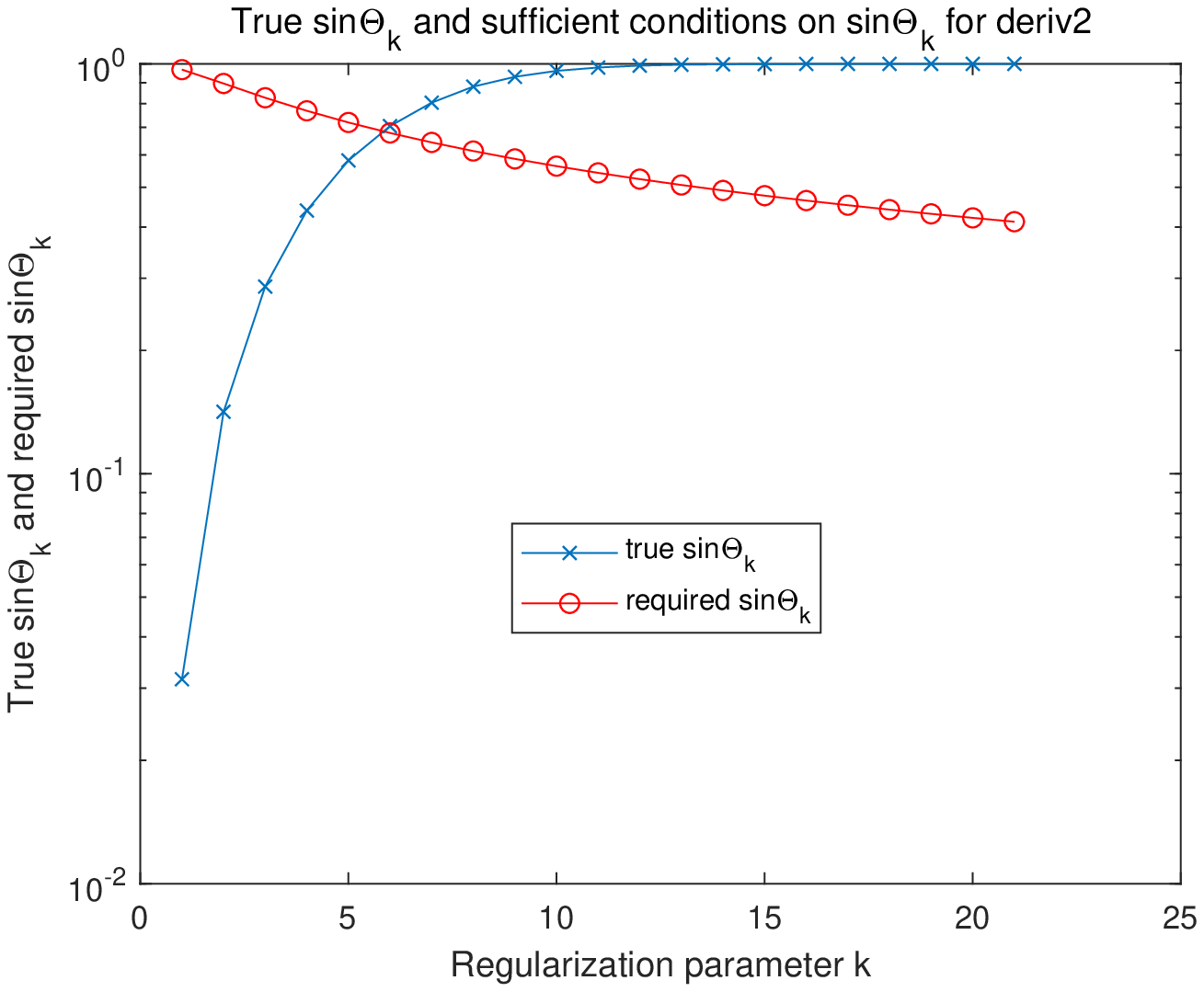}}
\centerline{(a)}
\end{minipage}
\hfill
\begin{minipage}{0.48\linewidth}
  \centerline{\includegraphics[width=6.0cm,height=3.5cm]{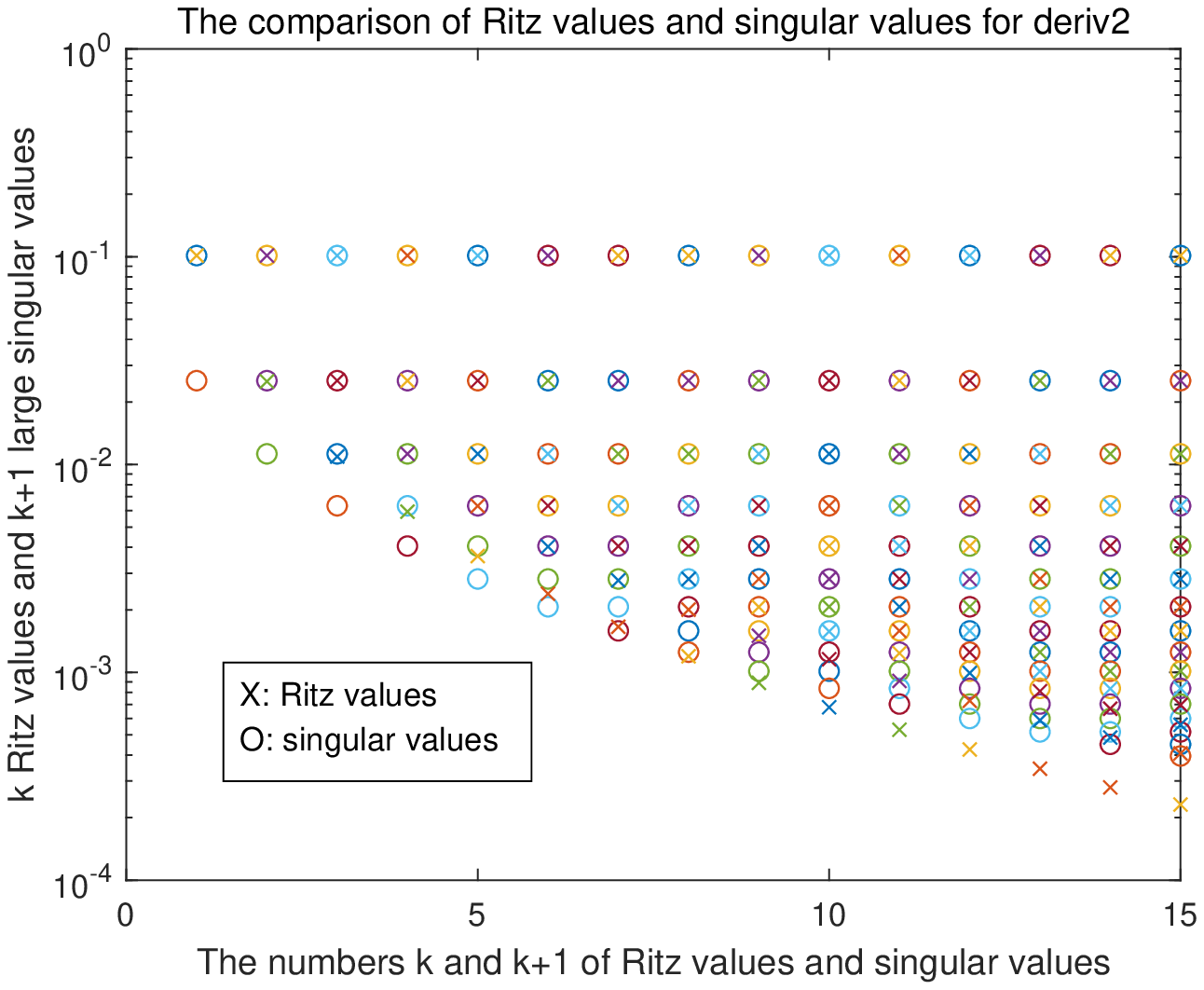}}
\centerline{(b)}
\end{minipage}
\caption{(a): The true $\|\sin\Theta(\mathcal{V}_k,\mathcal{V}_k^R)\|$
and the required sufficient conditions on them; (b): $k$ Ritz values
and the first $k+1$ large singular values
of {\sf deriv2}, $k=1,2,\ldots,15$.}
\label{figderiv2}
\end{figure}

\begin{figure}[!htp]
\begin{minipage}{0.48\linewidth}
  \centerline{\includegraphics[width=6.0cm,height=3.5cm]{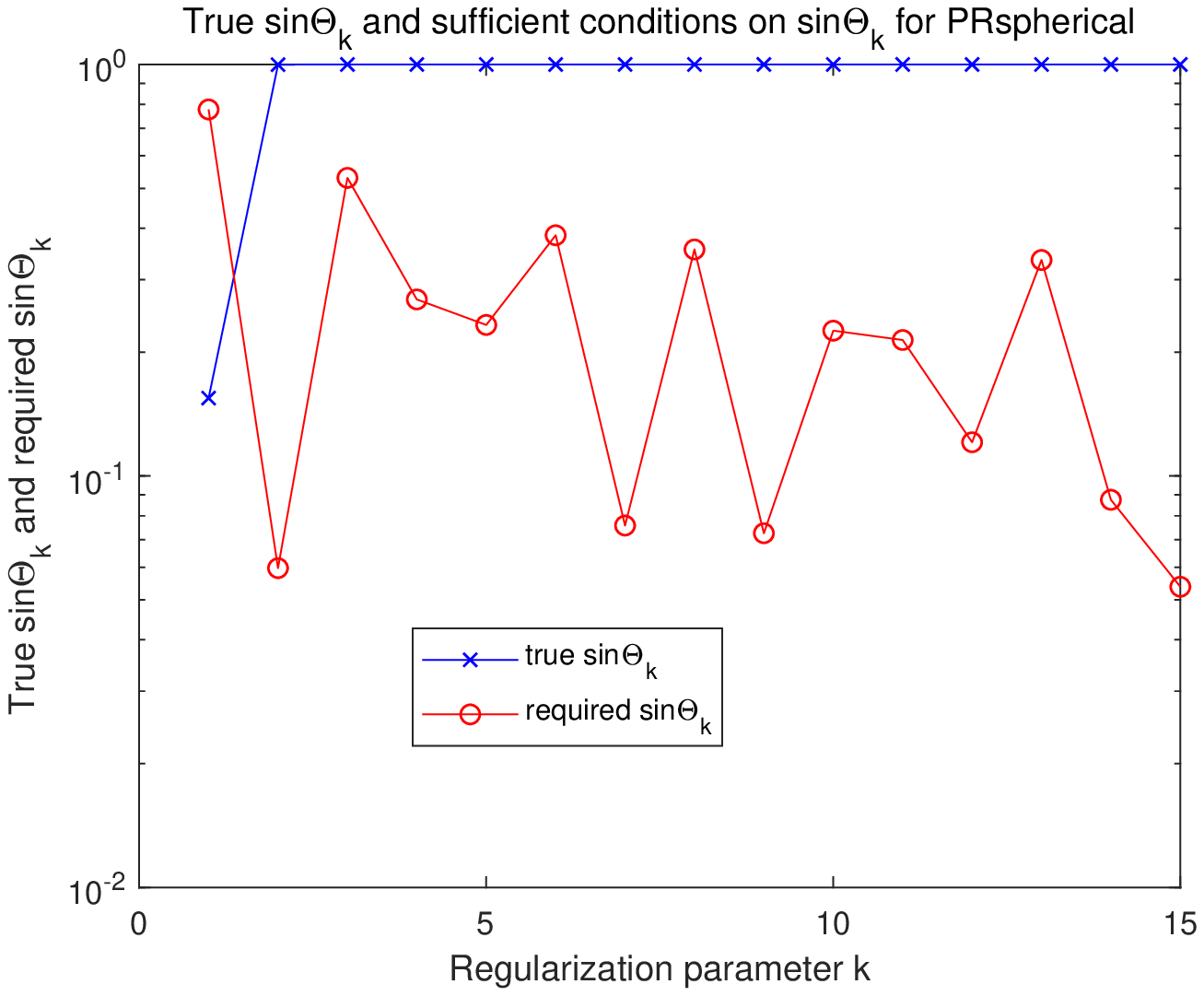}}
\centerline{(a)}
\end{minipage}
\hfill
\begin{minipage}{0.48\linewidth}
  \centerline{\includegraphics[width=6.0cm,height=3.5cm]{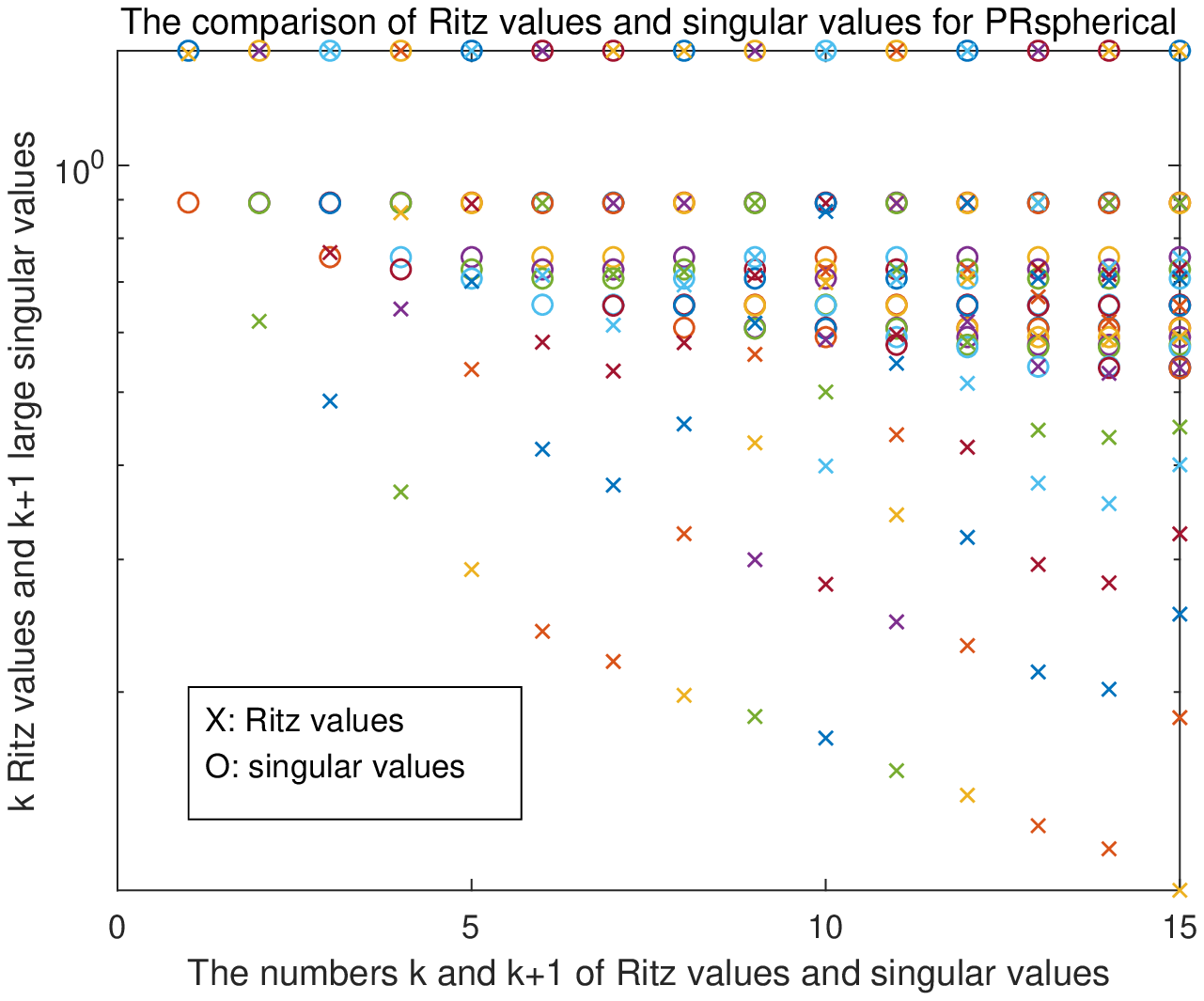}}
\centerline{(b)}
\end{minipage}
\caption{(a): The true $\|\sin\Theta(\mathcal{V}_k,\mathcal{V}_k^R)\|$
and the required sufficient conditions on them; (b): $k$ Ritz values
and the first $k+1$ large singular values
of {\sf PRspherical}, $k=1,2,\ldots,10$.}
\label{figspherical1}
\end{figure}

\begin{figure}[!htp]
\begin{minipage}{0.48\linewidth}
  \centerline{\includegraphics[width=6.0cm,height=3.5cm]{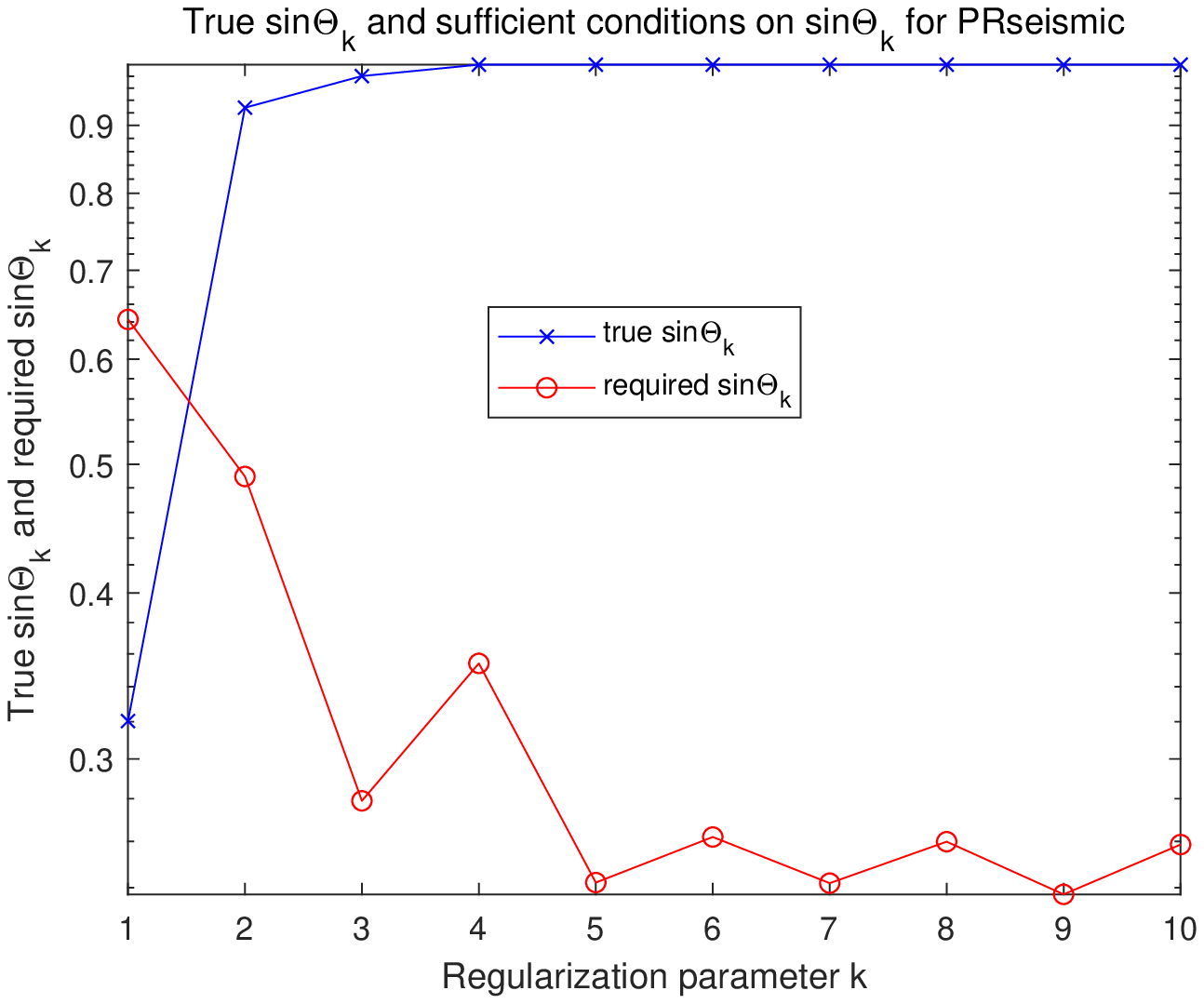}}
\centerline{(a)}
\end{minipage}
\hfill
\begin{minipage}{0.48\linewidth}
  \centerline{\includegraphics[width=6.0cm,height=3.5cm]{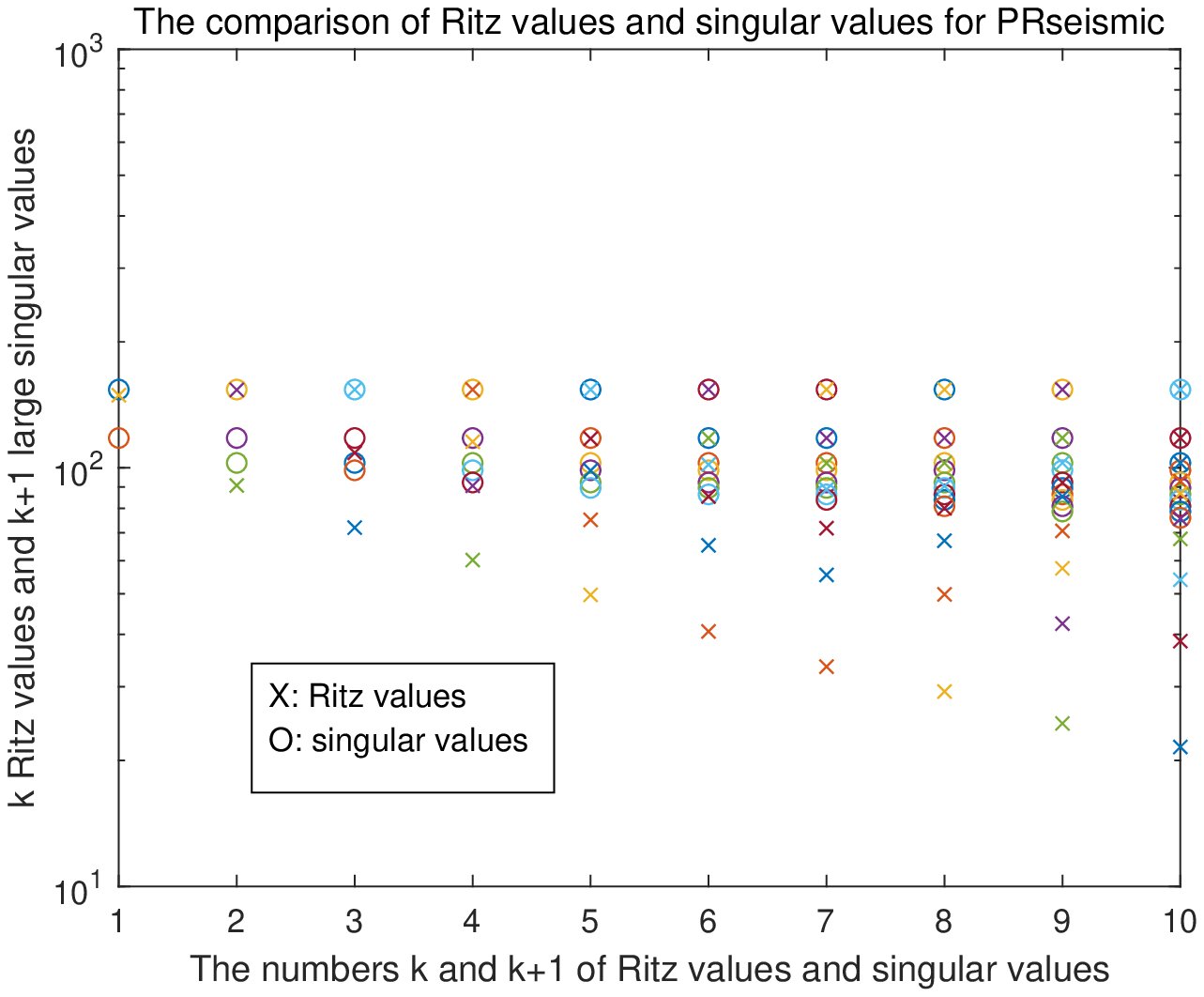}}
\centerline{(b)}
\end{minipage}
\caption{(a): The true $\|\sin\Theta(\mathcal{V}_k,\mathcal{V}_k^R)\|$
and the required sufficient conditions on them; (b): $k$ Ritz values
and the first $k+1$ large singular values
of {\sf PRseismic}, $k=1,2,\ldots,10$.}
\label{figprseismic1}
\end{figure}

For each of the six problems, we first investigate
the true $\|\sin\Theta(\mathcal{V}_k,\mathcal{V}_k^R)\|$ and
the required sufficient conditions
$\|\sin\Theta(\mathcal{V}_k,\mathcal{V}_k^R)\|=\sqrt{1-\varepsilon_k^2}$
that makes \eqref{est1} hold, from which,
\eqref{thetak} and \eqref{approxeq} it is known that $\theta_k^{(k)}>\sigma_{k+1}$
may hold. In the tests, we take $\varepsilon_k=\frac{\sigma_{k+1}}{\sigma_k}$
and compute $\|\sin\Theta(\mathcal{V}_k,\mathcal{V}_k^R)\|
=\sqrt{1-\varepsilon_k^2}$.
We check how the required sufficient conditions are met for each problem and
a given $k$. We depict the true $\|\sin\Theta(\mathcal{V}_k,\mathcal{V}_k^R)\|$
versus the required ones in
Figures~\ref{figshaw2}--\ref{figprseismic1} (a) and
draw the comparison diagrams of $k$ Ritz values $\theta_i^{(k)}$
and first $k+1$ large singular values $\sigma_i$ for each $k$
in Figures~\ref{figshaw2}--\ref{figprseismic1} (b).

Figure~\ref{figshaw2} (a) indicates that for {\sf shaw}
the required sufficient conditions are fulfilled in the first 20 iterations
except for $k=18$.
Figures~\ref{figheat2} (a) and Figures~\ref{figderiv2}--\ref{figprseismic1}
(a) show that for {\sf heat}, {\sf deriv2}, {\sf PRspherical} and
{\sf PRseismic}
the sufficient conditions on $\|\sin\Theta(\mathcal{V}_k,\mathcal{V}_k^R)\|$
are satisfied until $k=3$, $k=5$, $k=2$ and $k=1$, respectively, after which
the true $\|\sin\Theta(\mathcal{V}_k,\mathcal{V}_k^R)\|$ starts to increase and
approaches one quickly; for {\sf PRblurmotion}, it is even worse
that the required $\|\sin\Theta(\mathcal{V}_k,\mathcal{V}_k^R)\|$
are never met for any $k\geq 1$, as shown by Figure~\ref{figPRblur2} (a).
These results justify our theory that (i)
the required $\|\sin\Theta(\mathcal{V}_k,\mathcal{V}_k^R)\|$ are
met more easily for severely ill-posed problems than for moderately
and mildly ill-posed problems; (ii)
the required $\|\sin\Theta(\mathcal{V}_k,\mathcal{V}_k^R)\|$ are fairly
small for moderately and especially mildly ill-posed problems,
but the true $\|\sin\Theta(\mathcal{V}_k,\mathcal{V}_k^R)\|$ approach one as $k$
increases and they tend to one faster for mildly ill-posed problems;
(iii) for moderately and mildly ill-posed problems the true
$\|\sin\Theta(\mathcal{V}_k,\mathcal{V}_k^R)\|$ exhibit monotonically
increasing tendency and approach one with $k$, which confirms our
results.

Next we numerically investigate the behavior of
the smallest Ritz value $\theta_k^{(k)}$
and verify close relationships between it and the required sufficient
condition on $\|\sin\Theta(\mathcal{V}_k,\mathcal{V}_k^R)\|$ that
ensures $\theta_k^{(k)}>\sigma_{k+1}$.
For {\sf shaw}, we see from Figure~\ref{figshaw2} (b) that all the
$\theta_k^{(k)}$ are above $\sigma_{k+1}$ for $k=1,2,\ldots,20$,
including $k=18$ at which the sufficient condition fails to meet. This
indicates that the sufficient condition is not necessary
for ensuring $\theta_k^{(k)}>\sigma_{k+1}$.
For {\sf heat}, Figure~\ref{figheat2} (b) clearly shows
that the $k$ Ritz values $\theta_i^{(k)}$ approximate
the first $k$ large singular values $\sigma_i$ of {\sf heat}, which
includes $\theta_k^{(k)}>\sigma_{k+1}$, for $k=1,2,3$, at which
the required sufficient conditions are satisfied, and
$\theta_k^{(k)}<\sigma_{k+1}$ appears exactly from $k=4$ onwards.
This example illustrates that the required sufficient condition
is also necessary, for if they are not met then
$\theta_k^{(k)}<\sigma_{k+1}$. The numerical results on {\sf deriv2} are similar
to those on {\sf heat}; from Figure~\ref{figderiv2} (b), we see that
$\theta_k^{(k)}>\sigma_{k+1}$ until $k=6$, after which
the required sufficient condition fails to fulfill
and $\theta_k^{(k)}<\sigma_{k+1}$ occurs.

For the mildly ill-posed {\sf PRblurmotion},
notice that the required sufficient conditions on
$\|\sin\Theta(\mathcal{V}_k,\mathcal{V}_k^R)\|$ fail to meet for
all $k\geq 1$. By Theorem~\ref{initial}, the $k$ Ritz values $\theta_i^{(k)}$ may
not approximate the first $k$ singular values $\sigma_i$ in natural order
for all $k\geq 1$. This is indeed the case, as shown clearly by
Figure~\ref{figPRblur2} (b), which shows that  all the
$\theta_k^{(k)}$ are below $\sigma_{k+1}$.
Regarding the mildly ill-posed {\sf PRspherical} and
{\sf PRseismic}, the sufficient conditions are satisfied only for $k=1$,
as is seen from Figures~\ref{figspherical1}--\ref{figprseismic1} (a). The $k$
Ritz values $\theta_i^{(k)}$ interlace the first $k+1$ large singular
values $\sigma_i$ in natural order only for $k=1$, and
afterwards $\theta_k^{(k)}<\sigma_{k+1}$,
as indicated clearly by Figures~\ref{figspherical1}--\ref{figprseismic1} (b).
For these two problems, at iteration
$k=1$ the Ritz value $\theta_1^{(1)}$ lies between $\sigma_1$ and
$\sigma_2$ and is closer to $\sigma_1$.
Again, these results demonstrate that our sufficient conditions are tight.
Moreover, compared with the previous problems, we find that, generally,
the more slowly the singular values $\sigma_i$ decay, the harder the sufficient
condition is to fulfill, and the sooner $\theta_k^{(k)}<\sigma_{k+1}$
occurs.

\begin{figure}[!htp]
\begin{minipage}{0.48\linewidth}
  \centerline{\includegraphics[width=6.0cm,height=3.5cm]{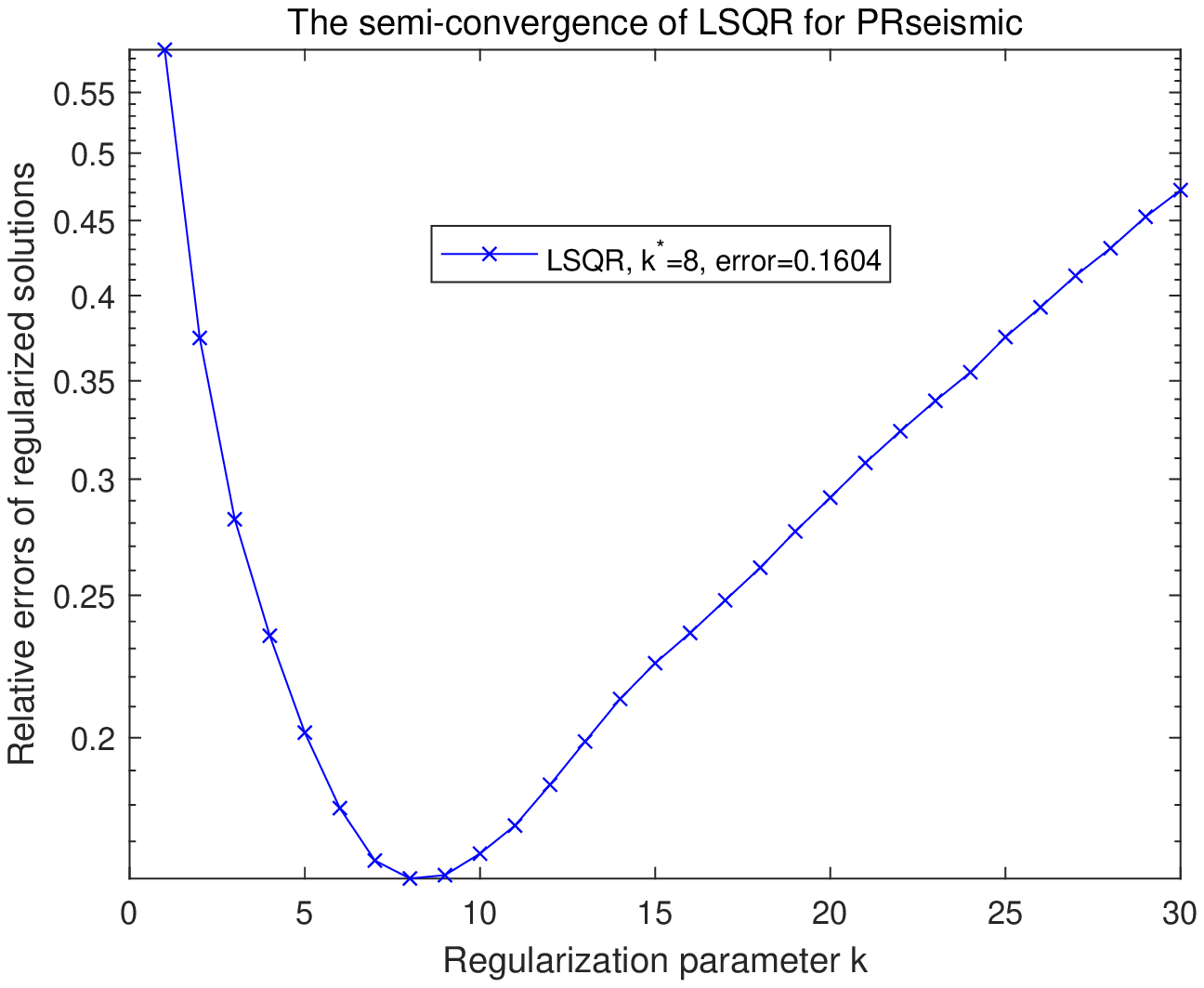}}
\centerline{(a)}
\end{minipage}
\hfill
\begin{minipage}{0.48\linewidth}
  \centerline{\includegraphics[width=6.0cm,height=3.5cm]{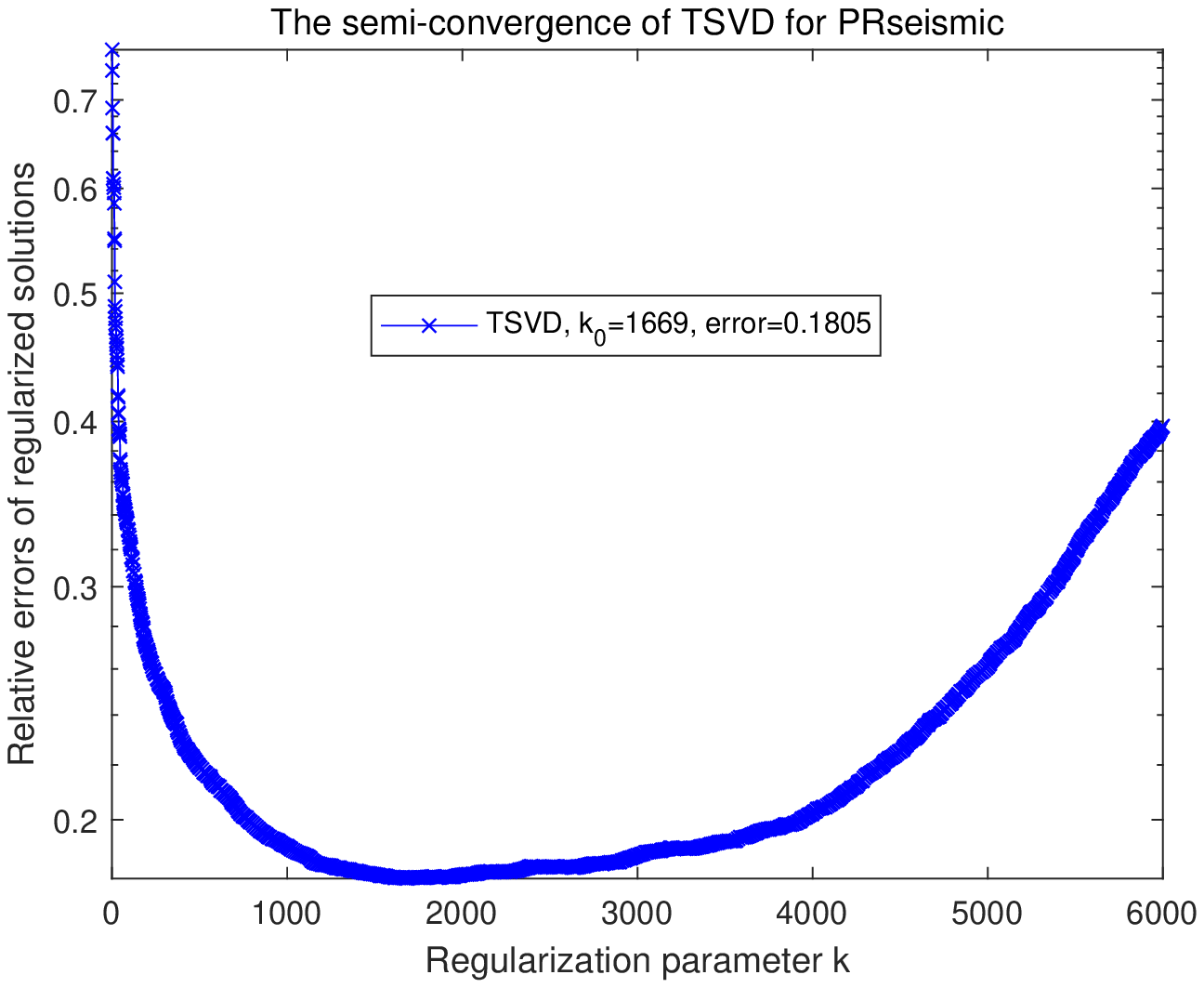}}
\centerline{(b)}
\end{minipage}
\caption{(a)-(b): The semi-convergence processes of LSQR and TSVD for
{\sf PRseismic}.}
\label{figprseismic2}
\end{figure}

\begin{figure}[!htp]
\begin{minipage}{0.48\linewidth}
  \centerline{\includegraphics[width=6.0cm,height=3.5cm]{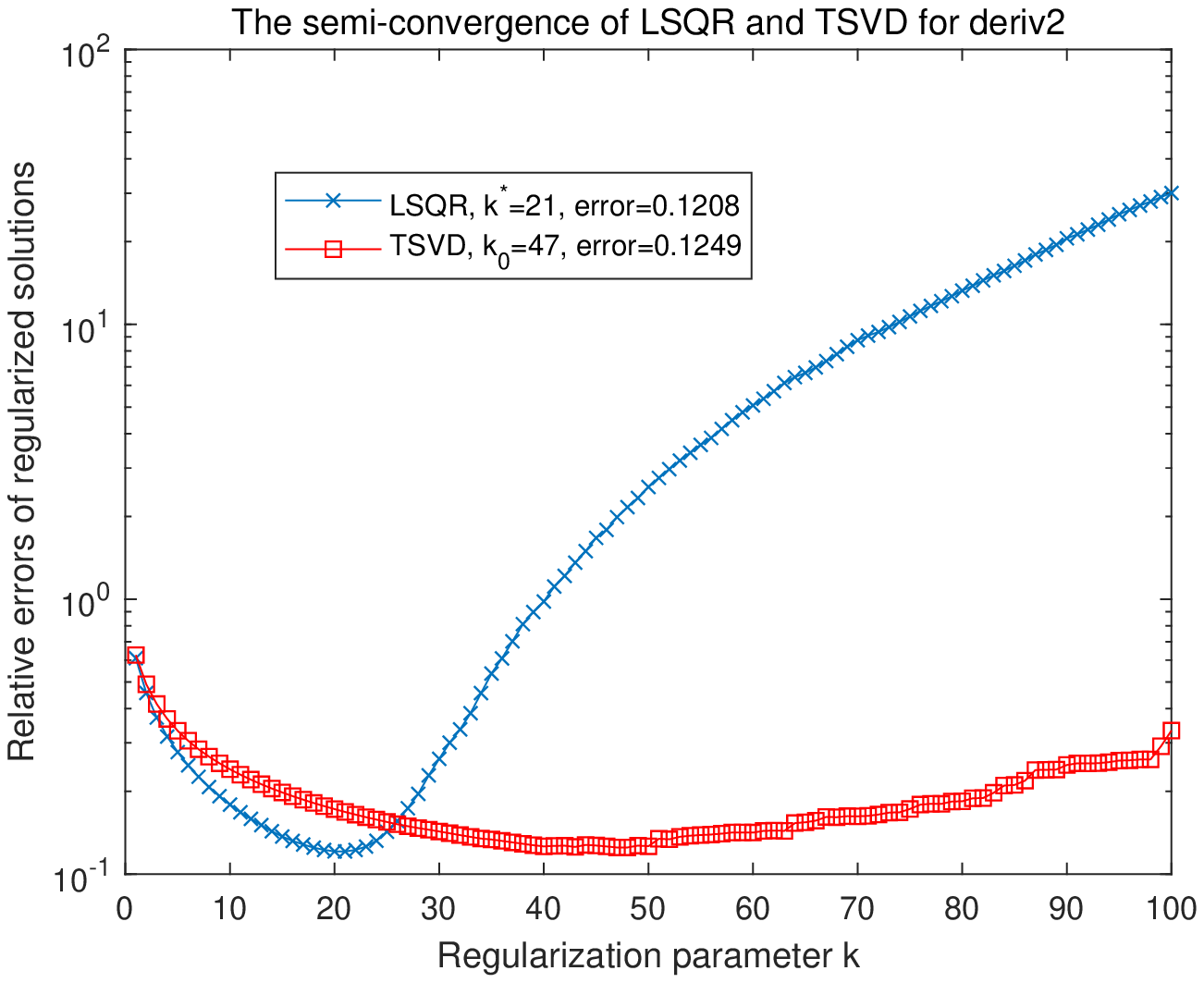}}
\centerline{(a)}
\end{minipage}
\hfill
\begin{minipage}{0.48\linewidth}
  \centerline{\includegraphics[width=6.0cm,height=3.5cm]{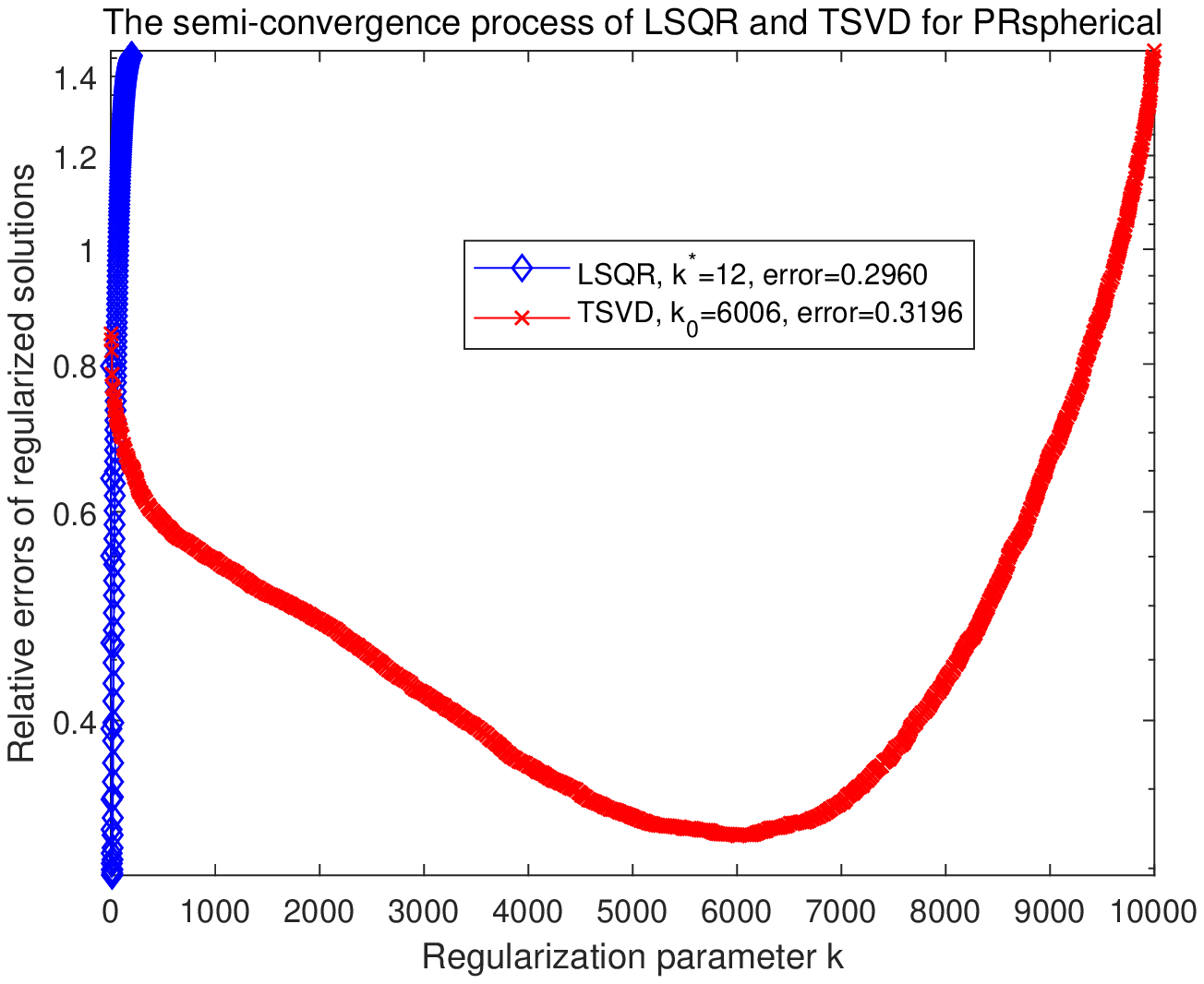}}
\centerline{(b)}
\end{minipage}
\caption{(a): the semi-convergence process of LSQR and TSVD for {\sf deriv2};
(b): the semi-convergence process
of LSQR and TSVD for {\sf PRspherical2}.}
\label{figspherical2}
\end{figure}

In Figure~\ref{figprseismic2}, we
depict the semi-convergence processes of LSQR and the TSVD method for
{\sf PRseismic}, where it is seen that $k^*=8$ and $k_0=1669$.
In Figure~\ref{figspherical2},
we depict the semi-convergence processes of LSQR and the TSVD method for
{\sf deriv2} and {\sf PRspherical}, where the semi-convergence point $k^*=21$
of LSQR and the transition point $k_0=47$ of the TSVD method
for {\sf deriv2}, and $k^*=12$ and $k_0=6006$ for {\sf PRspherical}.
For these three problems, we find
that $k^*\ll k_0$, especially for {\sf PRseismic} and {\sf PRspherical}.
However, it is clear from the figures that
the best regularized solutions by LSQR are at least
as accurate as those computed by the TSVD method. Again,
they results illustrate
that the approximations of $\theta_i^{(k)}$ to the large $\sigma_i$
in natural order until the occurrence of semi-convergence of LSQR are
{\em not necessary} conditions for the full regularization of LSQR.

\section{Conclusions}\label{concl}

For a general large-scale linear discrete ill-posed problem
\eqref{eq1}, the Krylov iterative solvers
LSQR and CGLS are most popularly used, and they, together with CGME and LSMR,
are deterministic 2-norm filtering iterative regularization methods.
These methods have general regularizing effects and exhibit
semi-convergence. For each of them, if the regularized solution at
semi-convergence are as accurate as the best one obtained by the TSVD method,
which has been known to find a 2-norm filtering best possible solution,
the method has the full regularization. In this case, for a given
problem, once the semi-convergence is practically recognized, we
have computed a best possible regularized solution and simply
stop the method. The determination of semi-convergence can,
in principle, be determined by a suitable parameter-choice
method, such as the L-curve criterion and the discrepancy principle.

In the simple singular value case,
as the first and fundamental step towards to understanding the regularization
of LSQR, CGME and LSMR, we have established
the $\sin\Theta$ theorem for the 2-norm
distance $\|\sin\Theta(\mathcal{V}_k,\mathcal{V}_k^R)\|$
between the underlying $k$ dimensional Krylov subspace $\mathcal{V}_k^R$ and
the $k$ dimensional dominant right singular subspace $\mathcal{V}_k$, and
derived accurate estimates on the distances for the three kinds of
ill-posed problems under the simplifying assumptions on the actual decay
of the singular values of $A$. Then we have manifested some
intrinsic relationships between the smallest Ritz values $\theta_k^{(k)}$
and $\|\sin\Theta(\mathcal{V}_k,\mathcal{V}_k^R)\|$. The results will
provide absolutely necessary background and ingredients for
studying the problems mentioned in the beginning of Section~\ref{sine}.

We have reported illuminating numerical examples to show that our estimates
are sharp and realistic, and have justified that our sufficient conditions on
$\theta_k^{(k)}>\sigma_{k+1}$ and $\theta_k^{(k)}<\sigma_{k+1}$ are tight.
Also, we have numerically confirmed some
important properties on the factors $|L_j^{(k)}(0)|$, $j=1,2,\ldots,k$.

Very surprisingly, we have found that for all the test problems, independent
of the degree of ill-posedness of \eqref{eq1}, the LSQR best
regularized solutions $x_{k^*}^{lsqr}$ are as accurate as
the TSVD best solutions $x_{k_0}^{tsvd}$. This indicates that LSQR has the full
regularization for the test problems, even though the Ritz values do not approximate
the large singular values of $A$ in natural before some $k\leq k^*$.
We have made numerical experiments on each of the test problems for different
noise levels $\varepsilon$ and have had the same findings.
Furthermore, we have observed that LSMR has the full
regularization too for the test problems, but it is
not the case for CGME, whose best regularized solutions
are considerably less accurate than and can hardly be as accurate as
those by LSQR, even for severely ill-posed problems.
Therefore, the full or partial regularization
of these Krylov solvers is much more complicated that one may have expected,
and deserves high attention and in-depth study.



\end{document}